\documentclass[12pt]{article}
\newcommand{\bfns}{\begin{footnotesize}}

\newcommand{\efns}{\end{footnotesize}}

\newtheorem{lemma}{Lemma}
\newtheorem{definition}{Definition}

\newtheorem{theorem}{Theorem}
\usepackage{latexsym}
\usepackage{amssymb}
\usepackage{multirow}
\usepackage{verbatim}
\usepackage{graphicx} 
\usepackage{color}
\usepackage{url}

\newcommand{\FF}{{\cal F}}

\begin{title}{Graphs of girth at least $5$ with orders \\  between $20$ and $32$}\end{title}
\author{Alice Miller$^{1}$ and  Michael Codish$^{2}$\\
\small{$^{1}$ School of Computing Science, University of Glasgow, UK,} \\\small{$^{2}$ Department of Computer Science, Ben-Gurion University of the Negev, Israel.}}

\begin{document}
\maketitle
\begin{abstract} We prove properties of extremal graphs of girth $5$ and  order $20\leq v \leq 32$. In each case we identify the possible minimum and maximum degrees, and in some cases prove the existence of (non-trivial) embedded stars. These proofs allow for tractable search for and identification of all non isomorphic cases.
\end{abstract}

\section{Introduction} In this paper we consider graphs of girth at least $5$, i.e. graphs which have no 3-cycles or 4-cycles. In \cite{Garnick93}, the maximum number of edges in graphs of girth at least $5$ with $v$ vertices, $f_4(v)$, were established for $v\leq 30$, and $v=50$, and the number of unique (up to isomorphism) graphs with $f_4(v)$ edges,  $F_4(v)$,  was caculated for $v\leq 21$ and $v=50$. Note that, once $f_4(v)$ has been established, it is straightforward to generate $F_4(v)$ for $n<20$ using nauty's geng tool \cite{nauty}.In \cite{comiprst1} and \cite{comiprst2} we introduce two symmetry breaking constraints for combinatorial search and illustrated their effectiveness by applying it to the problem of finding graphs with girth at least $5$ (and other combinatorial problems). We were able to reproduce the results of \cite{Garnick93} and find values  $F_4(v)$ that were previously unpublished. Our symmetry breaking constraints alone were not sufficient to crack some of the harder cases. We were however able to reduce our search space by applying new theoretical results for cases $v>22$. In many cases we were able to identify non-trivial {\it embedded stars} that must be present in a graph with the maximum number of edges and girth at least $5$, and fixing this embedded star made our search tractable. An outline of the proof of the existence of such stars is given in \cite{comiprst2}, but for space reasons did not include full proofs. We give the full proofs in this paper. All of the graphs for $F_{4}(v)$, where $20\leq v \leq 32$ are presented as incidence arrays in our appendices and a full list in {\it graph6} notation (as implemented in nauty \cite{nauty}) can be obtained by contacting the authors directly.

\section{Preliminary Definitions and Results}
In this section we give some preliminary definitions and results that will be useful in this paper. 

The girth of a graph $\Gamma=(V,E)$ is the size of the smallest cycle contained in
it. 
Let $\FF_k(v)$ denote the set of graphs with $v$ vertices and girth at
least $k+1$. Let $f_k(v)$ denote the maximum number of edges in a
graph in $\FF_k(v)$. A graph in $\FF_k(v)$ with $f_k(v)$ edges is
called \emph{extremal}. The number of non-isomorphic extremal graphs
in $\FF_k(v)$ is denoted $F_k(v)$.
Extremal graph problems involve discovering values of $f_k(v)$ and
$F_k(v)$ and finding witnesses.
In~\cite{AbajoD10} the authors attribute the discovery of values
$f_4(v)$ for $v\leq 24$ to \cite{Garnick93} and for $25\leq v\leq 30$
to \cite{Garnick92}.  In \cite{Garnick93} the authors report values of
$F_4(v)$ for $v\leq 21$. In \cite{Garnick93} and \cite{WangDM01}
algorithms are applied to compute lower bounds on $f_4(v)$ for $31\leq
v\leq 200$. Some of these lower bounds are improved in
\cite{AbajoBD10} and improved upper bounds for $33\leq
v\leq 42$ are proved in \cite{bong2017}.
Values of $f_4(v)$ for $v\leq 30$, and of $F_4(v)$ for $v\leq 21$ are
available as sequences \texttt{A006856} and \texttt{A159847} of the
On-Line Encyclopedia of Integer Sequences~\cite{oeis}. We have extended sequence  \texttt{A159847} to include values for $22\leq v \leq 32$, following the results presented in this paper.

In the remainder of this paper, as we are only interested in graphs with girth at least $5$, we abbreviate  $f_4(v)$ and $F_4(v)$ to $f(v)$ and $F(v)$ respectively. 

\begin{definition} A graph $\Gamma=(V,E)$ is said to be {\it extremal} if it has girth at least $5$ and any graph on $V$ vertices with more than $e=|E|$ edges has girth less than $5$. For any $v$, $f(v)$ and $F(v)$ denote the number of edges in an extremal graph on $v$ vertices and the number of unique extremal graphs on $v$ vertices respectively. For any extremal graph, $\delta$ and $\Delta$ denote the minimum and maximum degree respectively.
\end{definition}

The following lemma, and it's proof, is from \cite{Garnick93}.
\begin{lemma}\label{lem:garnickLemma}
If $\Gamma$ is an extremal graph of order $v$, with $e=f(v)$ edges, whose minimum and maximum degrees are $\delta$ and $\Delta$ respectively, then 
$$e-f(v-1)\leq \delta \leq \sqrt(v-1)\;{\rm and}\; 
\lceil 2e/v\rceil\leq  \Delta\leq  (v-1)/\delta$$
\end{lemma}
Values of $f(v)$ for $22\leq v \leq 30$ can be found in \cite{Garnick93}, and those for $v=31$ and $v=32$ are given in \cite{comiprst1}. By applying Lemma \ref{lem:garnickLemma}, and observing (by counting edges) that $\delta=\Delta$ is only possible when $v\Delta = 2f(v)$ we obtain possible values of $(\delta,\Delta)$, for $22\leq v \leq 33$. These are given in Table \ref{table:deltaTable}

\begin{table}
\begin{tabular}{l|l|l}
$v$ & $f(v)$ & $(\delta,\Delta)$\\
\hline
21&44&$(3,5)$, $(3,6)$, $(4,5)$\\
22&47&$(3,5)$, $(3,6)$, $(3,7)$, $(4,5)$\\
23&50&$(3,5)$, $(3,6)$, $(3,7)$, $(4,5)$\\
24&54&$(4,5)$\\
25&57&$(3,5)$, $(3,6)$, $(3,7)$, $(3,8)$, $(4,5)$, $(4,6)$\\
26&61&$(4,5)$, $(4,6)$\\
27&65&$(4,5)$, $(4,6)$\\
28&68&$(3,5)$, $(3,6)$, $(3,7)$, $(3,8)$, $(3,9)$, $(4,5)$, $(4,6)$\\
29&72&$(4,5)$, $(4,6)$, $(4,7)$\\
30&76& $(4,6)$, $(4,7)$\\
31&80& $(4,6)$, $(4,7)$, $(5,6)$\\
32&85&$(5,6)$\\
33&87&$(2,6)$, $(2,7)$, $(2,8)$, $(2,9)$, $(2,10)$, $(2,11)$, $(2,12)$,\\
&&        $(2,13)$, $(2,14)$, $(2,15)$, $(2,16)$, $(3,6)$, $(3,7)$, $(3,8)$,\\
&&        $(3,9)$, $(3,10)$, $(4,6)$, $(4,7)$, $(4,8)$, $(5,6)$\\ 
\hline
\end{tabular}
\caption{Possible values of $(\delta,\Delta)$ for $22\leq v \leq 32$\label{table:deltaTable}}
\end{table}

\begin{definition}\label{def:embedded}
 If $\Gamma$ is a graph of girth at least $5$ we say that $\Gamma$
has a embedded $S_{D,[d_{0}-1,d_{1}-1, \ldots, d_{D-1}-1]}$  star if $\Gamma$ has a vertex of degree $D$ whose children have degrees at least $d_{0},d_{1},\ldots,d_{D-1}$ respectively. 
\end{definition}
Note that the children and grandchildren of a vertex in a a graph of girth at least $5$ are distinct, and there are no edges between the children. If $\Gamma$ is such a graph and has minimum and maximum degrees $\delta$ and $\Delta$, then any vertex $x$ of degree $\Delta$ is the centre of a {\it trivial} embedded 
$S_{\Delta,[\delta-1,\delta-1, \ldots, \delta-1]}$ star.  In order to reduce the number of isomorphic graphs produced, any combinatorial search for extremal graphs will assume a fixed position of such a star (for example, with central node as vertex $0$). The purpose of this paper is to show that in some cases an extremal graph must contain larger (non-trivial) embedded stars, in some cases containing all of the vertices in the graph. This allows us to fix such a star and  reduce both our search space and the number of isomorphic graphs generated.

\begin{definition}\label{def:sinkNode}
A node $x$  of an extremal graph is said to be a {\it sink node} if it has maximum degree and  is the centre of an embedded star that contains  all of the vertices (i.e. $x$ is at distance at most $2$ from all other vertices). 
\end{definition}

From a result in \cite{Garnick93}, no two vertices in an extremal graph of girth at least $5$ are at distance greater than $3$ from each other. 

\begin{definition}\label{def:setsDifferent}
For a graph $\Gamma=(V,E)$, for any $m>=2$, $S_{\Gamma,m}$ is the set of sets of $m$ vertices that are at a distance of $3$ from each other. 
\end{definition}

\begin{lemma}\label{lem:sinknode} 
If $\Gamma$ is a graph of girth at least $5$ on $v$ vertices and contains a sink node $x$, then (i)  $\Gamma$ has an embedded   $S_{\Delta,[d_{0}-1,d_{1}-1,\ldots, d_{\Delta-1}-1]}$ star, where $\Delta+1+\Sigma_{i=0}^{\Delta-1}(d_{i}-1) = v$. (ii) for $m\geq 2$, no element of  $S_{\Gamma,m}$  contains $x$ or more than one child of $x$.
\end{lemma}

\vspace{.25cm}
\noindent
{\bf Proof} Follows immediately as $x$ is at distance at most $2$ from all other vertices and the children of $x$ are at  distance $2$ from each other. 

\begin{lemma}\label{lem:inductiveExtremal} If $\Gamma=(V,E)$ is an extremal graph and $|V|=v$ with  $\delta(\Gamma)=f(v)-f(v-1)$, then $\Gamma=(V,E)$ is constructed from an extremal graph $\Gamma^{\prime}$ of order $v-1$ by adding a new vertex $x$ of degree $\delta$ to $\Gamma^{\prime}$ and $\delta$ edges from $x$ to a set $S_{\Gamma^{\prime},\delta}$.
\end{lemma}

\begin{theorem} For each $20\leq v \leq 32$, an extremal graph $\Gamma$ must have $(\delta,\Delta)$ taking one of the pairs of values indicated in Table \ref{embeddedStarTable}. In each case, we show the number of edges ($f(v)$), the  largest star known to be embedded in $\Gamma$, the number of distinct graphs $(F(v))$, and the method of proof. In all cases Lemma \ref{lem:inductiveExtremal} is employed. The ``method'' column indicates whether the graphs were obtained from \cite{Garnick93} (G) or, if not, whether search is employed (based on the proven existence of embedded stars) (S), and whether some (or all) graphs are constructed by hand (H). Note that $S,H$ denotes that some of the graphs were found using search, and some by hand.
\end{theorem}

\begin{table}[t]
\label{embeddedStarTable}
\begin{center}
\setlength{\tabcolsep}{3pt}
\begin{tabular}{|c|l|l|l|l|l|}
\hline
v&$f(v)$&$(\delta,\Delta)$&Embedded star&$F(v)$&method\\
\hline
20&$41$&$(3,5)$&$S_{5,[3,3,3,3,2]}$&1&G\\
\hline
21&$44$&$(3,5)$&$S_{5,[3,3,3,3,3]}$&3&H\\
 & &$(4,5)$&$S_{5,[3,3,3,3,3]}$&&\\
\hline
22&$47$&$(3,5)$&$S_{5,[4,3,3,3,3]}$&3&S\\
 & &$(4,5)$&$S_{5,[4,3,3,3,3]}$&&\\
\hline
23&$50$&$(3,5)$&$S_{5,[4,4,3,3,3]}$&7&S\\
&  &$(4,5)$&$S_{5,[4,4,3,3,3]}$ or $S_{5,[4,3,3,3,3]}$&&\\
\hline
24&$54$&$(4,5)$&$S_{5,[4,4,4,3,3]}$&1&S\\
\hline
25&$57$&$(3,5)$&$S_{5,[4,4,4,4,3]}$ or $S_{5,[4,4,4,3,3]}$&6&S\\
&  &    $(4,5)$& $S_{5,[4,4,4,4,3]}$ or $S_{5,[4,4,4,3,3]}$&&\\
&  &    $(4,6)$&$S_{6,[3,3,3,3,3,3]}$&&\\
\hline
26& $61$&$(4,5)$&$S_{5,[4,4,4,4,4]}$&2&H\\
\hline
27&$65$& $(4,5)$&$S_{5,[4,4,4,4,4]}$&1&H\\
\hline
28&$68$& $(3,6)$ & $S_{6,[4,4,4,4,3,2]}$&4&S,H\\ 
&& $(4,5)$ &  $S_{5,[4,4,4,4,4]}$&&\\
&& $(4,6)$ &  $S_{6,[4,4,4,3,3,3]}$&&\\
\hline
29&$72$&$(4,6)$ & $S_{6,[4,4,4,4,3,3]}$&1&H\\
\hline
30&$76$&$(4,6)$ & $S_{6,[4,4,4,4,4,3]}$ and $S_{6,[5,4,4,4,3,3]}$&1&H\\
\hline
31&$80$&$(4,6)$ & $S_{6,[4,4,4,4,4,4]}$ and $S_{6,[5,4,4,4,4,3]}$&2&H,S\\
&&$(5,6)$ & $S_{6,[4,4,4,4,4,4]}$&&\\
\hline
32&$85$&$(5,6)$&$S_{6,[5,4,4,4,4,4]}$&1&H\\
\hline

\end{tabular}
\end{center}
\caption{Embedded stars\label{table:thetable} for 
           $20 \leq v \leq 25$, and $v=32$}
\end{table}

\vspace{.25cm}
\noindent
{\bf Proof} Follows from Theorems \ref{theorem:g20}, \ref{theorem:g21}, \ref{theorem:g22}, \ref{theorem:g23}, \ref{theorem:g24}, \ref{theorem:g25}, \ref{theorem:g26}, \ref{theorem:g27}, \ref{theorem:g28}, \ref{theorem:g29}, \ref{theorem:g30}, \ref{theorem:g31} and \ref{theorem:g32}.

\begin{definition}\label{defn:linearspace}
A linear space $\Lambda$ on a set of $n$ points $V$ is a collection $B=\{B_{1},\ldots,B_{b}\}$ of subsets of $V$ called blocks, such that every block has at least two points and each pair of points is in exactly one block. A prelinear space $\Lambda$ is a set of blocks on $V$ in which every block has at least two points and each pair of points is in at most one block.
\end{definition}

\vspace{.25cm}
\noindent
Any prelinear space can be extended to a linear space, by adding suitable blocks of size $2$, so we tend to use the term linear space and prelinear space interchangeably. A linear space with no blocks of size $2$ is called a {\it proper} linear space.

\begin{lemma}\label{lem:blocks} Let $\Gamma=(V,E)$ be a graph of girth at least $5$ and $$V_{1} = \{v_{0},v_{2},\ldots, v_{m-1}\}$$ a subset of $V$. If, for all $v_{i}\in V_{1}$, $b_{i}$ is the set of elements of $V\setminus V_{1}$ that are adjacent to $v_{i}$, then $B=\{b_{i} : 0\leq i \leq m-1\;{\rm and}\; |b_{i}|\geq 2\}$ is a (pre) linear space on $V\setminus V_{1}$ and no set $b_{i}$ contains two vertices that are adjacent in $\Gamma$.  
\end{lemma}

\vspace{.25cm}
\noindent
{\bf Proof} No pair of elements in $V\setminus V_{1}$ can appear in more than one set $b_{i}$ as $\Gamma$ contains no  $4$-cycles, so $B$ is a (pre) linear space. Similarly, no set $b_{i}$ contains two vertices that are adjacent in $\Gamma$, as $\Gamma$ contains no $3$-cycles. 

\vspace{.25cm}
\noindent
Note that in the context of linear spaces, we generally refer to the elements of the blocks as {\it points}, whereas in the context of graphs the elements of edges are vertices. When applying Lemma \ref{lem:blocks} this distinction is less clear.

\begin{definition}\label{definition:packing} A $(v,k,1)$-packing design is a linear space on a set $V$ where the blocks all have size $k$. For any $v,k$ function $P(k,v)$ is called the maximum packing number of $k$ and $v$, and is the maximum number of blocks in such a packing design.
\end{definition}
The following result follows from a simple counting argument (see \cite{johnson1}):
\begin{lemma} For any $k\leq v$, 
$$P(k,v)\leq \lfloor v\lfloor(v-1)/(k-1)\rfloor/k\rfloor$$
\end{lemma}
The upper bound is not achievable in some cases. Absolute values for $P(k,v)$ for $3\leq k\leq 6$ and $1\leq v\leq 28$ are given in Table \ref{table:packingTable}. These are obtained from \cite{spence}, \cite{brslsm}, and theoretical results in \cite{brouwer1}.

\begin{table}
\begin{tabular}{|l|l|l|l|l|}
\hline
$v$ & $P(3,v)$ &$P(4,v)$ & $P(5,v)$ & $P(6,v)$\\
\hline
4&1&1&0&0\\
5&2&1&1&0\\
6&4&1&1&1\\
7&7&2&1&1\\
8&8&2&1&1\\
9&12&3&2&1\\
10&13&5&2&1\\
11&17&6&2&2\\
12&20&9&3&2\\
13&26&13&3&2\\
14&28&14&4&2\\
15&35&15&6&3\\
16&37&20&6&3\\
17&44&20&7&3\\
18&48&22&9&4\\
19&57&25&12&4\\
20&60&30&16&5\\
21&70&31&21&7\\
22&73&37&21&7\\
23&83&40&23&8\\
24&88&42&24&9\\
25&100&50&30&10\\
26&104&52&30&13\\
27&117&54&30&14\\
28&121&63&33&16\\
\hline
\end{tabular}
\caption{Maximum packing numbers for $k=4$, $k=5$ and $k=6$\label{table:packingTable}}
\end{table}
The following Lemma is a direct consequence of the definition of $P(k,v)$:

\begin{lemma}\label{lemma:packingLemma} No linear space on $v$ points contains more than  $P(k,v)$ blocks of size $k$.  
\end{lemma}

In the lemma below, the {\it regularity} of a point  $p$ ($reg_{p}$) in a linear space is the number of blocks that contain it and, for $i>0$, $deg_{i}$ is the number of points that have regularity $i$. Note that, since all the other points in the blocks containing $p$ must be distinct, if all blocks have size $r$, then $reg_{p}\leq (v-1)/(r-1)$. 

\begin{lemma}\label{lemma:blockWeight} 
If $\Lambda=(V,B)$ is a linear space and, for any block $b\in B$, the weight of $b$, $wt(b)$,  is the sum of the regularities of the points in $b$,  then (i)  block $b$ intersects $wt(b)-|b|$ other blocks, and no block can have weight more than $|B|+|b|-1$ and, (ii)  $\Sigma_{b\in B}wt(b) = \Sigma_{i>0} i^{2}deg_{i}$.
\end{lemma}

\vspace{.25cm}
\noindent{\bf Proof}  (i) Any point $p$ in $b$ with regularity $reg_{p}$ intersects $r_{p}-1$ other blocks, and no other block contains more than one element of $b$, so $b$ intersects $\Sigma_{p\in b}reg_{p}-1$ other blocks, and the result follows. (ii) If we sum the weights of the blocks, for every vertex $p$, $reg_{p}$ is counted $reg_{p}$ times. So  $\Sigma_{b\in B}wt(b) = \Sigma_{p} (reg_{p}^{2})$. For every $i>0$ there are $deg_{i}$ vertices of regularity $i$, so $\Sigma_{b\in B}wt(b) = \Sigma_{i>0}deg_{i}.i^{2}$, as required. 

\begin{lemma}\label{lemma:12.1} If $\Lambda$ is a linear space on $12$ points with $8$ blocks of size $4$, then all of the blocks intersect.
\end{lemma}

\vspace{.25cm}
\noindent{\bf Proof} There are no points of regularity greater than $3$ (or the blocks containing such a point would contain $13$ distinct points), and at most $4$ of regularity less than $3$ (if there are $5$, then $22$ spaces must be filled by $7$ points, implying at least one point of regularity greater than $3$). Suppose that there are two non-intersecting blocks, then they must each have weight at most $10$ and each contain a pair of vertices of regularity at most $2$ (from Lemma \ref{lemma:blockWeight}). It follows that $(deg_{2},deg_{3})=(4,8)$. But then two of the blocks contain only points of regularity $3$ and have weight $12$ which contradicts  Lemma \ref{lemma:blockWeight}.

%
\begin{lemma}\label{lemma:13.1}  If $\Lambda$ is a linear space on $13$ points with $9$ blocks of size $4$, then (1) there is no set of  $3$ non-intersecting blocks. (2) There are at most $3$ pairs of non-intersecting blocks.
\end{lemma}
\vspace{.25cm}
\noindent{\bf Proof} 
First observe that if there is a set of $3$ parallel blocks or more than $1$ pair of parallel blocks then no point is in at most $1$ block, or removing the point leaves a linear space on at most $12$ points with $8$ blocks of size at least $4$ with a non-intersecting pair of blocks. There is no linear space on fewer than $12$ points with $8$ blocks of size at least $4$  \cite{bettenbetten4}, and if $v=12$ we have a contradiction to Lemma \ref{lemma:12.1}. So every point is in between $2$ and $4$  blocks and  $(deg_{2},deg_{3},deg_{4})= (3,10,0)$, $(4,8,1)$, $(5,6,2)$ or $(6,4,3)$. By Lemma \ref{lemma:blockWeight} no block can have weight more than $12$ and  $\Sigma_{b\in B}wt(b) = 4deg_{2}+9deg_{3}+16deg_{4}$.  

\vspace{.2cm}
\noindent
(1) If there is a set of $3$ non-intersecting blocks, then there are $3$ blocks of weight at most $10$, and so $\Sigma_{b\in B}wt(b) \leq 102$. This is only possible if  $(deg_{2},deg_{3},deg_{4})= (3,10,0)$. But the set of $3$ parallel blocks must each contain two of the vertices of degree $2$, and so must intersect, which is a contradiction. 
(2)  If there are $4$ pairs of non-intersecting blocks then there are $8$ blocks with weight at most $11$ and $\Sigma_{b\in B}wt(b) \leq 100$, which is impossible. 

\begin{lemma}\label{lemma:16.1}
If $\Lambda$ is a linear space on at most $16$ points with $6$ blocks of size $5$ then either $(deg_{1},deg_{2},deg_{3})=(3,12,1)$, $deg_{2}=15$, or $(deg_{1},deg_{2})=(2,14)$. 
\end{lemma}

\vspace{.25cm}
\noindent{\bf Proof} No point can be in $4$ blocks. If $p$ has regularity $3$ then for every block containing $p$ there is a point of regularity $1$ (as all blocks have weight at most $10$, by Lemma \ref{lemma:blockWeight}). So we must solve the following set of equations/inequalities:
$deg_{1}+2deg_{2}+3deg_{3}=30$; $deg_{1}+deg_{2}+deg_{3}\leq 16$, and $deg_{1}\geq 3deg_{3}$. There are $3$ solutions, giving the result required.

\vspace{.25cm}
\noindent
In the remainder of this paper we present all extremal graphs of orders $20\leq v \leq 32$. In all cases, for all graphs $\Gamma$ we provide the adjacency matrix, degree sequence, edge set, relevant sets $S_{\Gamma,n}$, and any other relevant material  one of the appendices (A-M) to this paper.

\section{Case $v=20$}

\begin{theorem}\label{theorem:g20} If $\Gamma$ is a graph of girth at least $5$ on $20$ vertices with $41$ edges, then $(\delta,\Delta)=(3,5)$ and $\Gamma$ has an embedded  $S_{5,[3,3,3,3,2]}$ star.
\end{theorem}

\begin{figure}
\centering
\includegraphics[width=0.5\textwidth]{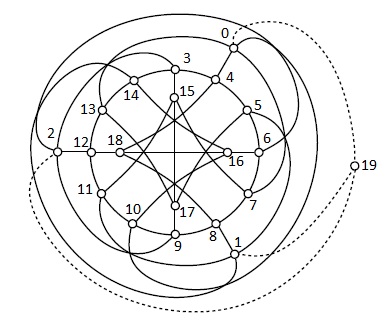}
\caption{Extremal graph on $20$ vertices with $41$ edges
\label{fig:extremal20}}
\end{figure}

\vspace{.25cm}
\noindent
{\bf Proof} There is precisely one graph on $20$ vertices with $41$ edges. This is stated in \cite{Garnick93} and the graph is given in \cite{Garnick92}. This graph is given in Figure \ref{fig:extremal20} (the numbering of the vertices is our own). There are $3$ embedded $S_{5,[3,3,3,3,2]}$ stars. Note that the sink nodes corresponding to these stars are precisely the vertices of degree 5, namely $0$, $1$ and $2$.

\begin{lemma}\label{lemma:g20}
If $\Gamma$ is a graph of girth at least $5$ on $20$ vertices with $41$ edges, then $(deg_{3},deg_{4},deg_{5}) =(1,16,3)$ and every element of $S_{\Gamma,3}$  contains $2$ vertices of degree $4$ and a child of a sink node, namely the  single vertex of degree $3$. 
\end{lemma}

\begin{figure}
\centering
\includegraphics[width=0.8\textwidth]{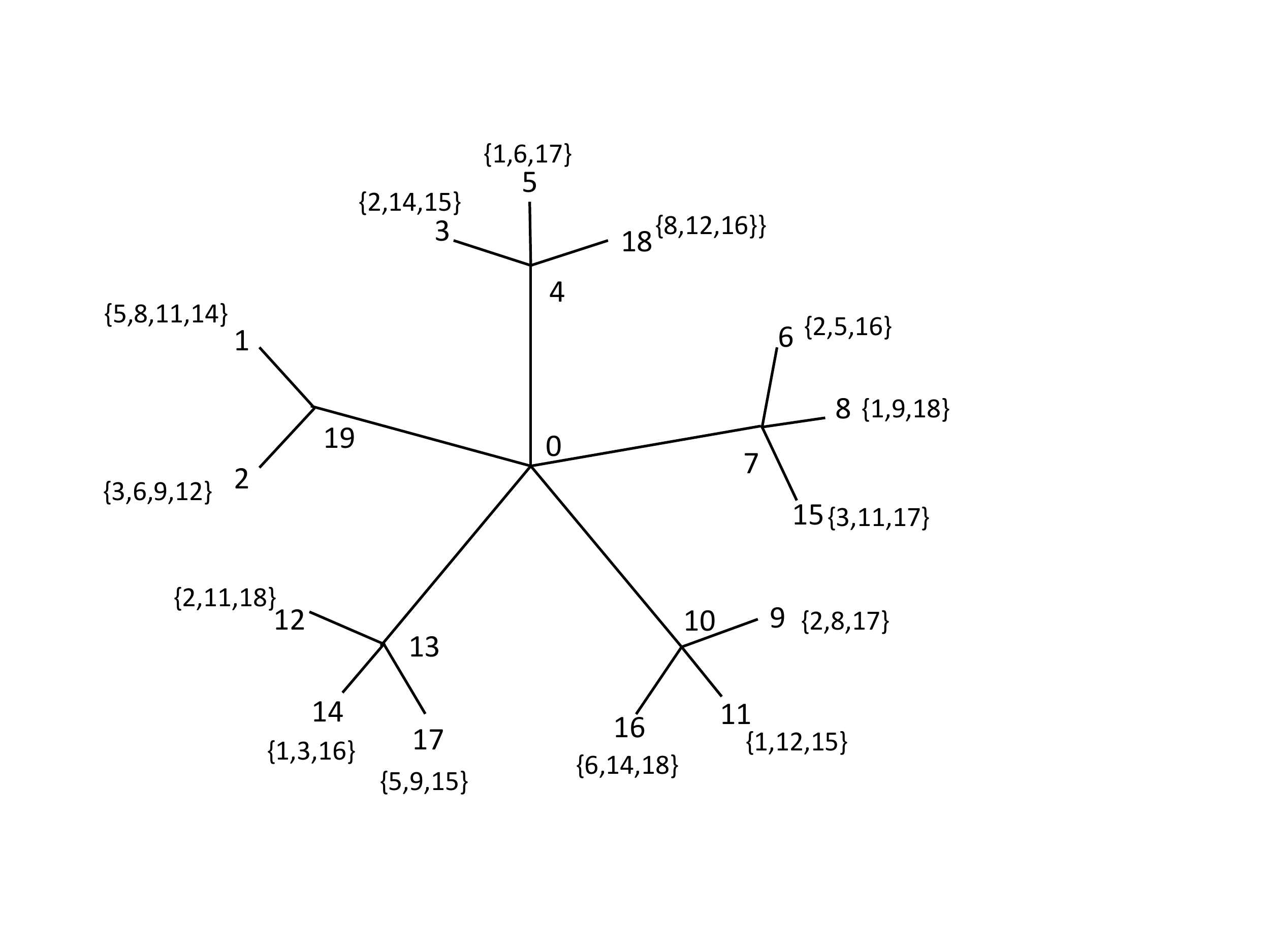}
\caption{Extremal graph on $20$ vertices with $41$ edges, as a star with root $0$\label{fig:extremal20A}}
\end{figure}

\vspace{.25cm}
\noindent
This graph is given in Figure \ref{fig:extremal20}, and is represented as an embedded star with sink node $0$ in Figure \ref{fig:extremal20A}. 
%
%
Note that $$S_{\Gamma, 3}=\{\{15,16,19\}, \{15,18,19\}, \{16,17,19\}, \{17,18,19\}\}$$
 Since vertices $15,16,17,18$ all have degree $4$, and $19$ is a child of all roots (i.e. $0$, $1$ and $2$), of degree $3$,  we are done. 

\section{Case $v=21$}

\begin{theorem}\label{theorem:g21} If $\Gamma$ is a graph of girth at least $5$ on $21$ vertices with $44$ edges, then either $(\delta,\Delta)=(3,5)$   or $(\delta,\Delta)=(4,5)$. In either case, $\Gamma$ has an embedded $S_{5,[3,3,3,3,3]}$ star.
\end{theorem}

\vspace{.25cm}
\noindent
{\bf Proof} Let $\Gamma$ be a graph of girth at least $5$ with $21$ vertices and $44$ edges. Then $\Gamma$ is extremal and, by Table  \ref{table:deltaTable} it follows that the possible values for the minimum and maximum degrees, $\delta$ and $\Delta$ are $(\delta,\Delta)=(3,5)$, $(3,6)$, and $(4,5)$.

Suppose that there is a vertex of degree $3$, $x$ say. Since $\delta=f(21)-f(20)$, by Lemma \ref{lem:inductiveExtremal}, $\Gamma$ can be constructed from the extremal graph on $20$ vertices described in Lemma \ref{lemma:g20}, $\Gamma^{\prime}$, by adding a vertex $x$ of degree $3$ and edges from $x$ to all vertices in an element of $S_{\Gamma^{\prime},3}$. Since no element of $S_{\Gamma^{\prime},3}$ contains vertices of degree $5$, it follows that $\Delta=5$. So we have $(\delta,\Delta)=(3,5)$. Since, by Lemma \ref{lemma:g20} every element of $S_{\Gamma^{\prime},3}$ contains  single vertex of $\Gamma^{\prime}$ of degree $3$ that is a child of every sink node, it follows that 
 every sink node of $\Gamma^{\prime}$ is a sink node of $\Gamma$
 of degree $5$ whose children all have degree $4$. Hence there is an embedded   $S_{5,[3,3,3,3,3]}$ star. 

If $(\delta,\Delta)=(4,5)$ then any node $x$ of degree $5$ must only have children of degree $4$ and the result follows.

\begin{figure}
\centering
\includegraphics[width=0.6\textwidth]{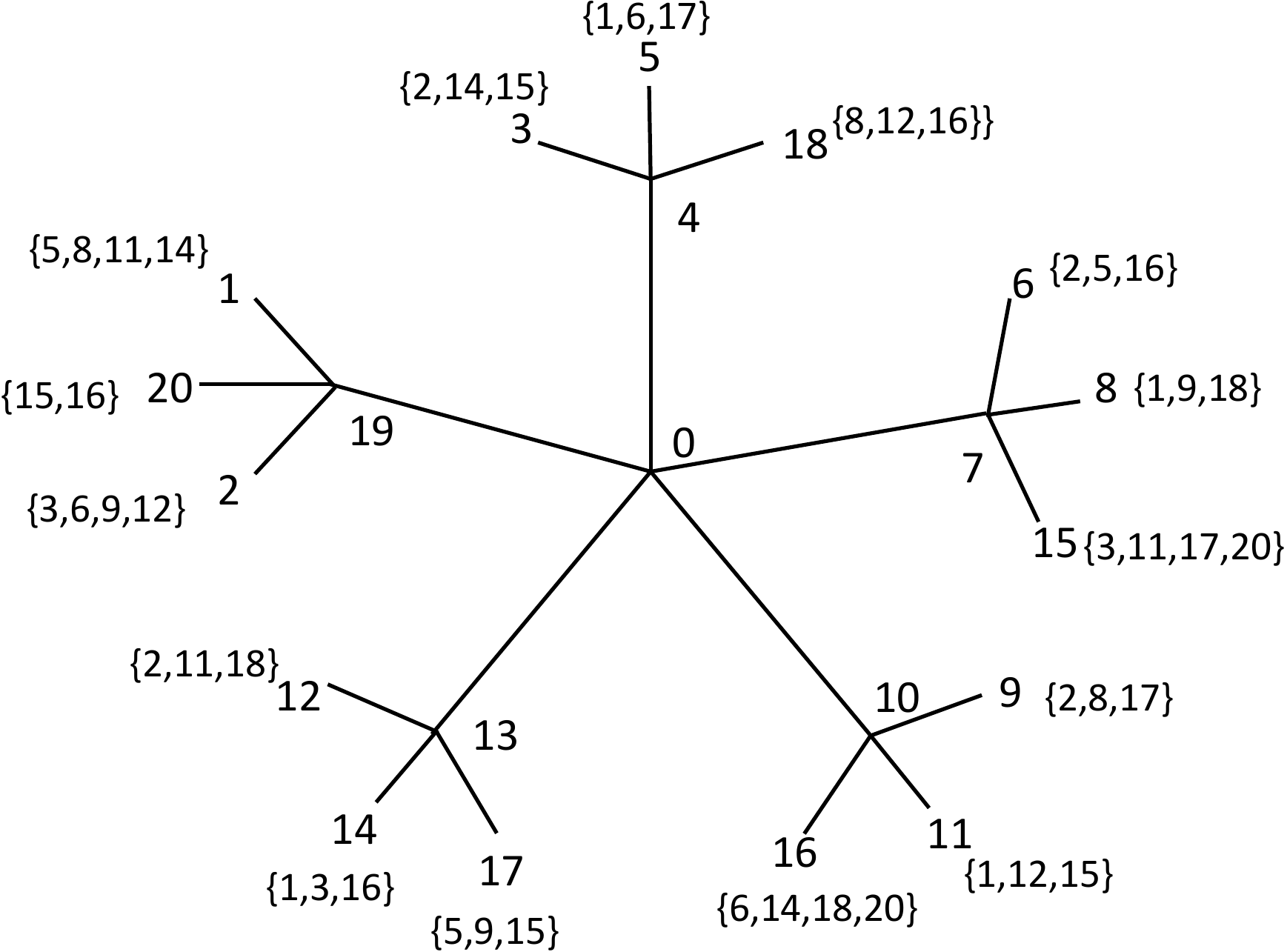}
\caption{Extremal graph on $21$ vertices with $44$ edges and $(\delta,\Delta)=(3,5)$, as a star with root $0$\label{fig:extremal21A}.}
\end{figure}

\begin{figure}
\centering
\includegraphics[width=1.0\textwidth]{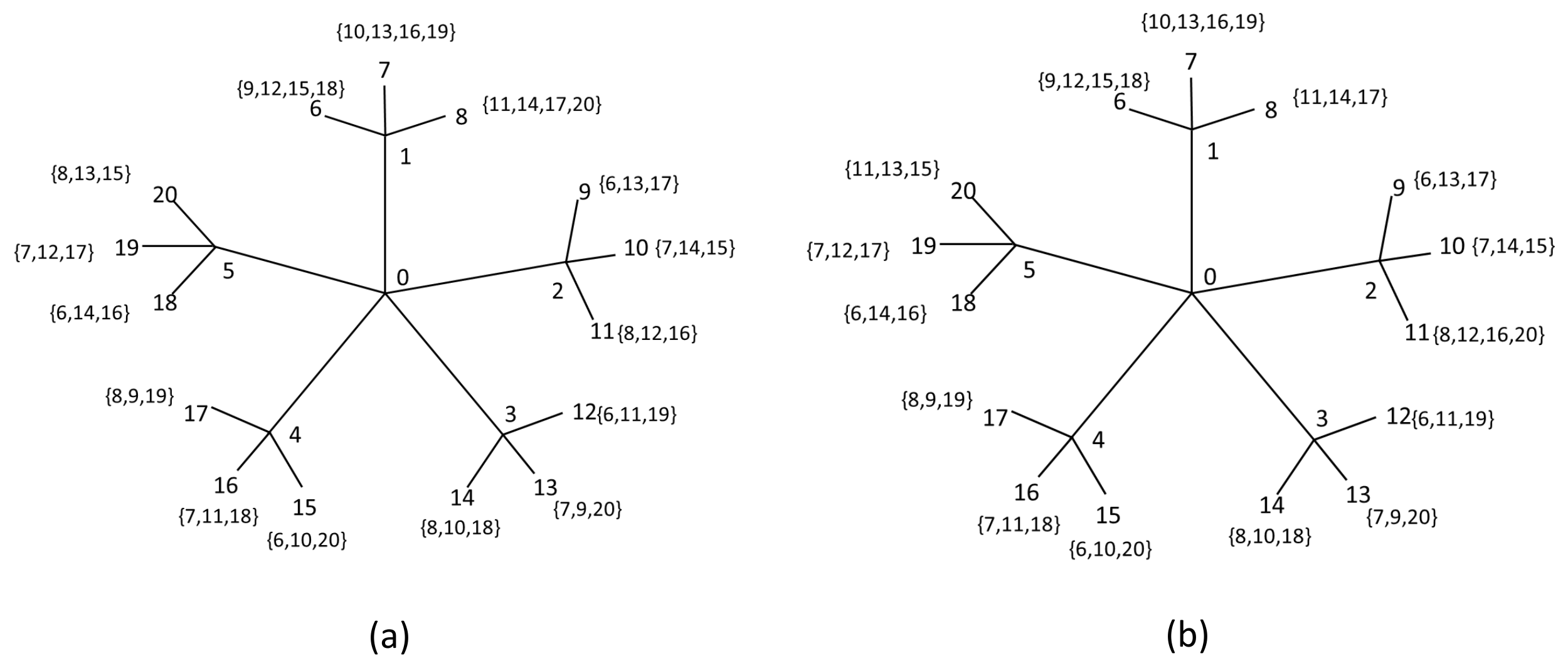}
\caption{Extremal graphs on $21$ vertices with $44$ edges and $(\delta,\Delta)=(4,5)$, as stars with root $0$\label{fig:extremal21B}}
\end{figure}

\begin{lemma}\label{lemma:g21}
If $\Gamma$ is a graph on $21$ vertices of girth at least $5$ with $44$ edges, then one of the following holds: (i) $(deg_{3},deg_{4},deg_{5})=(1,15,5)$
 and   $\Gamma$  has $3$ embedded  $S_{5,[3,3,3,3,3]}$ stars and  is isomorphic to the graph  shown in Figure \ref{fig:extremal21A}. In this case every element of $S_{\Gamma,3}$ contains $3$ vertices of degree less than $5$, and the child of a root of degree $4$. 
(ii) $(deg_{3},deg_{4},deg_{5})=(0,17,4)$, all vertices of degree $5$ are sink nodes and $\Gamma$ is isomorphic to one of the graphs shown in Figure \ref{fig:extremal21B}. In both cases $S_{\Gamma,3}$ is empty. 
\end{lemma}

\vspace{.25cm}
\noindent
{\bf Proof} From the proof of Theorem \ref{theorem:g21}, $(\delta,\Delta)=(3,5)$ or $(4,5)$. If $(\delta,\Delta)=(3,5)$ then $(deg_{3},deg_{4},deg_{5})=(1,15,5)$ and $\Gamma$ is obtained from the extremal graph $\Gamma^{\prime}$ on $20$ vertices shown in Figure \ref{fig:extremal20},  by adding an additional vertex $x$ and edges from $x$ to one of the sets $S_{\Gamma^{\prime},3}$ identified in the proof of Lemma \ref{lemma:g20}. 
Using nauty \cite{nauty}, it can be shown that all graphs constructed in this way are isomorphic to the graph  shown in Figure \ref{fig:extremal21A} (note that in the figure (and in Figure \ref{fig:extremal21B})  neighbours of leaf nodes are indicated as sets). 
 In this case the sink nodes are $0$, $1$ and $2$ and $S_{\Gamma,3}=\{ \{4,9,20\}, \{5,12,20\}, (8,13,20\}, \{17,18,19\}\}$. The result follows. 

Suppose then that there is no vertex of degree $3$. Then $(\delta,\Delta)=(4,5)$ and $(deg_{4},deg_{5})=(17,4)$. All vertices of degree $5$ are sink nodes.  By \cite{Garnick93}, we know that there are $3$ extremal graphs $(21,44)$, so there are two graphs of this type. They must be isomorphic to the two graphs shown in Figure \ref{fig:extremal21B}, which have been shown to be non-isomorphic. Note that graph (b) can be obtained from graph (a) by replacing edge $(8,20)$ with edge $(11,20)$ and in  graph (a), the $3$ grandchildren of degree $5$ have the same parent, and in graph (b) two of the grandchildren of degree $5$ have the same parent, and the other does not. In both cases $S_{\Gamma,3}$ is empty. 

\section{Case $v=22$} 
\begin{theorem}\label{theorem:g22} If $\Gamma$ is a graph of girth at least $5$ on $22$ vertices with $47$ edges, then either $(\delta,\Delta)=(3,5)$ or $(\delta,\Delta)=(4,5)$. In either case there is an embedded $S_{5,[4,3,3,3,3]}$ star. 
\end{theorem}

\vspace{.25cm}
\noindent{\bf Proof}  
Since $\Gamma$ is extremal, it follows from Table  \ref{table:deltaTable} that the values of $(\delta,\Delta)$ are either $ (3,5)$, $(3,6)$, $(3,7)$ or  $(4,5)$. 

If there is a vertex, $x$, of degree $3$, then since $\delta=f(22)-f(21)$, by Lemma \ref{lem:inductiveExtremal}, 
$\Gamma$ can be constructed from one of the extremal graphs on $21$ vertices described in Lemma \ref{lemma:g21}, $\Gamma^{\prime}$, by adding a vertex $x$ of degree $3$ and edges from $x$ to all vertices in an element of $S_{\Gamma^{\prime},3}$.
Since $S_{\Gamma^{\prime},3}$ is non-empty, it follows that  $\Gamma^{\prime}$ is the first graph defined in Lemma \ref{lemma:g21}, which has a sink node of degree $5$ with $4$ children of degree $4$. Since every element of $S_{\Gamma^{\prime},3}$ only contains elements of degree less than $5$, it follows that, in $\Gamma$, $\Delta=5$. So $(\delta,\Delta)=(3,5)$. Since every element of $S_{\Gamma^{\prime},3}$ contains a child of a sink node of $\Gamma^{\prime}$ of degree $4$, $\Gamma$ has an embedded $S_{5,[4,3,3,3,3]}$ star.

Suppose then that there are no vertices of degree $3$. Thus $(\delta,\Delta)=(4,5)$ and there are $16$ vertices of degree $4$ and
 $6$ of degree $5$.  

If there are no sink vertices, all vertices of degree $5$ have $5$ neighbours of degree $4$ (and
 $20$ grandchildren), and it follows from Lemma \ref{lem:blocks} that there is a linear space on the vertices of degree $4$. This linear space consists of $6$ blocks of size $5$. Let $B$ denote this set of blocks and $\Gamma_{4}$ the graph formed by the vertices of degree $4$ and the edges between them.  It follows from  Lemma \ref{lemma:16.1} that, since any vertex in $V_{4}$ that is in $i$ blocks in $B$ has degree $4-i$ in $\Gamma_{4}$, if $deg_{i}$ is the number of vertices of degree $i$ in $\Gamma_{4}$, $(deg_{1},deg_{2},deg_{3},deg_{4})=(1,12,3,0)$, $(0,15,0,1)$ or $(0,14,2,0)$.  If, in $\Gamma_{4}$, there is a chain of $4$ vertices of degree $2$,  $a-b-c-d$ say, then no two of $a$, $b$, $c$ can be in the same block of $B$, or there is a cycle of length at most $4$. Hence two of the blocks  of $B$ contain $a$, two contain $b$ and two contain $c$. Now similarly $d$ is not in a block in $B$ with $b$ or $c$, so must be in both blocks of $B$ containing $a$, which is not possible. So we assume there is no such chain and it follows that $(deg_{1},deg_{2},deg_{3},deg_{4})\neq (0,15,0,1)$.  By a similar argument, if $a$ is a vertex of degree $1$ then there is no chain $a-b-c$ where $b$ and $c$ have degree $2$.  Hence, if $(deg_{1},deg_{2},deg_{3},deg_{4})=(1,12,3,0)$ and $a$ the vertex of degree $1$, with neighbour $b$, either$b$ has degree $3$, or $b$ has degree $2$ and its other neighbour $c$ has degree $3$. In either case, there is a vertex of regularity $1$ in the blocks that is in no block with $a$, hence one of the blocks on $a$ contains no vertex of regularity $1$ which is not possible (see the proof of Lemma \ref{lemma:16.1}). 

 So $(deg_{1},deg_{2},deg_{3},deg_{4})=(0,14,2,0)$. 

Let $x$ and $y$ be the vertices of degree $3$. Since there can be no component consisting entirely of vertices of degree $2$ (or it would either contain a cycle of length $3$, or a chain of at least $4$ vertices of degree $2$) either $x$ and $y$ are in different components, or $\Gamma_{4}$ consists of a single component. If $x$ and $y$ are in different components of $\Gamma_{4}$ then their components contain a single vertex of odd degree and all other vertices of even degree, which is not possible (the sum of the degrees in a component must be even). So $\Gamma_{4}$ consists of a single component and either there is a single path (of length at most $4$) from $x$ to $y$ and cycles containing $x$ and $y$ respectively, or $3$ paths from $x$ and $y$, together containing all of the vertices of degree $2$. In either case, there is a chain of at least $4$ vertices of degree $2$.

\begin{lemma}\label{lemma:g22} There are $3$ nonisomorphic graphs of girth at least $5$ on $22$ vertices with $47$ edges. In one,  $(\delta,\Delta)=(4,5)$  and $(deg_{4},deg_{5})=(16,6)$.  In the other two, $(\delta,\Delta)=(3,5)$ and $(deg_{3},deg_{4},deg_{5})=(2,12,8)$ and $(1,14,7)$ respectively. In all cases, every element of $S_{\Gamma,3}$ contains $3$ vertices of degree less than $5$, and contains a child of a sink node of $\Gamma$. 
\end{lemma}

\vspace{.25cm}
\noindent
{\bf Proof} By Theorem \ref{theorem:g22} any such graph has $(\delta,\Delta)=(3,5)$  or $(\delta,\Delta)=(4,5)$ and there is an embedded   $S_{5,[4,3,3,3,3]}$ star. Restricting search to those cases produced a set of graphs. Isomorph elimination using nauty \cite{nauty} revealed there to be three such graphs, $\Gamma_{0}$ with $(\delta,\Delta)=(4,5)$, and $(deg_{4},deg_{5})=(16,6)$, and  $\Gamma_{1}$ and $\Gamma_{2}$ with $(\delta,\Delta)=(3,5)$ and $(deg_{3},deg_{4},deg_{5})=(2,12,8)$ and $(1,14,7)$ respectively.

\begin{figure}
\centering
\includegraphics[width=0.7\textwidth]{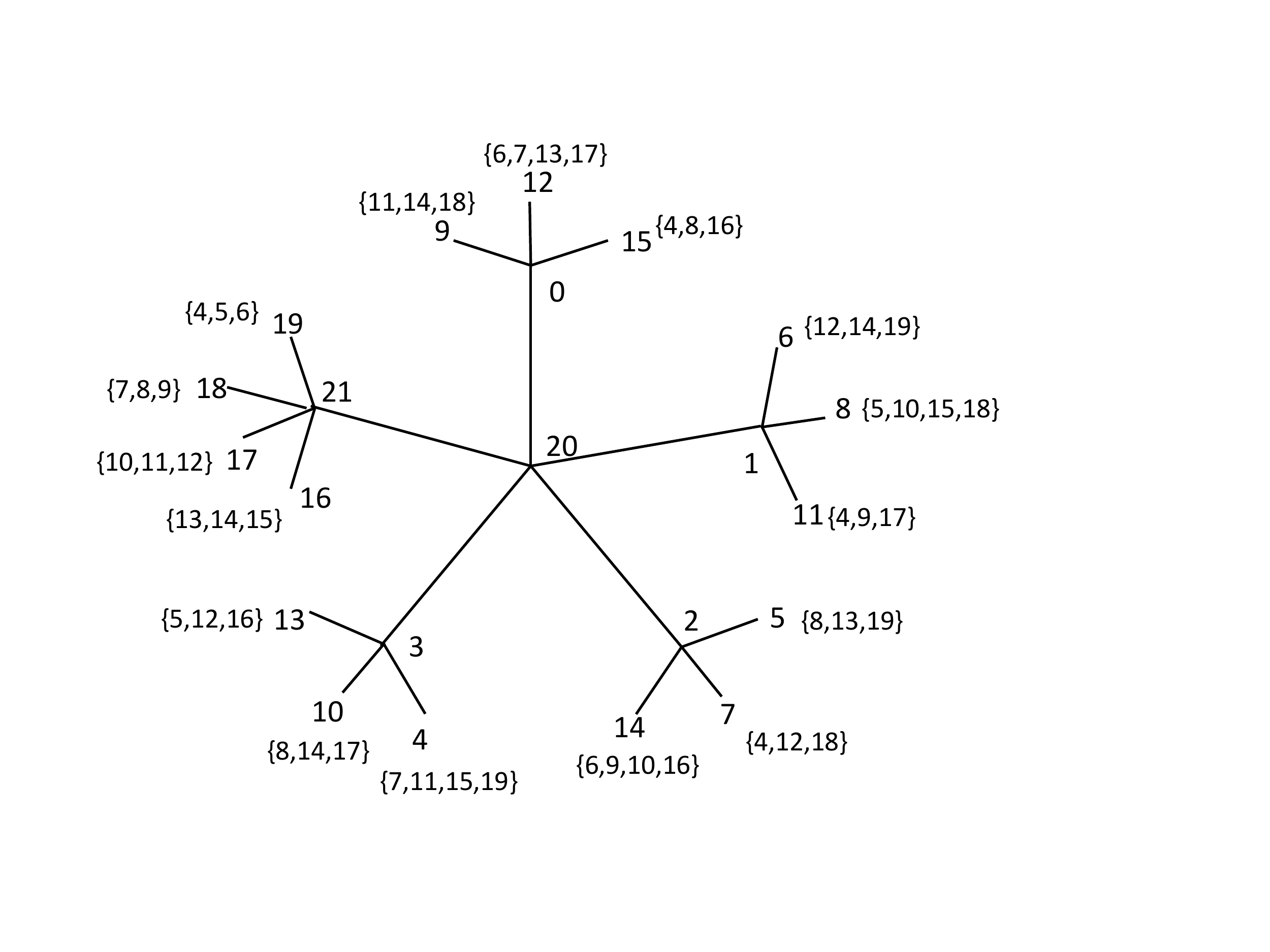}
\caption{Extremal graph $\Gamma_{0}$ on $22$ vertices with $44$ edges as a star, $(\delta,\Delta)=(4,5)$\label{fig:extremal22A0}}
\end{figure}

\begin{figure}
\centering
\includegraphics[width=0.7\textwidth]{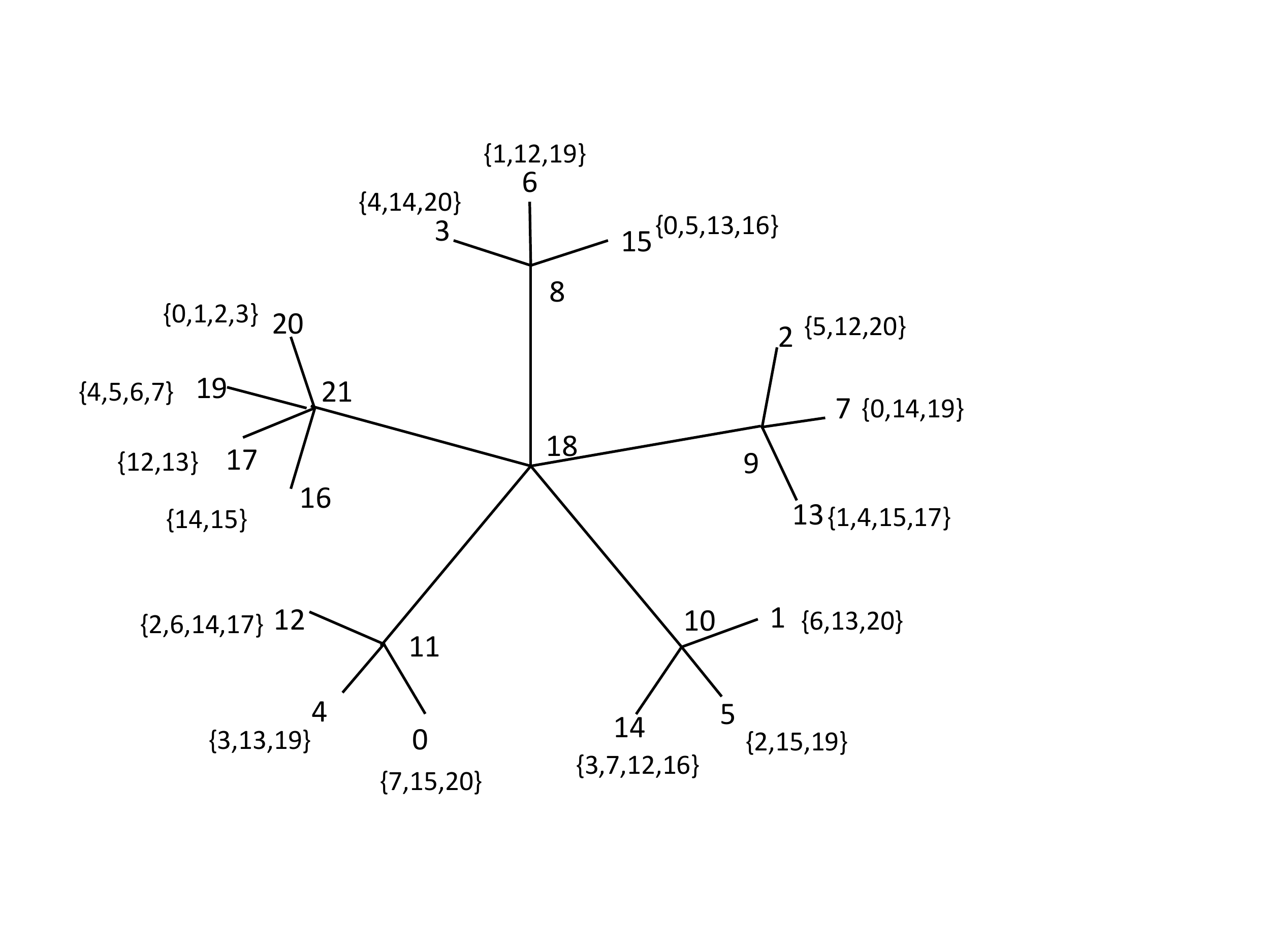}
\caption{Extremal graph $\Gamma_{1}$ on $22$ vertices with $44$ edges as a star, $(\delta,\Delta)=(3,5)$\label{fig:extremal22A1}}
\end{figure}

\begin{figure}
\centering
\includegraphics[width=0.7\textwidth]{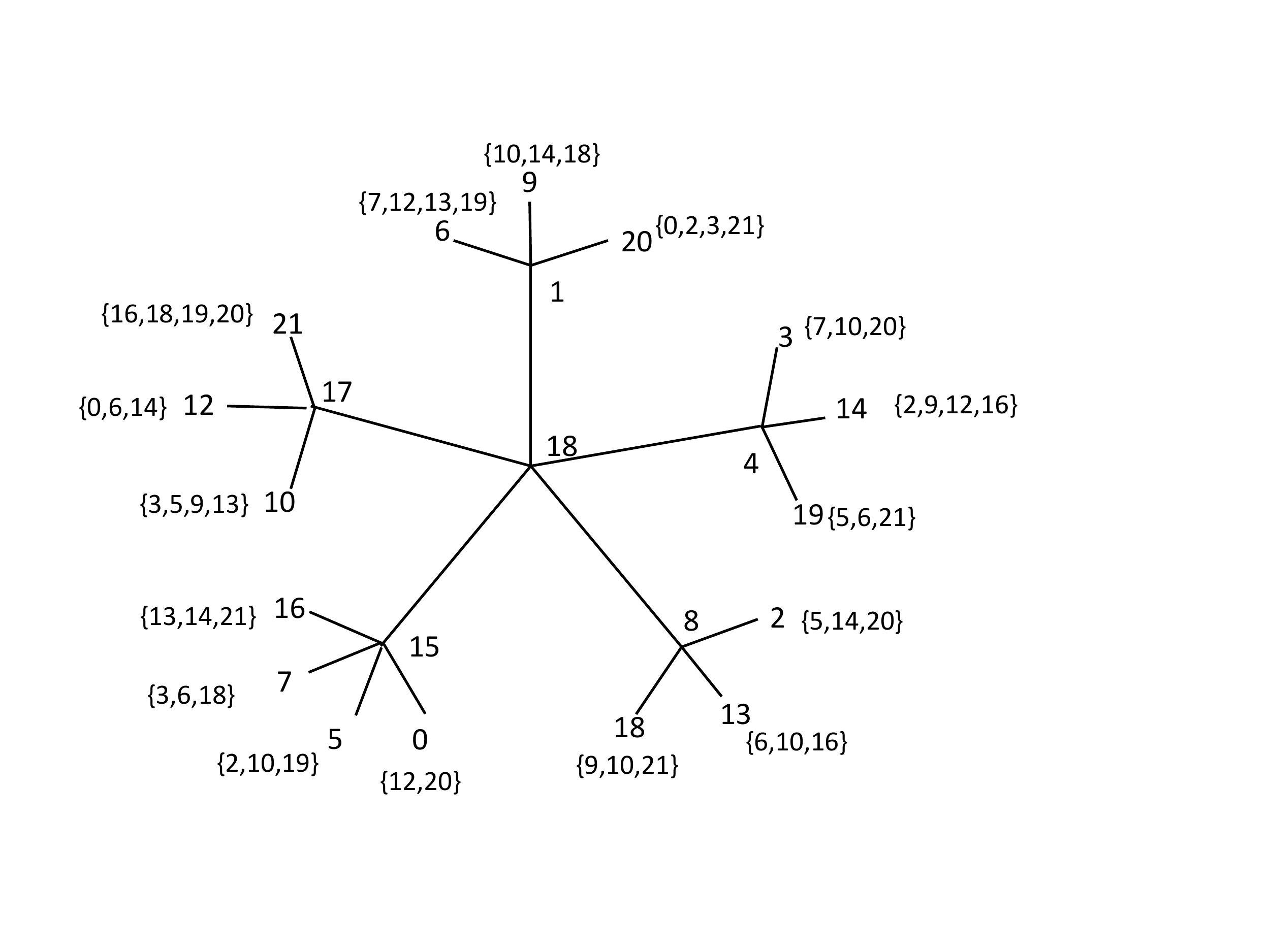}
\caption{Extremal graph $\Gamma_{2}$ on $22$ vertices with $44$ edges as a star, $(\delta,\Delta)=(3,5)$\label{fig:extremal22A2}}
\end{figure}

These graphs are shown in Figures \ref{fig:extremal22A0}-\ref{fig:extremal22A2}. In each case the graph is drawn as a star about a sink node with $4$ children of degree  $4$ and one of degree $5$. The vertices adjacent to each of the grandchildren is indicated as a set.

\noindent
Examination of these graphs using a computer reveals that in each case every element of $S_{\Gamma,3}$ contains no vertex of degree greater than $4$, and contains at least one vertex that is the child of a sink node. See Appendix C for details. 

\begin{lemma}\label{lemma:g22A} If $\Gamma$ is a graph of girth at least $5$  on $22$ vertices with $47$ edges then $S_{\Gamma,3}$ contains at most $4$ non-intersecting sets. 
If all vertices of  degree $3$ (if any) are contained in a single element of $S_{\Gamma,3}$ then $\Gamma$ is $\Gamma_{1}$ or $\Gamma_{2}$ of Lemma \ref{lemma:g22}. If in addition there are two vertices of degree $4$ which have fewer than $3-deg_{3}$ neighbours of degree $4$ between them, $\Gamma$ is $\Gamma_{2}$ . In this case there are three vertices of degree $4$ which have fewer than two neighbours of degree $4$ and they each have one such neighbour, $p$, $q$ and $r$ respectively. No two of $p$, $q$ and $r$ are in a common element of $S_{\Gamma,3}$ with the vertex of degree $3$. 

\end{lemma}

\vspace{.25cm}
\noindent
{\bf Proof} From examination of the graphs in Appendix C,  $S_{\Gamma,3}$ has at most $4$ non-intersecting sets and if $S_{\Gamma,3}$ contains  a set $X$ of $3$ non-intersecting sets then $\Gamma$ is graph $\Gamma_{0}$ or $\Gamma_{2}$. Examination of the vertices of degree $4$  shows that all vertices of degree $4$ have at least two neighbours of degree $4$ unless $\Gamma$ is $\Gamma_{2}$  and the vertex is $1$, $16$, or $17$, which have neighbours of degree $4$: $12$, $13$ and $9$ respectively. No two of these vertices are in an element of  $S_{\Gamma,3}$ with the vertex of degree $3$ (i.e. vertex $0$).

\section{Case $v=23$} 

\begin{theorem}\label{theorem:g23} If $\Gamma$ is a graph of girth at least $5$  on $23$ vertices with $50$ edges then $(\delta,\Delta)=(3,5)$ and there is an embedded  $S_{5,[4,4,3,3,3]}$ star or  $(\delta,\Delta)=(4,5)$ and either there is an embedded  $S_{5,[4,4,3,3,3]}$ star  or there is no embedded  $S_{5,[4,4,3,3,3]}$ star and there is an embedded $S_{5,[4,3,3,3,3]}$ star. 
\end{theorem}

\vspace{.25cm}
\noindent
{\bf Proof} 
Since $\Gamma$ is extremal, it follows from Table  \ref{table:deltaTable} that the values of $(\delta,\Delta)$ are either $ (3,5)$, $(3,6)$, $(3,7)$ or  $(4,5)$. 

If there is a vertex, $x$, of degree $3$, then since $\delta=f(23)-f(22)$, by Lemma \ref{lem:inductiveExtremal}, 
$\Gamma$ can be constructed from one of the extremal graph on $22$ vertices described in Lemma \ref{lemma:g22}, $\Gamma^{\prime}$, by adding a vertex $x$ of degree $3$ and edges from $x$ to all vertices in an element of $S_{\Gamma^{\prime},3}$. 

 Since in each case $S_{\Gamma^{\prime},3}$ contains no vertex of degree greater than $4$, (by Lemma \ref{lemma:g22}), in $\Gamma$ $\Delta=5$. Hence $(\delta, \Delta)=(3,5)$. In addition, since $S_{\Gamma^{\prime},3}$ contains a child of a sink node of $\Gamma^{\prime}$, it follows that $\Gamma$ has a sink node of degree $5$ with $2$ children of degree $5$ and $3$ of degree $4$. 

If there is no vertex of degree $3$ then $(\delta,\Delta)=(4,5)$ and $(deg_{4},deg_{5})=(15,8)$. If no vertex of degree $5$ has any neighbour of degree $5$ then there is a linear space on $15$ points with $8$ blocks of size $5$, contradicting Lemma \ref{lemma:packingLemma}. Since every vertex of degree $5$ has at least one neighbour of degree $5$, the result follows.

\begin{lemma}\label{lemma:g23}  There are $7$ non-isomorphic graphs on $23$ vertices with $50$ edges. All of the graphs contain embedded $S_{5,[4,4,3,3,3]}$ stars. For $3$ of the graphs, $S_{\Gamma,4}$ is empty. The remaining $4$ graphs consist of:
\begin{itemize}
\item $1$ graph with $(\delta,\Delta)=(4,5)$ and $(deg_{4},deg_{5})=(15,8)$. Every element of $S_{\Gamma,4}$ contains $4$ vertices of degree $4$, including at least one child of a sink node 
\item $2$ graphs with  $(\delta,\Delta)=(3,5)$ and $(deg_{3},deg_{4},deg_{5})=(1,13,9)$. Every element of $S_{\Gamma,4}$ contains the vertex of degree $3$ and $3$ vertices of degree $4$, including a  child of a sink node.
\item $1$ graph with  $(\delta,\Delta)=(3,5)$ and $(deg_{3},deg_{4},deg_{5})=(2,11,10)$. Every element of $S_{\Gamma,4}$ contains the two vertices of degree $3$ and $2$ vertices of degree $4$, including a child of a sink node. 
\end{itemize}
\end{lemma} 

\vspace{.25cm}
\noindent
{\bf Proof} (1) By Theorem \ref{theorem:g23} any such graph has $(\delta,\Delta)=(3,5)$  and there is an embedded  $S_{5,[4,4,3,3,3]}$ star or  $(\delta,\Delta)=(4,5)$ and either there is an embedded  $S_{5,[4,4,3,3,3]}$ star  or there is no embedded  $S_{5,[4,4,3,3,3]}$ star and there is an embedded $S_{5,[4,3,3,3,3]}$ star. Restricting search to those cases produced a set of
 graphs. Isomorph elimination using Nauty \cite{nauty} revealed there to be $7$ such graphs, two with 
$(\delta,\Delta)=(4,5)$, and $(deg_{3},deg_{4},deg_{5})=(0,15,8)$, and  $5$ with $(\delta,\Delta)=(3,5)$, two of which
 have  $(deg_{3},deg_{4},deg_{5})=(2,11,10)$ and three of which have $(deg_{3},deg_{4},deg_{5})=(1,13,9)$. In all cases there is an embedded  $S_{5,[4,4,3,3,3]}$ star, and using a computer we have verified the structure of $S_{\Gamma,4}$ in each case. See Appendix D for details. 

\begin{lemma}\label{lemma:g23A} Suppose that $\Gamma$ is a graph of girth at least $5$ on $23$ vertices with $(\delta,\Delta)=(4,5)$ and  $\Gamma_{4}=(V_{4},E_{4})$ the subgraph of $\Gamma$ on the vertices of degree $4$, then there is at most one vertex in $V_{4}$ that is in no edges in $E_{4}$. If two vertices have degrees $0$ and $1$ in $V_{4}$, then they have a common neighbour of degree $5$ in $\Gamma$. 
\end{lemma}

\vspace{.25cm}
\noindent
{\bf Proof} From the proof of Lemma \ref{lemma:g23} there are two such graphs with  $(\delta,\Delta)=(4,5)$. Only the first of these contains a vertex of degree $4$ that has no neighbours of degree $4$, (vertex $a$ say). There are two vertices that have degree $1$, and each of these are in an edge in $\Gamma$ with $a$. See Appendix D for details.

\section{Case $v=24$}

\begin{theorem}\label{theorem:g24} If $\Gamma$ is a graph of girth at least $5$  on $24$ vertices with $54$ edges then $(\delta,\Delta)=(4,5)$ and there is an embedded $S_{5,[4,4,4,3,3]}$ star.
\end{theorem}

\vspace{.25cm}
\noindent
{\bf Proof}
From Table  \ref{table:deltaTable} we know that in this case all vertices have degree $4$ or $5$. There are $12$ vertices of degree
 $4$ and $12$ vertices of degree $5$. Since $\delta=f(24)-f(23)$, by Lemma \ref{lem:inductiveExtremal}, 
$\Gamma$ can be constructed from one of the extremal graph on $23$ vertices described in Lemma \ref{lemma:g23}, $\Gamma^{\prime}$, by adding a vertex $x$ of degree $4$ and edges from $x$ to all vertices in an element of $S_{\Gamma^{\prime},4}$.
In all cases, $S_{\Gamma^{\prime},4}$ includes precisely one child of degree $4$ of any central node of an embedded $S_{5,[4,4,3,3,3]}$ star. It follows that $\Gamma$ contains an embedded $S_{5,[4,4,4,3,3]}$ star.

\begin{lemma}\label{lemma:g24} (1) There is $1$ graph $\Gamma$ of girth at least $5$  on $24$ vertices with $54$ edges. 
(2) Any element of $S_{\Gamma,3}$ contains only vertices of degree $4$. (3) There are no $4$ non-intersecting sets in $S_{\Gamma,3}$ for which there are are two pairs of sets with no edge between any two vertices in the same pair. (4) If $X$ is a set of $3$ non-intersecting elements of  $S_{\Gamma,3}$ then at least two of the three vertices of degree $4$ that are not in any set in $X$ are in an edge with an element in a set in $X$.
\end{lemma}

\begin{figure}
\centering
\includegraphics[width=0.8\textwidth]{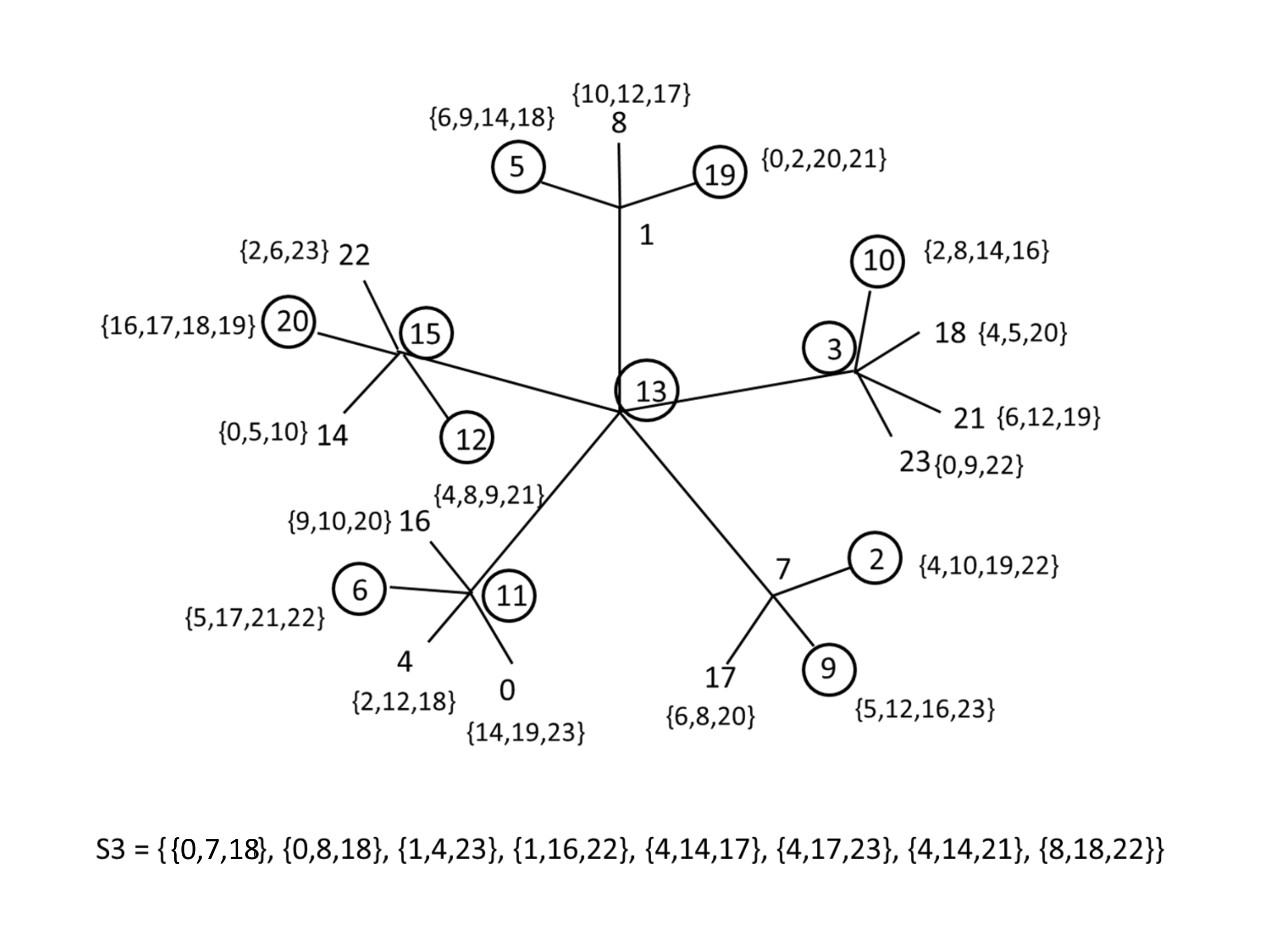}
\caption{Extremal graph on $24$ vertices
\label{fig:graph24}}
\end{figure}

\vspace{.25cm}
\noindent
{\bf Proof}(1) By Theorem \ref{theorem:g24}  $(\delta,\Delta)=(4,5)$ and there is an embedded $S_{5,[4,4,4,3,3]}$ star.  Restricting search to this case produced a set of graphs. Isomorph elimination using nauty \cite{nauty} revealed there to be one such graph, with $2$ sink nodes. This is shown in Figure \ref{fig:graph24} as a star about vertex $13$, together with set $S_{\Gamma,3}$. Clearly all elements of $S_{\Gamma,3}$ contain only vertices of degree $4$. (2) and (3) follow from observation, see Appendix E for details. 

\section{Case $v=25$}

\vspace{.25cm}
\begin{theorem}
\label{theorem:g25}
If $\Gamma$ is a graph of girth at least $5$ on $25$ vertices with $57$ edges then  $(\delta,\Delta)=(3,5)$, $(4,5)$ or $(4,6)$. If $(\delta,\Delta)=(3,5)$ or $(4,5)$ then $\Gamma$ has an embedded  $S_{5,[4,4,4,4,3]}$ star or there is no embedded $S_{5,[4,4,4,4,3]}$ star and there is an embedded $S_{5,[4,4,4,3,3]}$ star.  If  $(\delta,\Delta)=(4,6)$ then $1\leq deg_{6}\leq 2$ and there is an embedded $S_{6,[3,3,3,3,3,3]}$ star. 
\end{theorem}

\vspace{.25cm}
\noindent
{\bf Proof} We know from Table  \ref{table:deltaTable} that in this case the possible values for $(\delta,\Delta)$ are $(3,5)$, $(3,6)$, $(3,7)$, $(3,8)$, $(4,5)$ and $(4,6)$. If there is a vertex, $x$, of degree $3$, then since $\delta=f(25)-f(24)$, by Lemma \ref{lem:inductiveExtremal}, 
$\Gamma$ can be constructed from the extremal graph on $24$ vertices described in Lemma \ref{lemma:g24}, $\Gamma^{\prime}$, by adding a vertex $x$ of degree $3$ and edges from $x$ to all vertices in an element of $S_{\Gamma^{\prime},3}$.
Since $S_{\Gamma^{\prime},3}$ only contains vertices of degree $4$, it follows that $(\delta,\Delta)=(3,5)$. Since $\Gamma^{\prime}$ contains a star $[5,[4,4,4,3,3])$, $\Gamma$ does too, and the result for $\delta=3$ holds.

If $(\delta,\Delta)=(4,5)$ then $(deg_{4},deg_{5})=(11,14)$ and the sets of vertices of degree $4$ adjacent to each vertex of degree $5$ form a linear space (see Lemma \ref{lem:blocks}) $\Lambda$. If no vertex of degree $5$ has more than $2$ neighbours of degree $5$ then $\Lambda$ consists of $14$ blocks of size at least $3$. It follows from \cite{bettenbetten4} that $\Lambda$ has at least $11$ blocks of size $3$ (and, correspondingly, at least $11$ vertices of degree $5$ have two neighbours of degree $5$. 

\begin{figure}
\centering
\includegraphics[width=0.75\textwidth]{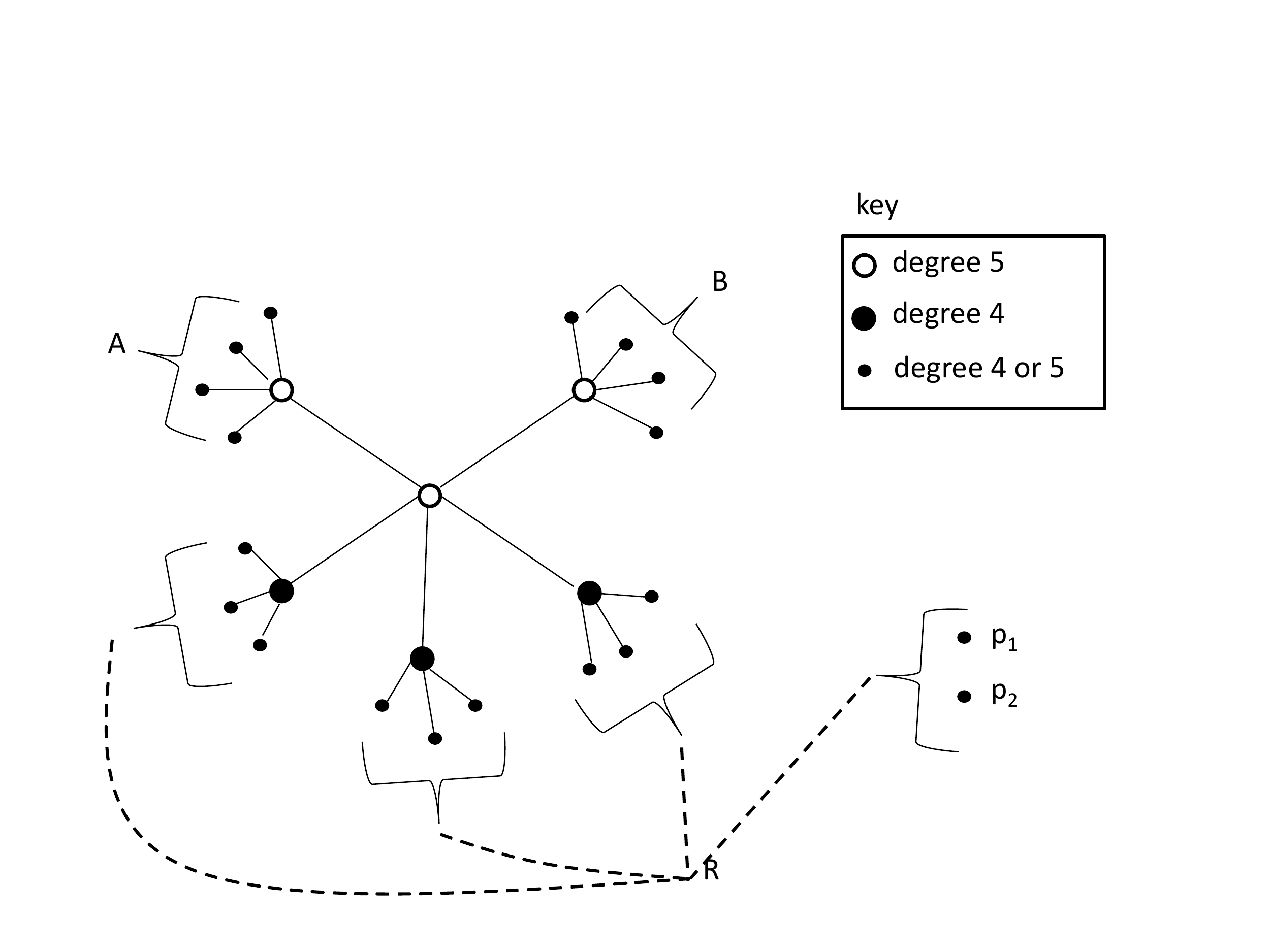}
\caption{$25$, embedded $S_{5,[4,4,3,3,3]}$ star
\label{fig:star25false}}
\end{figure}

Let $x$ be such a  vertex of degree $5$.  Let S be the embedded $S_{5,[4,4,3,3,3]}$ star with $x$ as root (see Figure \ref{fig:star25false}). Note that  $A$ and $B$ denote the sets of children of the two children of $x$ of degree $5$, and $R$ the set consisting of the remaining grandchildren of $x$, and the vertices $p_{1}$  and $p_{2}$ that are not contained in the star. Since we have assumed that all vertices of degree $5$ have at most one child of degree $5$, $A$ and $B$ contain at most $1$ vertex each of degree $5$ and $R$ contains at least $9$ vertices of degree $5$ and at most $2$ of degree $4$.  Now if $A$ contains a vertex of degree 5 then it must have at least 2 neighbours in $R$ that do not have degree $5$. So $R$ contains exactly $9$ vertices of degree $5$ and $2$ of degree 4. Then $B$ must contain a  vertex of degree $5$, which must be adjacent to the same two elements of $R$, and there is a square. It follows that $A$ and $B$ contain no vertices of degree $5$, and every element of $R$ has degree $5$. But then $p_{1}$ must be adjacent to $3$ elements of $A\cup B$ (as it can have at most $2$ neighbours of degree $5$), which is impossible, as there will then be a square. So there is at least one vertex of degree $5$ that has more than $2$ neighbours of degree $5$ and so $\Gamma$ either has an embedded 
$S_{5,[4,4,4,4,3]}$ star, or has no embedded $S_{5,[4,4,4,4,3]}$ star but has an embedded $S_{5,[4,4,4,3,3]}$ star.

If $(\delta,\Delta)=(4,6)$, any vertex of degree $6$ must have $6$ neighbours of degree $4$, so the result follows. 

\begin{lemma}\label{lemma:g25} There are $6$ graphs of girth at least $5$ on $25$ vertices with $57$ edges. Three of the graphs have $(\delta,\Delta)=(3,5)$ and $(deg_{3},deg_{4},deg_{5})=(1,19,5)$. One of these graphs has no sink node and set $S_{\Gamma,4}$ is empty. The other two have an embedded $S_{5,[4,4,4,4,3]}$ star and non-empty  sets $S_{\Gamma,4}$ where every element of $S_{\Gamma,4}$ contains a child of degree $4$ of a sink node, the vertex of degree $3$, and $2$ further vertices of degree $4$.
In the remaining three graphs, $(\delta,\Delta)=(4,5)$, $(deg_{4},deg_{5})=(11,14)$ and there is an embedded $S_{5,[4,4,4,4,3]}$ star. One of these graphs has empty set  $S_{\Gamma,4}$, and for the other two $S_{\Gamma,4}$ contains a child of degree $4$ of a sink node and $3$ other vertices of degree $4$. In addition, for each of $\Gamma_{3},\ldots, \Gamma_{5}$, $S_{\Gamma,3}$ does not contain $3$ non-intersecting elements consisting entirely of vertices of degree $4$.  
\end{lemma}

\vspace{.25cm}
\noindent
{\bf Proof} Restricting our search to graphs with values of $(\delta,\Delta)$, with the identified embedded stars from Theorem \ref{theorem:g25}  produced a set of graphs.   Isomorph elimination using nauty \cite{nauty} revealed there to be six graphs, $A_{0},\ldots, A_{5}$. Using a computer the graphs were shown to satisfy the stated properties. For details see Appendix F.

\begin{lemma}\label{lemma:g25A} If $\Gamma$ is a graph of girth at least $5$ on $25$ vertices with $57$ edges and 
$X_{1}$, $X_{2}$, $X_{3}$ parallel elements of of $S_{\Gamma,3}$ then
\begin{enumerate}
\item if $(\delta,\Delta)=(4,5)$ one of the sets must contain a vertex of degree $5$
\item it is not possible that $X_{1}$ contains all vertices of degree $3$ and there are no edges between $X_{1}$ and $X_{2}$, or between $X_{1}$ and $X_{3}$.
\end{enumerate}
\end{lemma}

\vspace{.25cm}
\noindent
{\bf Proof} Follows from examination of graphs $\Gamma_{0}\ldots \Gamma_{5}$ (and corresponding sets $S_{\Gamma,3}$) as shown in Appendix F. 

\section{Case $v=26$}

\vspace{.25cm}
\begin{theorem}
\label{theorem:g26}
If $\Gamma$ is a graph of girth at least $5$ on $26$ vertices with $61$ edges then  $(\delta,\Delta)=(4,5)$ and $\Gamma$ has an embedded $S_{5,[4,4,4,4,4]}$ star. 
\end{theorem}

\vspace{.25cm}
\noindent
{\bf Proof} We know from Table  \ref{table:deltaTable} that in this case the possible values for $(\delta,\Delta)$ are $(4,5)$ and $(4,6)$. Since in all cases $\delta= f(26)-f(25)$, 
by Lemma \ref{lem:inductiveExtremal}, $\Gamma$ can be constructed from one of the extremal graphs on $25$ vertices described in Lemma \ref{lemma:g25}, $\Gamma^{\prime}$, by adding a vertex $x$ of degree $4$ and edges from $x$ to all vertices in an element of $S_{\Gamma^{\prime},4}$.
By Lemma \ref{lemma:g25},  it follows that, in $\Gamma$, $(\delta,\Delta)=(4,5)$ and $\Gamma$ has an embedded $S_{5,[4,4,4,4,4]}$ star.

\begin{lemma}\label{lemma:g26} There are $2$ graphs of girth at least $5$ on $26$ vertices with $61$ edges, both have $(\delta,\Delta)=(4,5)$. In one of these graphs, set $S_{\Gamma,4}$ is empty. In the other,  $S_{\Gamma,4}$ contains $1$ set, which  consists of $4$ vertices of degree $4$, none of which is a child of a sink node. In one graph there is a unique set $X=\{X_{1},X_{2}\}$ containing two non-intersecting elements of $S_{\Gamma,3}$ that each consist of elements of degree $4$, and, for all $x_{1}\in X_{1}$, $x_{2}\in X_{2}$ there is no edge $(x_{1},x_{2}) \in \Gamma$. In the other graph there is no such set $X$. 
\end{lemma}

\vspace{.25cm}
\noindent
{\bf Proof} From Theorem \ref{theorem:g26}  all extremal graphs on $26$ vertices can be constructed by adding a new vertex $x$, and edges to a set of $4$ vertices that is a set in $S_{\Gamma^{\prime},4}$, where $\Gamma^{\prime}$ is one of the $6$ extremal graphs on $25$ vertices. This produced a set of $6$ graphs.   Isomorph elimination using nauty \cite{nauty} revealed there to be two distinct graphs, $A_{0}, A_{1}$. Using a computer the graphs were shown to satisfy the stated properties. See Appendix G for details. 

\section{Case $v=27$}
\vspace{.25cm}
\begin{theorem}
\label{theorem:g27}
If $\Gamma$ is a graph of girth at least $5$ on $27$ vertices with $65$ edges then  $(\delta,\Delta)=(4,5)$ and $\Gamma$ has an embedded  $S_{5,[4,4,4,4,4]}$ star. 
\end{theorem}

\vspace{.25cm}
\noindent
{\bf Proof} We know from Table  \ref{table:deltaTable} that in this case the possible values for $(\delta,\Delta)$ are $(4,5)$ and  $(4,6)$. Since in all cases $\delta= f(27)-f(26)$, 
by Lemma \ref{lem:inductiveExtremal}, $\Gamma$ can be constructed from an extremal graph on $26$ vertices described in Lemma \ref{lemma:g26}, $\Gamma^{\prime}$, by adding a vertex $x$ of degree $4$ and edges from $x$ to all vertices in an element of $S_{\Gamma^{\prime},4}$.
By Lemma \ref{lemma:g26} it follows that, in $\Gamma$, $(\delta,\Delta)=(4,5)$ and (since any non-empty set $S_{\Gamma^{\prime},4}$ does {\it not} contain the child of a sink node), $\Gamma$  has an embedded  $S_{5,[4,4,4,4,4]}$ star. 

\begin{lemma}\label{lemma:g27} There is  $1$ graph of girth at least $5$ on $27$ vertices with $65$ edges. In this graph $(\delta,\Delta)=(4,5)$ and $S_{\Gamma,3}$ contains $5$ sets, each of which contains $1$ vertex of degree $4$ and $2$ of degree $5$.  All vertices of degree $5$ in elements of $S_{\Gamma,3}$ are roots of an embedded  $S_{5,[4,4,4,4,3]}$ star. Set $S_{\Gamma,4}$ is empty.
\end{lemma}

\vspace{.25cm}
\noindent
{\bf Proof} From Theorem \ref{theorem:g27}  all extremal graphs on $27$ vertices can be constructed by adding a new vertex $x$, and edges to a set of $4$ vertices that is a set in $S_{\Gamma^{\prime},4}$, where $\Gamma^{\prime}$ is one of the $2$ extremal graphs on $26$ vertices. Since, by Lemma \ref{lemma:g26} only one of the extremal graphs on $26$ vertices has a non-empty set  $S_{\Gamma^{\prime},4}$, and this contains only one set, this produces a single graph.  Using a computer the graph was shown to satisfy the stated properties. See Appendix H for details. 

\section{Case $v=28$}
\vspace{.25cm}
\begin{theorem}
\label{theorem:g28}
If $\Gamma$ is a graph of girth at least $5$ on $28$ vertices with $68$ edges then  $(\delta,\Delta)=(3,6)$, $(deg_{3},deg_{4},deg_{5},deg_{6})=(1,4,21,2)$, and there is an embedded  $S_{6,[4,4,4,4,3,2]}$ star, or $(\delta,\Delta)=(4,5)$ or $(4,6)$. 

 If $(\delta,\Delta)=(4,5)$ then $(deg_{4},deg_{5})=(4,24)$ and there is an embedded $S_{5,[4,4,4,4,4]}$ star. If $(\delta,\Delta)=(4,6)$ then $(deg_{4},deg_{5},deg_{6})\in \{(i, 32-2i ,i-4):5\leq i \leq 13\}$ and there is an embedded $S_{6,[4,4,4,3,3,3]}$ star.
\end{theorem}

\vspace{.25cm}
\noindent
{\bf Proof} We know from Table  \ref{table:deltaTable} that in this case the possible values for $(\delta,\Delta)$ are $(3,5)$, $(3,6)$, $(3,7)$, $(3,8)$, $(3,9)$, $(4,5)$ or $(4,6)$. If  $\delta=3$,
by Lemma \ref{lem:inductiveExtremal}, $\Gamma$ can be constructed from the extremal graph on $27$ vertices described in Lemma \ref{lemma:g27}, $\Gamma^{\prime}$, by adding a vertex $x$ of degree $3$ and edges from $x$ to all vertices in an element of $S_{\Gamma^{\prime},3}$.
By Lemma \ref{lemma:g27} it follows that, in $\Gamma$, $(\delta,\Delta)=(3,6)$, $(deg_{3},deg_{4},deg_{5},deg_{6})=(1,4,21,2)$, and there is an embedded $S_{6,[4,4,4,4,3,2]}$ star. 

If $(\delta,\Delta)=(4,5)$ then $(deg_{4},deg_{5})=(4,24)$ and there are at most $16$ edges from the vertices of degree $4$ to the vertices of degree $5$. Hence at least one vertex of degree $5$ has no neighbours of degree $4$ and the result follows.

 If $(\delta,\Delta)=(4,6)$ then, since $deg_{4}>0$ and $deg_{6}>0$, $(deg_{4},deg_{5},deg_{6}) \in \{(i, 32-2i ,i-4):5\leq i \leq 16\}$. A vertex of degree $6$ can have at most one neighbour of degree $6$, and so has at least $5$ edges to the set $S$ containing the vertices of degrees $4$ and $5$. It follows that there is a linear space on $deg_{4}+deg_{5}$ points containing at least $deg_{6}$ blocks of size $5$. It follows from Lemma \ref{lemma:packingLemma} that this does not hold when $i\geq 14$.

 If there is no embedded $S_{6,[4,4,4,3,3,3]}$ star then every vertex of size $6$ has at least $4$ neighbours of degree $4$. This is only possible if there is a linear space on $deg_{4}$ points that has at least $deg_{6}$ blocks of size at least $4$.  By Lemma \ref{lemma:packingLemma}, this does not hold if $5<deg_{4}\leq 11$. So there is an embedded $S_{6,[4,4,4,3,3,3]}$ star except possibly where $deg_{4}=5$, $12$ or $13$.

Suppose that $deg_{4}=5$ then $(deg_{4},deg_{5},deg_{6})=(5,22,1)$.  If there is no embedded $S_{6,[4,4,4,3,3,3]}$ star the vertex $x$ of degree $6$ is
 adjacent to at least $4$ of the vertices of degree $4$. There are $3$ cases to consider: $x$ is adjacent to $5$ vertices of degree $4$ and $1$ of degree $5$;
 $x$ is adjacent to $4$ vertices of degree $4$ and there is a further edge among the vertices of degree $4$; and $x$ is adjacent to $4$ vertices of degree $4$
 and there is no further edge among the vertices of degree $4$. In all cases, let $T$ denote the set of edges from the vertices of degree $5$ to the vertices of degree $4$ and $6$. 
In the first two cases, $|T|=16$ and it follows that  the graph $\Gamma_{5}$ on the $22$ vertices of degree
 $5$ has $47$ edges and is extremal.

In the first case, each of the vertices of degree $4$ are in edges with $3$ of the vertices of degree $5$. Since all of the vertices of degree $4$ are adjacent to $x$, it follows that $S_{\Gamma_{5},3}$ contains $5$ non-intersecting sets.   This is not possible by Lemma \ref{lemma:g22A}.

In the second case, let $\Gamma_{5}$ be one of the extremal graphs $\Gamma_{0}$, $\Gamma_{1}$, $\Gamma_{2}$ described in Lemma \ref{lemma:g22}. 
Let $x$ be the vertex of degree $6$ and $\{a,b,c,d,e\}$ the vertices of degree $5$, where $x$ is adjacent to $a$,
 $b$, $c$ and $d$, and $d$ and $e$ are adjacent. Let $A$, $B$, $C$, $D$, $E$ and $X$ be the sets of vertices in
 $V_{5}$ that are adjacent to $a$, $b$, $c$, $d$, $e$ and $x$ in $\Gamma$. Note that each of these sets contain
 vertices that have degree at most $4$ in $\Gamma_{5}$, $|A|=|B|=|C|=|E|=3$, and $|D|=|X|=2$. Sets $|A|$, $|B|$, $|C|$, $|D|$
 and $|X|$ do not intersect, and any element that is in two of the sets is in $E$. It follows that $E$ contains all of the
 vertices of degree $3$ in $\Gamma_{5}$ and all of the vertices of degree $4$ in $\Gamma_{5}$ are in precisely one
 of the sets. Set $X$ must contain $2$ elements of degree $4$, $p$
 and $q$ say, that are not in edges with any element in sets $A$, $B$, $C$ or $D$ (or there is a square). So any
 neighbour of degree $4$ of $p$ or $q$ must be in set $E$, with the vertex of degree $3$ in $\Gamma_{5}$, and $p$ and $q$ have at most $2$ neighbours of degree $4$ between them. By
 Lemma \ref{lemma:g22A} this is not possible.

If $x$ is adjacent to $4$ vertices of degree $4$ and there is no
 further edge among the vertices of degree $4$, then let $\Sigma=\{a,b,c,d\}$ denote the set of vertices of degree $4$ adjacent to $x$, $e$ the remaining vertex of degree
 $4$, and $X$ the set of two vertices of degree $5$ that are adjacent to $x$. Let $\Gamma_{5,e}$ denote the graph formed by removing $x\cup \Sigma$ .  Then
 $\Gamma_{5,e}$ has $23$ vertices and $50$ edges,  is extremal and has vertices of degree $4$ or $5$. The vertices that have degree $4$ in $\Gamma_{5,e}$
 are those that are in an edge in $\Gamma$ with an element of $\Sigma$, vertices in $X$ and  $e$.  By an argument similar to the above, the elements of $X$ must have at most $1$ neighbour (i.e. $e$) of degree $4$ In $\Gamma_{5,e}$.  But Lemma \ref{lemma:g23A} the elements of $X$ have a common neighbour in
 $\Gamma_{5,e}$. Since, in $\Gamma$, both elements of $X$ are adjacent to $x$, there is a cycle of length $4$.  This is a contradiction, so if
 $(deg_{4},deg_{5},deg_{6})=(5,22,1)$ there is an embedded $S_{6,[4,4,4,3,3,3]}$ star.

Suppose that $(deg_{4},deg_{5},deg_{6})=(12, 8,8)$. If there is no embedded $S_{6,[4,4,4,3,3,3]}$ star every vertex of degree $6$ is adjacent to at least $4$ vertices of degree $4$. So there is a linear space on $12$ points containing $8$ blocks of size at least $4$. 
It follows from \cite{bettenbetten4} that there must be $12$ blocks of size exactly $4$. There is one such linear space, and all of the blocks intersect, by Lemma \ref{lemma:12.1}. It follows that no two vertices of degree $6$ can be adjacent to each other, or adjacent to a common vertex of degree $5$. Since every vertex is adjacent to two vertices of degree $5$ this is impossible. So if  $(deg_{4},deg_{5},deg_{6})=(12, 8,8)$ there is an embedded $S_{6,[4,4,4,3,3,3]}$ star.

Suppose that $(deg_{4},deg_{5},deg_{6})=(13, 6,9)$. If there is no embedded $S_{6,[4,4,4,3,3,3]}$ star every vertex of degree $6$ is adjacent to at least $4$ vertices of degree $4$. So there is a linear space $\Lambda$ on $13$ points containing $9$ blocks of size at least $4$. Since there is no linear space on $13$ points with a block of size $5$ and $8$ blocks of size at least $4$, there are exactly $36$ edges from the vertices of degree $6$ to the vertices of degree $4$. By Lemma \ref{lemma:13.1}, $\Lambda$ has no set of $3$ non-intersecting blocks so no vertex of degree $6$ has two neighbours of degree $6$, and so is in at least one edge with a vertex of degree $5$. Hence there are at least $9$ edges from the vertices of degree $6$ to the vertices of degree $5$. If a vertex of degree $5$ is in two edges with vertices of degree $6$, the corresponding blocks in $\Lambda$ must not intersect (or there is a square). By Lemma \ref{lemma:13.1}, $\Lambda$ has at most $3$ pairs of non-intersecting blocks and it follows that there are exactly $9$ edges from the vertices of degree $6$ to the vertices of degree $5$, and every vertex of degree $6$ is in one edge with another vertex of degree $6$. This is not possible (the sum of the degrees in the graph on the vertices of degree $6$ must be even).

\begin{lemma}\label{lemma:g28} There are $4$ graphs of girth at least $5$ on $28$ vertices with $68$ edges. 
\begin{itemize}
\item There is one graph for which  $(\delta,\Delta)=(3,6)$. Here $(deg_{3},deg_{4},deg_{5},deg_{6})=(1,4,21,2)$ and every element of 
$S_{\Gamma,4}$ contains the vertex of degree $3$, $1$ vertex of degree $4$ and $2$ of degree $5$.  Both vertices of degree $6$ are sink nodes, and the centre of embedded $S_{6,[4,4,4,4,3,2]}$ stars. 
\item There is one graph for which $(\delta,\Delta)=(4,5)$. Set $S_{\Gamma,4}$ is empty.
\item There are $2$ graphs for which $(\delta,\Delta)=(4,6)$. In each case there are embedded $S_{6,[444,333]}$ star. In the first graph  $(deg_{4},deg_{5},deg_{6})=(6,20,2)$ and every element of $S_{\Gamma,4}$ contains $2$ vertices of degree $4$ and $2$ of degree $5$, and in the second, $(deg_{4},deg_{5},deg_{6})=(7,18,3)$ and every element of $S_{\Gamma,4}$ contains $3$ vertices of degree $4$ and $1$ of degree $5$.  In both cases, for every sink node $p$, every element of $S_{\Gamma,4}$ contains exactly $1$ vertex of degree $4$ adjacent to $p$.  
\end{itemize}
\end{lemma}

\vspace{.25cm}
\noindent
{\bf Proof} (1) From Theorem \ref{theorem:g28}  all extremal graphs on $28$ vertices with $\delta=3$ can be constructed by adding a new vertex $x$, and edges to a set of $3$ vertices that is a set in $S_{\Gamma^{\prime},3}$, where $\Gamma^{\prime}$ is the extremal graph on $27$ vertices. Since, by Lemma \ref{lemma:g27} there is one graph $\Gamma^{\prime}$ and $5$ sets $S_{\Gamma^{\prime},3}$, there are $5$ graphs to check. Using nauty it was shown that all of these graphs are isomorphic, so there is only one graph in this case.   Using a computer the graph was shown to satisfy the stated property. 

\vspace{.25cm}
\noindent
(2) If  $(\delta,\Delta)=(4,5)$ then $(deg_{4},deg_{5})=(4,24)$. If there are $e$ edges amongst the vertices of degree $4$ then there are $16-2e$ between vertices of degree $4$ and vertices of degree $5$, and $52+e$ among the vertices of degree $5$. It follows that $e\leq 2$.  

If $e=2$ then the graph on the vertices of degree $5$ ($\Gamma^{\prime}$) has $54$ edges and is extremal. If $x$ is a vertex in $\Gamma$ of degree $4$ that is in no edges with other vertices of degree $4$, then the graph on the vertices of degree $5$ and $x$ has $25$ vertices and $58$ edges, which is not possible. So there are two parallel edges on the vertices of degree $4$. The sets of vertices of degree $5$ adjacent to each of the vertices of degree $4$ must not intersect (as $\Gamma^{\prime}$ contains no vertex of degree less than $4$) and for each pair of adjacent vertices of degree $4$ there is no edge between their neighbours of degree $5$ (or there is a square). It follows that $S_{\Gamma^{\prime},3}$ contains $4$ non-intersecting elements  for which there are are two pairs with no edge between either of the pairs. This is not possible by Lemma \ref{lemma:g24}.

\vspace{.25cm}
\noindent
If $e=0$ then restricting search to this case showed there to be no graphs. The search in this case took $13$ hours.  

\vspace{.25cm}
\noindent
If $e=1$ then suppose that $v_{1}$ and $v_{2}$ are adjacent. The graph $\Gamma^{\prime}$ with vertices $V\setminus\{v_{1},v_{2}\}$ has $26$ vertices and $61$ edges, so is maximal. The sets of three vertices of degree $5$ adjacent to $v_{1}$ and $v_{2}$ respectively must be distinct and there must be no edge between these sets. So $S_{\Gamma^{\prime},3}$ contains a pair of non-intersecting sets that have no edges between them. By Lemma \ref{lemma:g26} there is one maximal graph $\Gamma^{\prime}$ for which $S_{\Gamma^{\prime},3}$ contains such a pair, and one suitable pair (($X_{1},X_{2})$ say). Hence, since $\Gamma$ is constructed from $\Gamma^{\prime}$ by adding a new pair of adjacent vertices $v_{1}$ and $v_{2}$ with edges to  the vertices in $X_{1}$ and $X_{2}$, there is one graph $\Gamma$.  

\vspace{.25cm}
\noindent
(3) Let $\Gamma$ be a graph with $(\delta,\Delta)=(4,6)$.  By the proof of Theorem \ref{theorem:g28} $(deg_{4},deg_{5},deg_{6}) \in \{(i, 32-2i ,i-4):5\leq i \leq 16\}$. If there is an edge between two vertices, $v_{1}$ and $v_{2}$ of degree $4$, then removing the vertices  and all edges on them leaves an extremal graph on $26$ vertices. Since neither extremal graph on $26$ vertices has any vertices of degree $6$, it follows that the sets of vertices adjacent to $v_{1}$ (that are not $v_{2}$) and to $v_{2}$ (that are not $v_{1}$) are elements of $S_{\Gamma^{\prime},3}$ that are parallel, have no edges between them, and contain at least one vertex of degree $5$. Only one of the extremal graphs on $26$ vertices have such a pair of sets, and there are $3$ such pairs. These lead to two non-isomorphic extremal graphs $\Gamma$ which  have $(deg_{4},deg_{5},deg_{6})=(6,20,2)$ and $(deg_{4},deg_{5},deg_{6})=(7,18,3)$ respectively. 

Let us assume then that there are no edges amongst the vertices of degree $4$. If a vertex $v_{1}$ of degree $5$ has $2$ neighbours, $q_{1}$, $q_{2}$ say, of degree $4$, then removing $v_{1}$, $q_{1}$ and $q_{2}$ leaves an extremal graph on $25$ vertices (with $57$ edges), $\Gamma^{\prime}$ say. Let $X_{1}$, $X_{2}$ and $X_{3}$ be the remaining neighbours of $v_{1}$, $q_{1}$ and $q_{2}$ in $\Gamma$. Then $X_{1}$, $X_{2}$ and $X_{3}$ are parallel elements of $S_{\Gamma^{\prime},3}$ for which there are no edges between $X_{1}$ and $X_{2}$, or between $X_{1}$ and $X_{3}$. In addition, since $\Gamma$ has no vertices of degree $3$, any vertex of degree $3$ in $\Gamma^{\prime}$ must be in $X_{1}$, $X_{2}$ or $X_{3}$. Since there are no edges between vertices of degree $4$ in $\Gamma$, if follows that all of the vertices of degree $3$ in $\Gamma^{\prime}$ are in $X_{1}$. There is no suitable extremal graph $\Gamma^{\prime}$, by Lemma \ref{lemma:g25A}. So in $\Gamma$ no vertex of degree $5$ has more than $1$ neighbour of degree $4$.

Consider the embedded $S_{6,[4,4,4,3,3,3]}$ star, with root $p$. Three of its children have degree $4$ and all other vertices of degree $4$ must be grandchildren, with parents having degree $5$ (as there are no edges between vertices of degree $4$). No two of these grandchildren can have the same parent (as no vertex of degree $5$ has more than one child of degree $4$). It follows that $i\leq 6$.

If $i=5$ then all of the grandchildren have degree $4$ or $5$. Any grandchild $x$ of degree $4$ must have $3$ edges to vertices of degree $5$ which don't have a parent of degree $4$. It follows $x$ has edges to two vertices with the same parent, and there is a square. 

\begin{figure}
\centering
\includegraphics[width=0.5\textwidth]{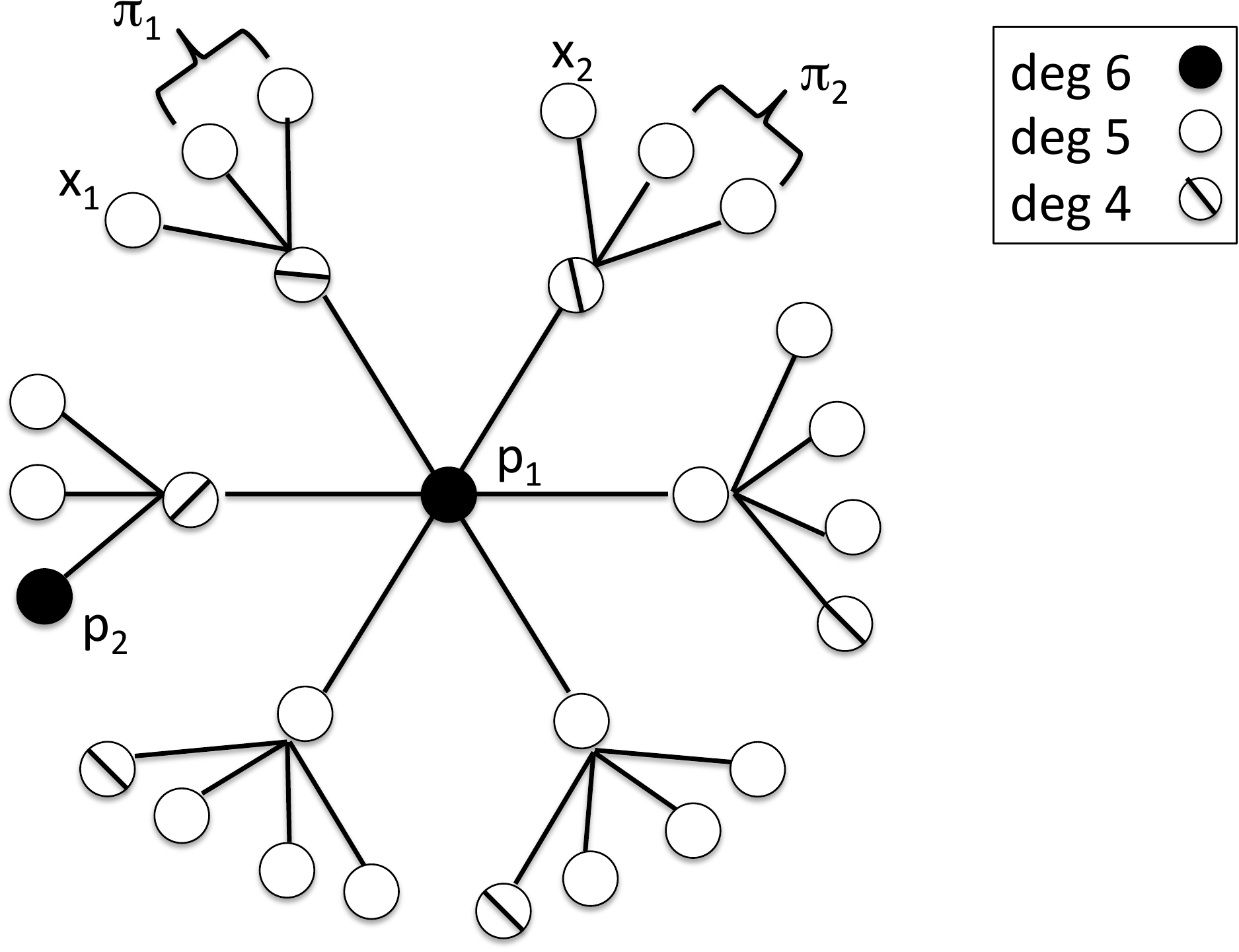}
\caption{Embedded star, $i=6$, no edges between vertices of degree $4$
\label{fig:28star}}
\end{figure}

If $i=6$ a similar argument shows that all of the grandchildren of degree $4$ must have edges to the grandchild $p_{2}$ of degree $6$, and $p_{2}$ has a parent of degree $4$ (or a similar argument applies).      So we have the situation illustrated in Figure \ref{fig:28star}. Note that $p_{2}$ is adjacent to all of the grandchildren of degree $4$ and two of the grandchildren of degree $5$, $x_{1}$ and $x_{2}$. Now the grandchildren of degree $4$ are each adjacent to $p_{2}$ and two of the grandchildren of degree $5$ that have a parent of degree $5$. Thus the vertices in $\pi_{1}$ and $\pi_{2}$ are at a distance greater than $2$ from $p_{2}$ unless they are adjacent to $x_{1}$ or $x_{2}$. At most two elements of $\pi_{1}\cup \pi_{2}$ can have edges to $x_{1}$ or $x_{2}$, so at least two vertices are at a distance greater than $2$ from $p_{2}$, which is impossible as $p_{2}$ is the root of an embedded $S_{6,[4,4,4,3,3,3]}$ star (which contains $27$ of the vertices).

\section{Case $v=29$}
\vspace{.25cm}
\begin{theorem}
\label{theorem:g29}
If $\Gamma$ is a graph of girth at least $5$ on $29$ vertices with $72$ edges then  $(\delta,\Delta)=(4,6)$ and $\Gamma$ has an embedded  $S_{6,[4,4,4,4,3,3]}$ star. 
\end{theorem}

\vspace{.25cm}
\noindent
{\bf Proof} We know from Table  \ref{table:deltaTable} that in this case the possible values for $(\delta,\Delta)$ are $(4,5)$, $(4,6)$ and $(4,7)$. Since in all cases $\delta= f(29)-f(28)$, $\Gamma$ is constructed from an extremal graph on $28$ vertices, $\Gamma^{\prime}$,  by adding a new vertex $x$ of degree $4$ and $4$ edges from $x$ to a set $S_{\Gamma^{\prime},4}$.  By Lemma \ref{lemma:g28} there are $4$ possibilities for $\Gamma^{\prime}$. (i) $\Gamma^{\prime}$ has $(deg_{3},deg_{4},deg_{5},deg_{6})=(1,4,21,2)$. In this case, all elements of $S_{\Gamma^{\prime},4}$ contain the vertex of degree $3$, which is a child of the root of a $S_{6,[4,4,4,4,3,2]}$ star, and no vertex of degree $6$ -- the result follows. (ii) $\Gamma^{\prime}$ has $(deg_{4},deg_{5})=(4,24)$. In this case $S_{\Gamma^{\prime},4}$ is empty, so this case is impossible. (iii) $\Gamma^{\prime}$ has  $(deg_{4},deg_{5},deg_{6})=(6,20,2)$ or $(7, 18,3)$ then all elements of $S_{\Gamma^{\prime},4}$ contain a child of degree $4$ of a root of a $S_{6,[4,4,4,3,3,3]}$ star and no vertex of degree $6$, and the result follows.  

\begin{lemma}\label{lemma:g29} There is  $1$ graph of girth at least $5$ on $29$ vertices with $72$ edges. In this graph $(\delta,\Delta)=(4,6)$ and $(deg_{4},deg_{5},deg_{6})=(5,20,4)$. Set $S_{\Gamma,4}$ contains $2$ sets, each of which contains $2$ vertices of degree $4$ and $2$ of degree $5$.  
\end{lemma}

\vspace{.25cm}
\noindent
{\bf Proof} From Theorem \ref{theorem:g29}  any extremal graph $\Gamma$ on $29$ vertices can be constructed by adding a new vertex $x$, and edges to a set of $4$ vertices that is a set in $S_{\Gamma^{\prime},4}$, where $\Gamma^{\prime}$ is one of the $4$ extremal graphs on $28$ vertices. There are $7$ potential graphs $\Gamma$.  All of these graphs are shown to be isomorphic, using nauty, and the single graph was shown to satisfy the stated property. 

\section{Case $v=30$}
\vspace{.25cm}
\begin{theorem}
\label{theorem:g30}
If $\Gamma$ is a graph of girth at least $5$ on $30$ vertices with $76$ edges then  $(\delta,\Delta)=(4,6)$ and $\Gamma$ has an embedded  $S_{6,[4,4,4,4,4,3]}$ star and an embedded $S_{6,[5,4,4,4,3,3]}$ star.
\end{theorem}

\vspace{.25cm}
\noindent
{\bf Proof} We know from Table  \ref{table:deltaTable} that in this case the possible values for $(\delta,\Delta)$ are $(4,6)$ and $(4,7)$. Since in all cases $\delta= f(30)-f(29)$, $\Gamma$ is constructed from the extremal graph on $29$ vertices, $\Gamma^{\prime}$, by adding a new vertex $x$ of degree $4$ and $4$ edges from $x$ to a set in $S_{\Gamma^{\prime},4}$.  By Lemma \ref{lemma:g29}, $\Gamma^{\prime}$ has $(deg_{4},deg_{5},deg_{6})=(5,20,4)$. In this case, all elements of $S_{\Gamma^{\prime},4}$ contain two vertices of degree $4$ and two of degree $5$. So clearly, in $\Gamma$, $(\delta,\Delta)=(4,6)$ and $(deg_{4},deg_{5},deg_{6})=(4,20,6)$ . Both elements of $S_{\Gamma^{\prime},4}$ contain the child of an embedded $S_{6,[4,4,4,4,3,3]}$ star of $\Gamma^{\prime}$ of degree $4$, and the child of another embedded $S_{6,[4,4,4,4,3,3]}$ star of $\Gamma^{\prime}$ of degree $5$. The result follows. 

\begin{lemma}\label{lemma:g30} There is  $1$ graph of girth at least $5$ on $30$ vertices with $76$ edges. In this graph $(\delta,\Delta)=(4,6)$ and $(deg_{4},deg_{5},deg_{6})=(4,20,6)$. Set $S_{\Gamma,4}$ contains $2$ sets, each of which contains $2$ vertices of degree $4$ and $2$ of degree $5$, and $S_{\Gamma,5}$ is empty.   
\end{lemma}

\vspace{.25cm}
\noindent
{\bf Proof} From Theorem \ref{theorem:g30}  any extremal graph $\Gamma$ on $30$ vertices can be constructed by adding a new vertex $x$, and edges to a set of $4$ vertices that is a set in $S_{\Gamma^{\prime},4}$, where $\Gamma^{\prime}$ is the extremal graphs on $29$ vertices. There are $2$ potential graphs $\Gamma$, which are shown to be isomorphic, using nauty, and the single graph was shown to satisfy the stated properties. 

\section{Case $v=31$}
\vspace{.25cm}
\begin{theorem}
\label{theorem:g31}
If $\Gamma$ is a graph of girth at least $5$ on $31$ vertices then $\Gamma$ has at most $80$ edges. If $\Gamma$ has $80$ edges then  either $(\delta,\Delta)=(4,6)$ and $\Gamma$ has an embedded $S_{6,[4,4,4,4,4,4]}$ star and an embedded $S_{6,[5,4,4,4,4,3]}$ star, or $(\delta,\Delta)=(5,6)$ and $\Gamma$ has an embedded $S_{6,[4,4,4,4,4,4]}$ star.
\end{theorem}

\vspace{.25cm}
\noindent
{\bf Proof} Suppose that $\Gamma$ had $81$ edges, then we must have $\delta=5$. Since the extremal graph $\Gamma^{\prime}$ on $29$ vertices with $76$ edges has empty set $S_{\Gamma^{\prime},5}$, this is not possible. So there are at most $80$ edges. We know from Table  \ref{table:deltaTable} that in this case the possible values for $(\delta,\Delta)$ are $(4,6)$, $(4,7)$ and $(5,6)$. 

If $\delta=4$ then since $\delta= f(31)-f(30)$, $\Gamma$ is constructed from the extremal graph on $30$ vertices, $\Gamma^{\prime}$, by adding a new vertex $x$ of degree $4$ and $4$ edges from $x$ to a set from $S_{\Gamma^{\prime},4}$.  By Lemma \ref{lemma:g30}  $\Gamma^{\prime}$ has $(deg_{4},deg_{5},deg_{6})=(4,20,6)$. 
All elements of $S_{\Gamma^{\prime},4}$ contain two vertices of degree $4$ and two of degree $5$. So clearly, in $\Gamma$, $(\delta,\Delta)=(4,6)$  and $(deg_{4},deg_{5},deg_{6})=(3,20,8)$. 
By Appendix $K$, $\Gamma^{\prime}$ has two  embedded $S_{6,[4,4,4,4,4,3]}$  stars, each of which has children from both elements, $s_{1}$ and $s_{2}$ say, of $S_{\Gamma^{\prime},4}$. In the first star the child from $s_{1}$ has degree $4$ and the child from $s_{2}$ has degree $5$, and in the second star, the reverse is true.  The result follows.

If $(\delta,\Delta)=(5,6)$ then $(deg_{5},deg_{6})=(26,5)$. Every vertex of degree $6$ has $6$ neighbours of degree $5$. 

\begin{lemma}\label{lemma:g31} There are $2$  graphs of girth at least $5$ on $31$ vertices with $80$ edges. In the first  graph $(\delta,\Delta)=(4,6)$, $(deg_{4},deg_{5},deg_{6})=(3,20,8)$ and set $S_{\Gamma,5}$ contains one element, which consists of $2$ vertices of degree $5$ and $3$ of degree $4$.  In the second, $(\delta,\Delta)=(5,6)$ $(deg_{5},deg_{6})=(26,5)$ and  $S_{\Gamma,5}$ contains one element, which consists of $5$ vertices of degree $5$.
\end{lemma}

\vspace{.25cm}
\noindent
{\bf Proof} From Theorem \ref{theorem:g31}  any extremal graph $\Gamma$ on $31$ vertices can either be constructed by adding a new vertex $x$, and edges to a set of $4$ vertices that is a set in $S_{\Gamma^{\prime},4}$, where $\Gamma^{\prime}$ is the extremal graphs on $30$ vertices, or has $(\delta,\Delta)=(5,6)$ and an embedded $S_{6,[4,4,4,4,4,4]}$ star. In the first case there are $2$ potential graphs $\Gamma$, which are shown to be isomorphic, using nauty, and the single graph was shown to satisfy the stated property. When $(\delta,\Delta)=(5,6)$ restricting search to this case, and fixing the embedded $S_{6,[4,4,4,4,4,4]}$   produced $720$ graphs. Applying nauty revealed there to be $1$ unique graph, which satisfies the stated property.

\section{Case $v=32$}

\begin{theorem}
\label{theorem:g32}
If $\Gamma$ is a graph of girth at least $5$ on $32$ vertices with $85$ edges then there is an embedded $S_{6,[5,4,4,4,4,4]}$ star. 
\end{theorem}

 \vspace{.25cm}
\noindent{\bf Proof} We know from Table  \ref{table:deltaTable} that in this case all vertices have degree $5$ or $6$. There are $22$ vertices of degree $5$ and $10$ vertices of degree $6$. 

Suppose that no vertex of degree $6$ has a neighbour of degree $6$. Then every vertex of degree $6$ has $6$ neighbours of degree $5$. The sets of neighbours of each of the vertices of degree $6$ form a linear space on $22$ points with $10$ blocks of size $6$, contradicting Lemma \ref{lemma:packingLemma}. 

\begin{lemma}\label{lemma:g32} There is  $1$ graph of girth at least $5$ on $32$ vertices with $85$ edges. In this graph $(\delta,\Delta)=(5,6)$, $(deg_{5},deg_{6})=(22,10)$ and $S_{\Gamma,3}$ is empty. Every element of $S_{\Gamma,2}$ contains $2$ vertices of degree $5$, at least one of which is the child of an embedded $S_{6,[5,4,4,4,4,4]}$ star. 
\end{lemma}

\vspace{.25cm}
\noindent
{\bf Proof} From Theorem \ref{theorem:g32}  all extremal graphs on $32$ vertices has $(\delta,\Delta)=(5,6)$ and an embedded $S_{6,[5,4,4,4,4,4]}$ star. Restricting search to this case and applying nauty revealed there to be $1$ graph. Alternatively, since $\delta=f(v)-f(v-1)$ we can construct all graphs by  adding a new vertex $x$ and edges to a set of $5$ vertices that is a set in $S_{\Gamma^{\prime},5}$, where $\Gamma^{\prime}$ is an extremal graph on $30$ vertices. There are two such graphs, and in each case $S_{\Gamma^{\prime},5}$ contains $1$ element. So there are $2$ potential graphs. These are shown to be isomorphic.

\section*{Appendix A}\label{appendix20}
The extremal graph $(20,41)$. Note that $\Gamma_{0}$ has $(deg_{3},deg_{4},deg_{5})=(1,16,3)$. 

\vspace{.25cm}
\noindent
$\Gamma_{0}$ has incidence array:

\vspace{.25cm}
\noindent
\begin{tiny}
$
\begin{array}{ccccc ccccc ccccc ccccc}
0&0&0&0&1&0&0&1&0&0&1&0&0&1&0&0&0&0&0&1\\
0&0&0&0&0&1&0&0&1&0&0&1&0&0&1&0&0&0&0&1\\
0&0&0&1&0&0&1&0&0&1&0&0&1&0&0&0&0&0&0&1\\
0&0&1&0&1&0&0&0&0&0&0&0&0&0&1&1&0&0&0&0\\
1&0&0&1&0&1&0&0&0&0&0&0&0&0&0&0&0&0&1&0\\
0&1&0&0&1&0&1&0&0&0&0&0&0&0&0&0&0&1&0&0\\
0&0&1&0&0&1&0&1&0&0&0&0&0&0&0&0&1&0&0&0\\
1&0&0&0&0&0&1&0&1&0&0&0&0&0&0&1&0&0&0&0\\
0&1&0&0&0&0&0&1&0&1&0&0&0&0&0&0&0&0&1&0\\
0&0&1&0&0&0&0&0&1&0&1&0&0&0&0&0&0&1&0&0\\
1&0&0&0&0&0&0&0&0&1&0&1&0&0&0&0&1&0&0&0\\
0&1&0&0&0&0&0&0&0&0&1&0&1&0&0&1&0&0&0&0\\
0&0&1&0&0&0&0&0&0&0&0&1&0&1&0&0&0&0&1&0\\
1&0&0&0&0&0&0&0&0&0&0&0&1&0&1&0&0&1&0&0\\
0&1&0&1&0&0&0&0&0&0&0&0&0&1&0&0&1&0&0&0\\
0&0&0&1&0&0&0&1&0&0&0&1&0&0&0&0&0&1&0&0\\
0&0&0&0&0&0&1&0&0&0&1&0&0&0&1&0&0&0&1&0\\
0&0&0&0&0&1&0&0&0&1&0&0&0&1&0&1&0&0&0&0\\
0&0&0&0&1&0&0&0&1&0&0&0&1&0&0&0&1&0&0&0\\
1&1&1&0&0&0&0&0&0&0&0&0&0&0&0&0&0&0&0&0\\
\end{array}
$
\end{tiny}

\vspace{.25cm}
\noindent
and degree sequence: $$(5,5,5,4,4,4,4,4,4,4,4,4,4,4,4,4,4,4,4,3)$$
Edge set is:
\begin{eqnarray*}
(0,4), (0,7), (0,10), (0,13), (0,19), (1,5), (1,8),
 (1,11),\\
 (1,14), (1,19), (2,3), (2,6), (2,9), (2,12), (2,19),\\ (3,4),
 (3,14), (3,15), (4,5), (4,18), (5,6), (5,17),\\ (6,7),
 (6,16), (7,8), (7,15), (8,9), (8,18), (9,10), \\
(9,17), (10,11), (10,16), (11,12), (11,15), (12,13), \\
(12,18), (13,14), (13,17), (14,16), (15,17), (16,18)
\end{eqnarray*}

\vspace{.25cm}
\noindent
Set $S_{\Gamma_{0},3}$  is given below. All elements contain vertex $19$ of degree $3$, and two other vertices of degree $4$. 
\begin{eqnarray*}
S_{\Gamma_{0},3} &=&\{(15,16,19) (15,18,19) (16,17,19) (17,18,19)\}\\
\end{eqnarray*}


\section*{Appendix B}\label{appendix21}
The $3$ extremal graphs $(21,44)$. Note that graph $\Gamma_{0}$ has $(deg_{3},deg_{4},deg_{5})=(1,15,5)$, and graphs $\Gamma_{1}$ and $\Gamma_{2}$ have $(deg_{4},deg_{5})=(17,4)$. 

\vspace{.25cm}
\noindent
$\Gamma_{0}$ has incidence array:

\vspace{.25cm}
\noindent
\begin{tiny}
$
\begin{array}{ccccc ccccc ccccc ccccc  c}
0&0&0&0&1&0&0&1&0&0&1&0&0&1&0&0&0&0&0&1&0\\
0&0&0&0&0&1&0&0&1&0&0&1&0&0&1&0&0&0&0&1&0\\
0&0&0&1&0&0&1&0&0&1&0&0&1&0&0&0&0&0&0&1&0\\
0&0&1&0&1&0&0&0&0&0&0&0&0&0&1&1&0&0&0&0&0\\
1&0&0&1&0&1&0&0&0&0&0&0&0&0&0&0&0&0&1&0&0\\
0&1&0&0&1&0&1&0&0&0&0&0&0&0&0&0&0&1&0&0&0\\
0&0&1&0&0&1&0&1&0&0&0&0&0&0&0&0&1&0&0&0&0\\
1&0&0&0&0&0&1&0&1&0&0&0&0&0&0&1&0&0&0&0&0\\
0&1&0&0&0&0&0&1&0&1&0&0&0&0&0&0&0&0&1&0&0\\
0&0&1&0&0&0&0&0&1&0&1&0&0&0&0&0&0&1&0&0&0\\
1&0&0&0&0&0&0&0&0&1&0&1&0&0&0&0&1&0&0&0&0\\
0&1&0&0&0&0&0&0&0&0&1&0&1&0&0&1&0&0&0&0&0\\
0&0&1&0&0&0&0&0&0&0&0&1&0&1&0&0&0&0&1&0&0\\
1&0&0&0&0&0&0&0&0&0&0&0&1&0&1&0&0&1&0&0&0\\
0&1&0&1&0&0&0&0&0&0&0&0&0&1&0&0&1&0&0&0&0\\
0&0&0&1&0&0&0&1&0&0&0&1&0&0&0&0&0&1&0&0&1\\
0&0&0&0&0&0&1&0&0&0&1&0&0&0&1&0&0&0&1&0&1\\
0&0&0&0&0&1&0&0&0&1&0&0&0&1&0&1&0&0&0&0&0\\
0&0&0&0&1&0&0&0&1&0&0&0&1&0&0&0&1&0&0&0&0\\
1&1&1&0&0&0&0&0&0&0&0&0&0&0&0&0&0&0&0&0&1\\
0&0&0&0&0&0&0&0&0&0&0&0&0&0&0&1&1&0&0&1&0\\
\end{array}
$
\end{tiny}

\vspace{.25cm}
\noindent
and degree sequence: $$(5,5,5,4,4,4,4,4,4,4,4,4,4,4,4,5,5,4,4,4,3)$$ Edge set is:
\begin{eqnarray*}
(0,4), (0,7), (0,10), (0,13), (0,19), (1,5), (1,8), (1,11),\\ (1,14), 
(1,19), (2,3), (2,6), (2,9), (2,12), (2,19), (3,4),\\ (3,14), (3,15), (4,5), 
(4,18), (5,6), (5,17), (6,7),\\ (6,16), (7,8), (7,15), (8,9), (8,18), (9,10), 
(9,17),\\ (10,11), (10,16), (11,12), (11,15), (12,13), (12,18), (13,14),\\ (13,17), (14,16), 
(15,17), (15,20), (16,18), (16,20), (19,20)
\end{eqnarray*}

\vspace{.25cm}
\noindent
Set $S_{\Gamma_{0},3}$ is given below. All elements contain vertices of degree $3$ or $4$.
\begin{eqnarray*}
S_{\Gamma_{0},3} &=&\{(4,9,20) (5,12,20) (8,13,20) (17,18,19) \}\\
\end{eqnarray*}

\vspace{.25cm}
\noindent
$\Gamma_{1}$ has incidence array:

\vspace{.25cm}
\noindent
\begin{tiny}
$
\begin{array}{ccccc ccccc ccccc ccccc  c}
0&1&1&1&1&1&0&0&0&0&0&0&0&0&0&0&0&0&0&0&0\\
1&0&0&0&0&0&1&1&1&0&0&0&0&0&0&0&0&0&0&0&0\\
1&0&0&0&0&0&0&0&0&1&1&1&0&0&0&0&0&0&0&0&0\\
1&0&0&0&0&0&0&0&0&0&0&0&1&1&1&0&0&0&0&0&0\\
1&0&0&0&0&0&0&0&0&0&0&0&0&0&0&1&1&1&0&0&0\\
1&0&0&0&0&0&0&0&0&0&0&0&0&0&0&0&0&0&1&1&1\\
0&1&0&0&0&0&0&0&0&1&0&0&1&0&0&1&0&0&1&0&0\\
0&1&0&0&0&0&0&0&0&0&1&0&0&1&0&0&1&0&0&1&0\\
0&1&0&0&0&0&0&0&0&0&0&1&0&0&1&0&0&1&0&0&1\\
0&0&1&0&0&0&1&0&0&0&0&0&0&1&0&0&0&1&0&0&0\\
0&0&1&0&0&0&0&1&0&0&0&0&0&0&1&1&0&0&0&0&0\\
0&0&1&0&0&0&0&0&1&0&0&0&1&0&0&0&1&0&0&0&0\\
0&0&0&1&0&0&1&0&0&0&0&1&0&0&0&0&0&0&0&1&0\\
0&0&0&1&0&0&0&1&0&1&0&0&0&0&0&0&0&0&0&0&1\\
0&0&0&1&0&0&0&0&1&0&1&0&0&0&0&0&0&0&1&0&0\\
0&0&0&0&1&0&1&0&0&0&1&0&0&0&0&0&0&0&0&0&1\\
0&0&0&0&1&0&0&1&0&0&0&1&0&0&0&0&0&0&1&0&0\\
0&0&0&0&1&0&0&0&1&1&0&0&0&0&0&0&0&0&0&1&0\\
0&0&0&0&0&1&1&0&0&0&0&0&0&0&1&0&1&0&0&0&0\\
0&0&0&0&0&1&0&1&0&0&0&0&1&0&0&0&0&1&0&0&0\\
0&0&0&0&0&1&0&0&1&0&0&0&0&1&0&1&0&0&0&0&0\\
\end{array}
$
\end{tiny}

\vspace{.25cm}
\noindent
and degree sequence: $$(5,4,4,4,4,4,5,5,5,4,4,4,4,4,4,4,4,4,4,4,4)$$ Edge set is:
\begin{eqnarray*}
(0,1), (0,2), (0,3), (0,4), (0,5), (1,6), (1,7), (1,8),\\ (2,9), 
(2,10), (2,11), (3,12), (3,13), (3,14), (4,15), (4,16), \\(4,17), (5,18), (5,19), 
(5,20), (6,9), (6,12), (6,15),\\ (6,18), (7,10), (7,13), (7,16), (7,19), (8,11), 
(8,14),\\ (8,17), (8,20), (9,13), (9,17), (10,14), (10,15), (11,12),\\ (11,16), (12,19), 
(13,20), (14,18), (15,20), (16,18), (17,19)
\end{eqnarray*}

\vspace{.25cm}
\noindent
Set $S_{\Gamma_{1},3}$ is empty.

\vspace{.25cm}
\noindent
$\Gamma_{2}$ has incidence array:

\vspace{.25cm}
\noindent
\begin{tiny}
$
\begin{array}{ccccc ccccc ccccc ccccc ccccc}
0&1&1&1&1&1&0&0&0&0&0&0&0&0&0&0&0&0&0&0&0\\
1&0&0&0&0&0&1&1&1&0&0&0&0&0&0&0&0&0&0&0&0\\
1&0&0&0&0&0&0&0&0&1&1&1&0&0&0&0&0&0&0&0&0\\
1&0&0&0&0&0&0&0&0&0&0&0&1&1&1&0&0&0&0&0&0\\
1&0&0&0&0&0&0&0&0&0&0&0&0&0&0&1&1&1&0&0&0\\
1&0&0&0&0&0&0&0&0&0&0&0&0&0&0&0&0&0&1&1&1\\
0&1&0&0&0&0&0&0&0&1&0&0&1&0&0&1&0&0&1&0&0\\
0&1&0&0&0&0&0&0&0&0&1&0&0&1&0&0&1&0&0&1&0\\
0&1&0&0&0&0&0&0&0&0&0&1&0&0&1&0&0&1&0&0&0\\
0&0&1&0&0&0&1&0&0&0&0&0&0&1&0&0&0&1&0&0&0\\
0&0&1&0&0&0&0&1&0&0&0&0&0&0&1&1&0&0&0&0&0\\
0&0&1&0&0&0&0&0&1&0&0&0&1&0&0&0&1&0&0&0&1\\
0&0&0&1&0&0&1&0&0&0&0&1&0&0&0&0&0&0&0&1&0\\
0&0&0&1&0&0&0&1&0&1&0&0&0&0&0&0&0&0&0&0&1\\
0&0&0&1&0&0&0&0&1&0&1&0&0&0&0&0&0&0&1&0&0\\
0&0&0&0&1&0&1&0&0&0&1&0&0&0&0&0&0&0&0&0&1\\
0&0&0&0&1&0&0&1&0&0&0&1&0&0&0&0&0&0&1&0&0\\
0&0&0&0&1&0&0&0&1&1&0&0&0&0&0&0&0&0&0&1&0\\
0&0&0&0&0&1&1&0&0&0&0&0&0&0&1&0&1&0&0&0&0\\
0&0&0&0&0&1&0&1&0&0&0&0&1&0&0&0&0&1&0&0&0\\
0&0&0&0&0&1&0&0&0&0&0&1&0&1&0&1&0&0&0&0&0\\
\end{array}
$
\end{tiny}

\vspace{.25cm}
\noindent
and degree sequence: $$(5,4,4,4,4,4,5,5,4,4,4,5,4,4,4,4,4,4,4,4,4)$$ Edge set is:
\begin{eqnarray*}
(0,1), (0,2), (0,3), (0,4), (0,5), (1,6), (1,7), (1,8), \\(2,9), 
(2,10), (2,11), (3,12), (3,13), (3,14), (4,15), (4,16), \\(4,17), (5,18), (5,19), 
(5,20), (6,9), (6,12), (6,15),\\ (6,18), (7,10), (7,13), (7,16), (7,19), (8,11), 
(8,14),\\ (8,17), (9,13),(9,17), (10,14), (10,15), (11,12), (11,16),\\ (11,20), (12,19), 
(13,20), (14,18), (15,20), (16,18), (17,19)
\end{eqnarray*}

\vspace{.25cm}
\noindent
Sets $S_{\Gamma,3}$ is empty.

\section*{Appendix C}\label{appendix22}
The $3$ extremal graphs $(22,47)$. Note that graph $\Gamma_{0}$ has $(deg_{4},deg_{5})=(16,6)$, graph $\Gamma_{1}$ has $(deg_{3},deg_{4},deg_{5})=(2,12,8)$ and graph $\Gamma_{2}$ has  $(deg_{3},deg_{4},deg_{5})=(1,14,7)$. 

\vspace{.25cm}
\noindent
$\Gamma_{0}$ has incidence array:

\vspace{.25cm}
\noindent
\begin{tiny}
$
\begin{array}{cccccccccccccccccccccccc}
0&0&0&0&0&0&0&0&0&1&0&0&1&0&0&1&0&0&0&0&1&0\\
0&0&0&0&0&0&1&0&1&0&0&1&0&0&0&0&0&0&0&0&1&0\\
0&0&0&0&0&1&0&1&0&0&0&0&0&0&1&0&0&0&0&0&1&0\\
0&0&0&0&1&0&0&0&0&0&1&0&0&1&0&0&0&0&0&0&1&0\\
0&0&0&1&0&0&0&1&0&0&0&1&0&0&0&1&0&0&0&1&0&0\\
0&0&1&0&0&0&0&0&1&0&0&0&0&1&0&0&0&0&0&1&0&0\\
0&1&0&0&0&0&0&0&0&0&0&0&1&0&1&0&0&0&0&1&0&0\\
0&0&1&0&1&0&0&0&0&0&0&0&1&0&0&0&0&0&1&0&0&0\\
0&1&0&0&0&1&0&0&0&0&1&0&0&0&0&1&0&0&1&0&0&0\\
1&0&0&0&0&0&0&0&0&0&0&1&0&0&1&0&0&0&1&0&0&0\\
0&0&0&1&0&0&0&0&1&0&0&0&0&0&1&0&0&1&0&0&0&0\\
0&1&0&0&1&0&0&0&0&1&0&0&0&0&0&0&0&1&0&0&0&0\\
1&0&0&0&0&0&1&1&0&0&0&0&0&1&0&0&0&1&0&0&0&0\\
0&0&0&1&0&1&0&0&0&0&0&0&1&0&0&0&1&0&0&0&0&0\\
0&0&1&0&0&0&1&0&0&1&1&0&0&0&0&0&1&0&0&0&0&0\\
1&0&0&0&1&0&0&0&1&0&0&0&0&0&0&0&1&0&0&0&0&0\\
0&0&0&0&0&0&0&0&0&0&0&0&0&1&1&1&0&0&0&0&0&1\\
0&0&0&0&0&0&0&0&0&0&1&1&1&0&0&0&0&0&0&0&0&1\\
0&0&0&0&0&0&0&1&1&1&0&0&0&0&0&0&0&0&0&0&0&1\\
0&0&0&0&1&1&1&0&0&0&0&0&0&0&0&0&0&0&0&0&0&1\\
1&1&1&1&0&0&0&0&0&0&0&0&0&0&0&0&0&0&0&0&0&1\\
0&0&0&0&0&0&0&0&0&0&0&0&0&0&0&0&1&1&1&1&1&0\\
\end{array}
$
\end{tiny}

\vspace{.25cm}
\noindent
and degree sequence: $$(4,4,4,4,5,4,4,4,5,4,4,4,5,4,5,4,4,4,4,4,5,5)$$ Edge set is:
\begin{eqnarray*}
(0,9), (0,12), (0,15), (0,20), (1,6), (1,8), (1,11), (1,20),\\ (2,5), 
(2,7), (2,14), (2,20), (3,4), (3,10), (3,13), (3,20),\\ (4,7), (4,11), (4,15), 
(4,19), (5,8), (5,13), (5,19), (6,12),\\ (6,14), (6,19), (7,12), (7,18), (8,10), 
(8,15), (8,18), (9,11),\\ (9,14), (9,18), (10,14), (10,17), (11,17), (12,13), (12,17), 
(13,16),\\ (14,16), (15,16), (16,21), (17,21), (18,21), (19,21), (20,21)
\end{eqnarray*}

\vspace{.25cm}
\noindent
The sink nodes are vertices $20$ and $21$, which have sets of children$\{0,1,2,3,21\}$ and $\{16,17,18,19,20\}$.

\vspace{.25cm}
\noindent
Set $S_{\Gamma_{0},3}$  is given below. All elements of $S_{\Gamma_{0},3}$ contain no vertex of degree $5$, and  contain at least one child of a sink node. 
%
%
%
%
\begin{eqnarray*}
S_{\Gamma_{0},3} &=&\{(0,10,19), (1,7,16), (2,15,17), (3,6,18)\}\\
\end{eqnarray*}

\vspace{.25cm}
\noindent 
Note the elements of $S_{\Gamma_{0},3}$ do not intersect. 

The neighbourhoods of each vertex of degree $4$ are given below, with vertices of degree $5$ shown in bold. It can be seen that every vertex of degree $4$ has at least two neighbours of degree $4$.
\begin{eqnarray*}
Nb(0) &=& 9, {\bf 12}, 15, {\bf 20}\\
Nb(1) &=& 6, {\bf 8}, 11, {\bf 20}\\
Nb(2) &=& 5, 7, {\bf 14}, {\bf 20}\\
Nb(3) &=& {\bf 4}, 10, 13, {\bf 20}\\
Nb(5) &=& 2, {\bf 8}, 13, 19\\ 
Nb(6) &=& 1, {\bf 12}, {\bf 14}, 19\\
Nb(7) &=& 2, {\bf 4}, {\bf 12}, 18\\
Nb(9) &=& 0, 11, {\bf 14}, 18\\
Nb(10) &=& 3, {\bf 8}, {\bf 14}, 17\\
Nb(11) &=& 1, {\bf 4}, 9, 17\\
Nb(13) &=& 3,5,{\bf 12}, 16\\
Nb(15) &=& 0, {\bf 4}, {\bf 8}, 16\\
Nb(16) &=&13, {\bf 14}, 15, {\bf 21}\\
Nb(17) &=& 10, 11, {\bf 12}, {\bf 21}\\
Nb(18) &=& 7,{\bf 8},9,{\bf 21}\\
Nb(19) &=&{\bf 4},5,6,{\bf 21}
\end{eqnarray*}

\vspace{.25cm}
\noindent
$\Gamma_{1}$ has incidence array:

\vspace{.25cm}
\noindent
\begin{tiny}
$
\begin{array}{cccccccccccccccccccccccc}
0&0&0&0&0&0&0&1&0&0&0&1&0&0&0&1&0&0&0&0&1&0\\
0&0&0&0&0&0&1&0&0&0&1&0&0&1&0&0&0&0&0&0&1&0\\
0&0&0&0&0&1&0&0&0&1&0&0&1&0&0&0&0&0&0&0&1&0\\
0&0&0&0&1&0&0&0&1&0&0&0&0&0&1&0&0&0&0&0&1&0\\
0&0&0&1&0&0&0&0&0&0&0&1&0&1&0&0&0&0&0&1&0&0\\
0&0&1&0&0&0&0&0&0&0&1&0&0&0&0&1&0&0&0&1&0&0\\
0&1&0&0&0&0&0&0&1&0&0&0&1&0&0&0&0&0&0&1&0&0\\
1&0&0&0&0&0&0&0&0&1&0&0&0&0&1&0&0&0&0&1&0&0\\
0&0&0&1&0&0&1&0&0&0&0&0&0&0&0&1&0&0&1&0&0&0\\
0&0&1&0&0&0&0&1&0&0&0&0&0&1&0&0&0&0&1&0&0&0\\
0&1&0&0&0&1&0&0&0&0&0&0&0&0&1&0&0&0&1&0&0&0\\
1&0&0&0&1&0&0&0&0&0&0&0&1&0&0&0&0&0&1&0&0&0\\
0&0&1&0&0&0&1&0&0&0&0&1&0&0&1&0&0&1&0&0&0&0\\
0&1&0&0&1&0&0&0&0&1&0&0&0&0&0&1&0&1&0&0&0&0\\
0&0&0&1&0&0&0&1&0&0&1&0&1&0&0&0&1&0&0&0&0&0\\
1&0&0&0&0&1&0&0&1&0&0&0&0&1&0&0&1&0&0&0&0&0\\
0&0&0&0&0&0&0&0&0&0&0&0&0&0&1&1&0&0&0&0&0&1\\
0&0&0&0&0&0&0&0&0&0&0&0&1&1&0&0&0&0&0&0&0&1\\
0&0&0&0&0&0&0&0&1&1&1&1&0&0&0&0&0&0&0&0&0&1\\
0&0&0&0&1&1&1&1&0&0&0&0&0&0&0&0&0&0&0&0&0&1\\
1&1&1&1&0&0&0&0&0&0&0&0&0&0&0&0&0&0&0&0&0&1\\
0&0&0&0&0&0&0&0&0&0&0&0&0&0&0&0&1&1&1&1&1&0\\
\end{array}
$
\end{tiny}

\vspace{.25cm}
\noindent
and degree sequence: $$(4,4,4,4,4,4,4,4,4,4,4,4,5,5,5,5,3,3,5,5,5,5)$$ Edge set is:
\begin{eqnarray*}
(0,7), (0,11), (0,15), (0,20), (1,6), (1,10), (1,13), (1,20),\\ (2,5), 
(2,9), (2,12), (2,20), (3,4), (3,8), (3,14), (3,20),\\ (4,11), (4,13), (4,19), 
(5,10), (5,15), (5,19), (6,8), (6,12),\\ (6,19), (7,9), (7,14), (7,19), (8,15), 
(8,18), (9,13), (9,18),\\ (10,14), (10,18), (11,12), (11,18), (12,14), (12,17), (13,15), 
(13,17),\\ (14,16), (15,16), (16,21), (17,21), (18,21), (19,21), (20,21)
\end{eqnarray*}

\vspace{.25cm}
\noindent
The sink nodes are vertices $18,19,20$ and $21$, which have sets of children$\{8,9,10,11,21\}$, $\{4,5,6,7,21\}$, $\{0,1,2,3,21\}$ and $\{16,17,18,19,20\}$. respectively.

\vspace{.25cm}
\noindent
Set $S_{\Gamma_{1},3}$ is given below. All elements contain no vertex of degree $5$ and contains at least one child of a sink node. No element of $S_{\Gamma_{1},3}$ contains both vertices of degree $3$ ($16$ and $17$). 
\begin{eqnarray*}
S_{\Gamma_{1},3} &=&\{(0,10,17), (1,11,16), (2,4,16), (3,5,17), (6,9,16), (7,8,17)\}
\end{eqnarray*}

\vspace{.25cm}
\noindent
$\Gamma_{2}$ has incidence array:

\vspace{.25cm}
\noindent
\begin{tiny}
$
\begin{array}{cccccccccccccccccccccccc}
0&0&0&0&0&0&0&0&0&0&0&0&1&0&0&1&0&0&0&0&1&0\\
0&0&0&0&0&0&1&0&0&1&0&1&0&0&0&0&0&0&0&0&1&0\\
0&0&0&0&0&1&0&0&1&0&0&0&0&0&1&0&0&0&0&0&1&0\\
0&0&0&0&1&0&0&1&0&0&1&0&0&0&0&0&0&0&0&0&1&0\\
0&0&0&1&0&0&0&0&0&0&0&1&0&0&1&0&0&0&0&1&0&0\\
0&0&1&0&0&0&0&0&0&0&1&0&0&0&0&1&0&0&0&1&0&0\\
0&1&0&0&0&0&0&1&0&0&0&0&1&1&0&0&0&0&0&1&0&0\\
0&0&0&1&0&0&1&0&0&0&0&0&0&0&0&1&0&0&1&0&0&0\\
0&0&1&0&0&0&0&0&0&0&0&1&0&1&0&0&0&0&1&0&0&0\\
0&1&0&0&0&0&0&0&0&0&1&0&0&0&1&0&0&0&1&0&0&0\\
0&0&0&1&0&1&0&0&0&1&0&0&0&1&0&0&0&1&0&0&0&0\\
0&1&0&0&1&0&0&0&1&0&0&0&0&0&0&1&0&1&0&0&0&0\\
1&0&0&0&0&0&1&0&0&0&0&0&0&0&1&0&0&1&0&0&0&0\\
0&0&0&0&0&0&1&0&1&0&1&0&0&0&0&0&1&0&0&0&0&0\\
0&0&1&0&1&0&0&0&0&1&0&0&1&0&0&0&1&0&0&0&0&0\\
1&0&0&0&0&1&0&1&0&0&0&1&0&0&0&0&1&0&0&0&0&0\\
0&0&0&0&0&0&0&0&0&0&0&0&0&1&1&1&0&0&0&0&0&1\\
0&0&0&0&0&0&0&0&0&0&1&1&1&0&0&0&0&0&0&0&0&1\\
0&0&0&0&0&0&0&1&1&1&0&0&0&0&0&0&0&0&0&0&0&1\\
0&0&0&0&1&1&1&0&0&0&0&0&0&0&0&0&0&0&0&0&0&1\\
1&1&1&1&0&0&0&0&0&0&0&0&0&0&0&0&0&0&0&0&0&1\\
0&0&0&0&0&0&0&0&0&0&0&0&0&0&0&0&1&1&1&1&1&0\\
\end{array}
$
\end{tiny}

\vspace{.25cm}
\noindent
and degree sequence: $$(3,4,4,4,4,4,5,4,4,4,5,5,4,4,5,5,4,4,4,4,5,5)$$ Edge set is:
\begin{eqnarray*}
(0,12), (0,15), (0,20), (1,6), (1,9), (1,11), (1,20), (2,5), \\(2,8), 
(2,14), (2,20), (3,4), (3,7), (3,10), (3,20), (4,11), \\(4,14), (4,19), (5,10), 
(5,15), (5,19), (6,7), (6,12), (6,13), \\(6,19), (7,15), (7,18), (8,11), (8,13), 
(8,18), (9,10), (9,14),\\ (9,18), (10,13), (10,17), (11,15), (11,17), (12,14), (12,17), 
(13,16),\\ (14,16), (15,16), (16,21), (17,21), (18,21), (19,21), (20,21)
\end{eqnarray*}

\vspace{.25cm}
\noindent
The sink nodes are vertices $11$ and $21$, which have sets of children$\{1,4,8,15,17\}$ and $\{16,17,18,19,20\}$.

\vspace{.25cm}
\noindent
Sets $S_{\Gamma_{2},3}$ is given below. All elements of degree $5$ are highlighted in bold. Every element of $S_{\Gamma_{2},3}$ contains no vertex of degree $5$ and contains at least one child of a sink node.

\begin{eqnarray*}
S_{\Gamma_{2},3} &=&\{(0,4,13), (0,4,18), (0,8,19), (0,9,19), \\&&(2,7,17), (3,8,12), (5,12,18)\}
\end{eqnarray*}
%

The neighbourhoods of each vertex of degree $4$ are given below, with vertices of degree $5$ shown in bold. It can be seen that every vertex of degree $4$ has at least two neighbours of degree $4$ or $3$, except vertices $1$, $16$ and $17$.
\begin{eqnarray*}
Nb(1) &=& {\bf 6}, 9, {\bf 11}, {\bf 20}\\
Nb(2) &=&5,8,{\bf 14},{\bf 20} \\
Nb(3) &=&4,7,{\bf 10},{\bf 20}\\
Nb(4) &=&3,{\bf 11},{\bf 14},19  \\
Nb(5) &=&2,{\bf 10},{\bf 15},19 \\
Nb(7) &=&3,{\bf 6},{\bf 15},18 \\
Nb(8) &=&2,{\bf 11},13,18 \\
Nb(9) &=&1,{\bf 10},{\bf 14},18 \\
Nb(12) &=&0,{\bf 6},{\bf 14},17 \\
Nb(13) &=&{\bf 6},8,{\bf 10},16  \\
Nb(16) &=&13,{\bf 14},{\bf 15},{\bf 21} \\
Nb(17) &=&{\bf 10},{\bf 11},12,{\bf 21}\\
Nb(18) &=&7,8,9,{\bf 21}\\
Nb(19) &=&4,5,{\bf 6},{\bf 21}
\end{eqnarray*}

The set of vertices of degree less than $5$ that are adjacent to $1$, $16$ and $17$ is $\{9,13,12\}$. No two of these vertices are in an element of $S_{\Gamma_{2},3}$ with the vertex of degree $3$ ($0$). 

\section*{Appendix D}\label{appendix23}
The $7$ extremal graphs $(23,50)$. Note that graphs $\Gamma_{0}-\Gamma_{1}$ have $(deg_{3},deg_{4},deg_{5})=(0,15,8)$, graphs $\Gamma_{2}-\Gamma_{4}$ have  $(deg_{3},deg_{4},deg_{5})=(1,13,9)$ and graphs $\Gamma_{5}-\Gamma_{6}$ have $(deg_{3},deg_{4},deg_{5})=(2,11,10)$.

\vspace{.25cm}
\noindent
$\Gamma_{0}$ has incidence array:

\vspace{.25cm}
\noindent
\begin{tiny}
$
\begin{array}{ccccccccccccccccccccccccc}
0&0&0&0&0&0&0&1&0&0&1&0&0&0&0&0&1&0&0&0&0&1&0\\
0&0&0&0&0&0&1&0&0&1&0&0&0&0&0&1&0&0&0&0&0&1&0\\
0&0&0&0&0&1&0&0&0&0&0&0&0&1&1&0&0&0&0&0&0&1&0\\
0&0&0&0&1&0&0&0&1&0&0&0&1&0&0&0&0&0&0&0&0&1&0\\
0&0&0&1&0&0&0&0&0&0&1&0&0&1&0&0&0&0&0&0&1&0&0\\
0&0&1&0&0&0&0&0&1&0&0&0&0&0&0&1&0&0&0&0&1&0&0\\
0&1&0&0&0&0&0&0&0&0&0&1&0&0&0&0&1&0&0&0&1&0&0\\
1&0&0&0&0&0&0&0&0&1&0&0&1&0&1&0&0&0&0&0&1&0&0\\
0&0&0&1&0&1&0&0&0&0&0&0&0&0&0&0&1&0&0&1&0&0&0\\
0&1&0&0&0&0&0&1&0&0&0&0&0&1&0&0&0&0&0&1&0&0&0\\
1&0&0&0&1&0&0&0&0&0&0&1&0&0&0&1&0&0&0&1&0&0&0\\
0&0&0&0&0&0&1&0&0&0&1&0&0&0&1&0&0&0&1&0&0&0&0\\
0&0&0&1&0&0&0&1&0&0&0&0&0&0&0&1&0&0&1&0&0&0&0\\
0&0&1&0&1&0&0&0&0&1&0&0&0&0&0&0&1&0&1&0&0&0&0\\
0&0&1&0&0&0&0&1&0&0&0&1&0&0&0&0&0&1&0&0&0&0&0\\
0&1&0&0&0&1&0&0&0&0&1&0&1&0&0&0&0&1&0&0&0&0&0\\
1&0&0&0&0&0&1&0&1&0&0&0&0&1&0&0&0&1&0&0&0&0&0\\
0&0&0&0&0&0&0&0&0&0&0&0&0&0&1&1&1&0&0&0&0&0&1\\
0&0&0&0&0&0&0&0&0&0&0&1&1&1&0&0&0&0&0&0&0&0&1\\
0&0&0&0&0&0&0&0&1&1&1&0&0&0&0&0&0&0&0&0&0&0&1\\
0&0&0&0&1&1&1&1&0&0&0&0&0&0&0&0&0&0&0&0&0&0&1\\
1&1&1&1&0&0&0&0&0&0&0&0&0&0&0&0&0&0&0&0&0&0&1\\
0&0&0&0&0&0&0&0&0&0&0&0&0&0&0&0&0&1&1&1&1&1&0\\
\end{array}
$
\end{tiny}

\vspace{.25cm}
\noindent
and degree sequence: $$(4,4,4,4,4,4,4,5,4,4,5,4,4,5,4,5,5,4,4,4,5,5,5)$$ Edge set is:
\begin{eqnarray*}
(0,7), (0,10), (0,16), (0,21), (1,6), (1,9), (1,15), (1,21), (2,5), \\
(2,13), (2,14), (2,21), (3,4), (3,8), (3,12), (3,21), (4,10), (4,13),\\
 (4,20), (5,8), (5,15), (5,20), (6,11), (6,16), (6,20), (7,9), (7,12),\\ 
(7,14), (7,20), (8,16), (8,19), (9,13), (9,19), (10,11), (10,15), (10,19),\\
 (11,14), (11,18), (12,15), (12,18), (13,16), (13,18), (14,17),\\
 (15,17), (16,17), (17,22), (18,22), (19,22), (20,22), (21,22)
\end{eqnarray*}

\vspace{.25cm}
\noindent
The sink nodes are vertices $20$ and $22$, which have sets of children$\{4,5,6,7,22\}$ and $\{17,18,19,20,21\}$.

\vspace{.25cm}
\noindent
Sets $S_{\Gamma_{0},3}$ and $S_{\Gamma_{0},4}$ are given below. 
 Neither set has elements containing vertices of degree $5$. Every element of $S_{\Gamma_{0},4}$ contains at least one child of a sink node.
\begin{eqnarray*}
S_{\Gamma_{0},3} &=&\{(0,5,18), (1,4,14), (1,8,14), (1,8,18), (2,6,12), (2,6,19),\\&& (2,12,19), (3,9,11), (3,9,17), (5,9,11), (6,12,19)\}\\
S_{\Gamma_{0},4}&= &\{(2,6,12,19)\}
\end{eqnarray*}

\vspace{.25cm} The subgraph on vertices of degree $4$ has vertex set $\{0,1,2,3,4,5,6,8,9,11,12,14,17,18,19\}$ and degree sequence $$(0,2,2,3,1,2,2,3,2,3,2,3,1,2,2)$$. The only vertex of degree $0$ in this subgraph is vertex $0$. This vertex has a common neighbour of degree $5$ with the vertices of degree $1$ in the subgraph (vertices $4$ and $17$), namely $10$ and $17$ respectively. 

\vspace{.25cm}
\noindent
$\Gamma_{1}$ has incidence array:

\vspace{.25cm}
\noindent
\begin{tiny}
$
\begin{array}{ccccccccccccccccccccccccc}
0&0&0&0&0&0&0&1&0&0&0&0&0&1&0&0&1&0&0&0&0&1&0\\
0&0&0&0&0&0&1&0&0&0&1&0&1&0&0&0&0&0&0&0&0&1&0\\
0&0&0&0&0&1&0&0&0&1&0&0&0&0&0&1&0&0&0&0&0&1&0\\
0&0&0&0&1&0&0&0&1&0&0&1&0&0&1&0&0&0&0&0&0&1&0\\
0&0&0&1&0&0&0&0&0&0&1&0&0&1&0&1&0&0&0&0&1&0&0\\
0&0&1&0&0&0&0&0&1&0&0&0&0&0&0&0&1&0&0&0&1&0&0\\
0&1&0&0&0&0&0&0&0&1&0&1&0&0&0&0&0&0&0&0&1&0&0\\
1&0&0&0&0&0&0&0&0&0&0&0&1&0&1&0&0&0&0&0&1&0&0\\
0&0&0&1&0&1&0&0&0&0&0&0&1&0&0&0&0&0&0&1&0&0&0\\
0&0&1&0&0&0&1&0&0&0&0&0&0&1&1&0&0&0&0&1&0&0&0\\
0&1&0&0&1&0&0&0&0&0&0&0&0&0&0&0&1&0&0&1&0&0&0\\
0&0&0&1&0&0&1&0&0&0&0&0&0&0&0&0&1&0&1&0&0&0&0\\
0&1&0&0&0&0&0&1&1&0&0&0&0&0&0&1&0&0&1&0&0&0&0\\
1&0&0&0&1&0&0&0&0&1&0&0&0&0&0&0&0&0&1&0&0&0&0\\
0&0&0&1&0&0&0&1&0&1&0&0&0&0&0&0&0&1&0&0&0&0&0\\
0&0&1&0&1&0&0&0&0&0&0&0&1&0&0&0&0&1&0&0&0&0&0\\
1&0&0&0&0&1&0&0&0&0&1&1&0&0&0&0&0&1&0&0&0&0&0\\
0&0&0&0&0&0&0&0&0&0&0&0&0&0&1&1&1&0&0&0&0&0&1\\
0&0&0&0&0&0&0&0&0&0&0&1&1&1&0&0&0&0&0&0&0&0&1\\
0&0&0&0&0&0&0&0&1&1&1&0&0&0&0&0&0&0&0&0&0&0&1\\
0&0&0&0&1&1&1&1&0&0&0&0&0&0&0&0&0&0&0&0&0&0&1\\
1&1&1&1&0&0&0&0&0&0&0&0&0&0&0&0&0&0&0&0&0&0&1\\
0&0&0&0&0&0&0&0&0&0&0&0&0&0&0&0&0&1&1&1&1&1&0\\
\end{array}
$
\end{tiny}

\vspace{.25cm}
\noindent
and degree sequence: $$(4,4,4,5,5,4,4,4,4,5,4,4,5,4,4,4,5,4,4,4,5,5,5)$$ Edge set is:
\begin{eqnarray*}
(0,7), (0,13), (0,16), (0,21), (1,6), (1,10), (1,12), (1,21), (2,5), \\
(2,9), (2,15), (2,21), (3,4), (3,8), (3,11), (3,14), (3,21), (4,10),\\ (4,13), 
(4,15), (4,20), (5,8), (5,16), (5,20), (6,9), (6,11), (6,20),\\ (7,12), (7,14), 
(7,20), (8,12), (8,19), (9,13), (9,14), (9,19), (10,16),\\ (10,19), (11,16), (11,18), 
(12,15), (12,18), (13,18), (14,17),\\ (15,17), (16,17), (17,22), (18,22), (19,22), (20,22), 
(21,22)
\end{eqnarray*}

\vspace{.25cm}
\noindent
Sets $S_{\Gamma_{1},3}$ and $S_{\Gamma_{1},4}$ are given below. Vertices of degree $5$ are written in bold.
\begin{eqnarray*}
S_{\Gamma_{1},3} &=&\{ (0,6,8), (0,6,15), (0,15,19), (1,5,13), (1,5,14), (1,13,17), \\
&&(2,7,10), (2,7,11), (2,10,18), (5,14,18), (6,8,17), (7,11,19),\\
&& (8,13,17), ({\bf 9},{\bf 12},{\bf 16}), (10,14,18), (11,15,19), \}\\
S_{\Gamma_{1},4}&= &\{\}
\end{eqnarray*}

\vspace{.25cm} The subgraph on vertices of degree $4$ has vertex set $\{0,1,2,5,6,7,8,10,11,13,14,15,17,18,19\}$ and degree sequence $(2,2,2,2,2,2,2,2,2,2,2,2,2,2,2)$.

\vspace{.25cm}
\noindent
$\Gamma_{2}$ has incidence array:

\vspace{.25cm}
\noindent
\begin{tiny}
$
\begin{array}{ccccccccccccccccccccccccc}
0&0&0&0&0&0&0&1&0&0&1&0&0&0&0&0&1&0&0&0&0&1&0\\
0&0&0&0&0&0&1&0&0&1&0&0&0&1&0&0&0&0&0&0&0&1&0\\
0&0&0&0&0&1&0&0&0&0&0&0&1&0&0&1&0&0&0&0&0&1&0\\
0&0&0&0&1&0&0&0&1&0&0&0&0&0&1&0&0&0&0&0&0&1&0\\
0&0&0&1&0&0&0&0&0&0&1&0&1&0&0&0&0&0&0&0&1&0&0\\
0&0&1&0&0&0&0&0&1&0&0&1&0&0&0&0&0&0&0&0&1&0&0\\
0&1&0&0&0&0&0&0&0&0&0&0&0&0&0&0&1&0&0&0&1&0&0\\
1&0&0&0&0&0&0&0&0&0&0&0&0&1&1&0&0&0&0&0&1&0&0\\
0&0&0&1&0&1&0&0&0&0&0&0&0&1&0&0&1&0&0&1&0&0&0\\
0&1&0&0&0&0&0&0&0&0&0&0&1&0&1&0&0&0&0&1&0&0&0\\
1&0&0&0&1&0&0&0&0&0&0&1&0&0&0&1&0&0&0&1&0&0&0\\
0&0&0&0&0&1&0&0&0&0&1&0&0&0&1&0&0&0&1&0&0&0&0\\
0&0&1&0&1&0&0&0&0&1&0&0&0&0&0&0&1&0&1&0&0&0&0\\
0&1&0&0&0&0&0&1&1&0&0&0&0&0&0&1&0&0&1&0&0&0&0\\
0&0&0&1&0&0&0&1&0&1&0&1&0&0&0&0&0&1&0&0&0&0&0\\
0&0&1&0&0&0&0&0&0&0&1&0&0&1&0&0&0&1&0&0&0&0&0\\
1&0&0&0&0&0&1&0&1&0&0&0&1&0&0&0&0&1&0&0&0&0&0\\
0&0&0&0&0&0&0&0&0&0&0&0&0&0&1&1&1&0&0&0&0&0&1\\
0&0&0&0&0&0&0&0&0&0&0&1&1&1&0&0&0&0&0&0&0&0&1\\
0&0&0&0&0&0&0&0&1&1&1&0&0&0&0&0&0&0&0&0&0&0&1\\
0&0&0&0&1&1&1&1&0&0&0&0&0&0&0&0&0&0&0&0&0&0&1\\
1&1&1&1&0&0&0&0&0&0&0&0&0&0&0&0&0&0&0&0&0&0&1\\
0&0&0&0&0&0&0&0&0&0&0&0&0&0&0&0&0&1&1&1&1&1&0\\
\end{array}
$
\end{tiny}

\vspace{.25cm}
\noindent
and degree sequence: $$(4,4,4,4,4,4,3,4,5,4,5,4,5,5,5,4,5,4,4,4,5,5,5)$$ Edge set is:
\begin{eqnarray*}
(0,7), (0,10), (0,16), (0,21), (1,6), (1,9), (1,13), (1,21), (2,5), \\
(2,12), (2,15), (2,21), (3,4), (3,8), (3,14), (3,21), (4,10), (4,12),\\ (4,20), 
(5,8), (5,11), (5,20), (6,16), (6,20), (7,13), (7,14), (7,20),\\ (8,13), (8,16), 
(8,19), (9,12), (9,14), (9,19), (10,11), (10,15), (10,19),\\ (11,14), (11,18), (12,16), 
(12,18), (13,15), (13,18), (14,17),\\ (15,17), (16,17), (17,22), (18,22), (19,22), (20,22), 
(21,22)
\end{eqnarray*}

\vspace{.25cm}
\noindent
Sets $S_{\Gamma_{2},3}$ and $S_{\Gamma_{2},4}$ are given below. Vertices of degree $5$ are written in bold.
\begin{eqnarray*}
S_{\Gamma_{2},3} &=&\{(0,5,9), (1,4,17), (1,5,17), (2,6,{\bf 14}), (2,6,19), \\&&(2,7,19), (3,6,15), (3,6,18), (9,15,{\bf 20})\}\\
S_{\Gamma_{2},4}&= &\{\}
\end{eqnarray*}

\vspace{.25cm}
\noindent
$\Gamma_{3}$ has incidence array:

\vspace{.25cm}
\noindent
\begin{tiny}
$
\begin{array}{ccccccccccccccccccccccccc}
0&0&0&0&0&0&0&1&0&0&1&0&0&1&0&0&0&0&0&0&0&1&0\\
0&0&0&0&0&0&1&0&0&1&0&0&0&0&0&0&1&0&0&0&0&1&0\\
0&0&0&0&0&1&0&0&1&0&0&0&1&0&0&1&0&0&0&0&0&1&0\\
0&0&0&0&1&0&0&0&0&0&0&1&0&0&1&0&0&0&0&0&0&1&0\\
0&0&0&1&0&0&0&0&0&0&0&0&0&1&0&0&1&0&0&0&1&0&0\\
0&0&1&0&0&0&0&0&0&0&1&0&0&0&0&0&0&0&0&0&1&0&0\\
0&1&0&0&0&0&0&0&1&0&0&1&0&0&0&0&0&0&0&0&1&0&0\\
1&0&0&0&0&0&0&0&0&1&0&0&1&0&1&0&0&0&0&0&1&0&0\\
0&0&1&0&0&0&1&0&0&0&0&0&0&1&1&0&0&0&0&1&0&0&0\\
0&1&0&0&0&0&0&1&0&0&0&0&0&0&0&1&0&0&0&1&0&0&0\\
1&0&0&0&0&1&0&0&0&0&0&1&0&0&0&0&1&0&0&1&0&0&0\\
0&0&0&1&0&0&1&0&0&0&1&0&0&0&0&1&0&0&1&0&0&0&0\\
0&0&1&0&0&0&0&1&0&0&0&0&0&0&0&0&1&0&1&0&0&0&0\\
1&0&0&0&1&0&0&0&1&0&0&0&0&0&0&0&0&0&1&0&0&0&0\\
0&0&0&1&0&0&0&1&1&0&0&0&0&0&0&0&0&1&0&0&0&0&0\\
0&0&1&0&0&0&0&0&0&1&0&1&0&0&0&0&0&1&0&0&0&0&0\\
0&1&0&0&1&0&0&0&0&0&1&0&1&0&0&0&0&1&0&0&0&0&0\\
0&0&0&0&0&0&0&0&0&0&0&0&0&0&1&1&1&0&0&0&0&0&1\\
0&0&0&0&0&0&0&0&0&0&0&1&1&1&0&0&0&0&0&0&0&0&1\\
0&0&0&0&0&0&0&0&1&1&1&0&0&0&0&0&0&0&0&0&0&0&1\\
0&0&0&0&1&1&1&1&0&0&0&0&0&0&0&0&0&0&0&0&0&0&1\\
1&1&1&1&0&0&0&0&0&0&0&0&0&0&0&0&0&0&0&0&0&0&1\\
0&0&0&0&0&0&0&0&0&0&0&0&0&0&0&0&0&1&1&1&1&1&0\\
\end{array}
$
\end{tiny}

\vspace{.25cm}
\noindent
and degree sequence: $$(4,4,5,4,4,3,4,5,5,4,5,5,4,4,4,4,5,4,4,4,5,5,5)$$ Edge set is:
\begin{eqnarray*}
(0,7), (0,10), (0,13), (0,21), (1,6), (1,9), (1,16), (1,21), (2,5),\\ 
(2,8), (2,12), (2,15), (2,21), (3,4), (3,11), (3,14), (3,21), (4,13),\\ (4,16), 
(4,20), (5,10), (5,20), (6,8), (6,11), (6,20), (7,9), (7,12),\\ (7,14), (7,20), 
(8,13), (8,14), (8,19), (9,15), (9,19), (10,11), (10,16),\\ (10,19), (11,15), (11,18), 
(12,16), (12,18), (13,18), (14,17),\\ (15,17), (16,17), (17,22), (18,22), (19,22), (20,22), 
(21,22)
\end{eqnarray*}

\vspace{.25cm}
\noindent
The sink nodes are vertices $21$ and $22$, which have sets of children$\{0,1,2,3,22\}$ and $\{17,18,19,20,21\}$

\vspace{.25cm}
\noindent
Sets $S_{\Gamma_{3},3}$ and $S_{\Gamma_{3},4}$ are given below.  Neither set has elements containing vertices of degree $5$. Every element of $S_{\Gamma_{3},4}$ contains the vertex of degree $3$ at least one vertex of degree $4$ that is the child of a sink node.
\begin{eqnarray*}
S_{\Gamma_{3},3} &=&\{(0,6,17), (1,5,13), (1,5,14), (1,5,18), (1,14,18), (3,5,9),\\&& (3,12,19), (5,9,13), (5,9,18), (5,13,17), (5,14,18)\}\\
S_{\Gamma_{3},4}&= &\{(1,5,14,18)\}
\end{eqnarray*}

\vspace{.25cm}
\noindent
$\Gamma_{4}$ has incidence array:

\vspace{.25cm}
\noindent
\begin{tiny}
$
\begin{array}{ccccccccccccccccccccccccc}
0&0&0&0&0&0&0&1&0&0&1&0&0&1&0&0&0&0&0&0&0&1&0\\
0&0&0&0&0&0&1&0&0&1&0&0&1&0&0&0&1&0&0&0&0&1&0\\
0&0&0&0&0&1&0&0&1&0&0&1&0&0&0&0&0&0&0&0&0&1&0\\
0&0&0&0&1&0&0&0&0&0&0&0&0&0&0&1&0&0&0&0&0&1&0\\
0&0&0&1&0&0&0&0&0&1&0&0&0&1&0&0&0&0&0&0&1&0&0\\
0&0&1&0&0&0&0&0&0&0&1&0&0&0&0&0&1&0&0&0&1&0&0\\
0&1&0&0&0&0&0&0&1&0&0&0&0&0&0&1&0&0&0&0&1&0&0\\
1&0&0&0&0&0&0&0&0&0&0&1&0&0&1&0&0&0&0&0&1&0&0\\
0&0&1&0&0&0&1&0&0&0&0&0&0&1&1&0&0&0&0&1&0&0&0\\
0&1&0&0&1&0&0&0&0&0&0&1&0&0&0&0&0&0&0&1&0&0&0\\
1&0&0&0&0&1&0&0&0&0&0&0&1&0&0&1&0&0&0&1&0&0&0\\
0&0&1&0&0&0&0&1&0&1&0&0&0&0&0&1&0&0&1&0&0&0&0\\
0&1&0&0&0&0&0&0&0&0&1&0&0&0&1&0&0&0&1&0&0&0&0\\
1&0&0&0&1&0&0&0&1&0&0&0&0&0&0&0&1&0&1&0&0&0&0\\
0&0&0&0&0&0&0&1&1&0&0&0&1&0&0&0&0&1&0&0&0&0&0\\
0&0&0&1&0&0&1&0&0&0&1&1&0&0&0&0&0&1&0&0&0&0&0\\
0&1&0&0&0&1&0&0&0&0&0&0&0&1&0&0&0&1&0&0&0&0&0\\
0&0&0&0&0&0&0&0&0&0&0&0&0&0&1&1&1&0&0&0&0&0&1\\
0&0&0&0&0&0&0&0&0&0&0&1&1&1&0&0&0&0&0&0&0&0&1\\
0&0&0&0&0&0&0&0&1&1&1&0&0&0&0&0&0&0&0&0&0&0&1\\
0&0&0&0&1&1&1&1&0&0&0&0&0&0&0&0&0&0&0&0&0&0&1\\
1&1&1&1&0&0&0&0&0&0&0&0&0&0&0&0&0&0&0&0&0&0&1\\
0&0&0&0&0&0&0&0&0&0&0&0&0&0&0&0&0&1&1&1&1&1&0\\
\end{array}
$
\end{tiny}

\vspace{.25cm}
\noindent
and degree sequence: $$(4,5,4,3,4,4,4,4,5,4,5,5,4,5,4,5,4,4,4,4,5,5,5)$$ Edge set is:
\begin{eqnarray*}
(0,7), (0,10), (0,13), (0,21), (1,6), (1,9), (1,12), (1,16), (1,21), \\
(2,5), (2,8), (2,11), (2,21), (3,4), (3,15), (3,21), (4,9), (4,13),\\ (4,20), 
(5,10), (5,16), (5,20), (6,8), (6,15), (6,20), (7,11), (7,14),\\ (7,20), (8,13), 
(8,14), (8,19), (9,11), (9,19), (10,12), (10,15), (10,19),\\ (11,15), (11,18), (12,14), 
(12,18), (13,16), (13,18), (14,17),\\ (15,17), (16,17), (17,22), (18,22), (19,22), (20,22), 
(21,22)
\end{eqnarray*}

\vspace{.25cm}
\noindent
There is one sink node, vertex $22$ which has set of children $\{17,18,19,20,21\}$.

\vspace{.25cm}
\noindent
Sets $S_{\Gamma_{4},3}$ and $S_{\Gamma_{4},4}$ are given below.  Neither set has elements containing vertices of degree $5$. All elements of $S_{\Gamma_{4},4}$ contain the vertex of degree $3$ at least vertex of degree $4$ that is child of a sink node $4$.
\begin{eqnarray*}
S_{\Gamma_{4},3} &=&\{(0,9,17), (2,4,12), (2,4,17), (3,5,14), (3,5,18),\\&& (3,7,16), (3,7,19), (3,16,19), (5,9,14), (7,16,19)\}\\
S_{\Gamma_{4},4}&= &\{(3,7,16,19)\}
\end{eqnarray*}

\vspace{.25cm}
\noindent
$\Gamma_{5}$ has incidence array:

\vspace{.25cm}
\noindent
\begin{tiny}
$
\begin{array}{ccccccccccccccccccccccccc}
0&0&0&0&0&0&0&0&0&0&1&0&0&1&0&0&0&0&0&0&0&1&0\\
0&0&0&0&0&0&0&1&0&1&0&0&0&0&0&0&1&0&0&0&0&1&0\\
0&0&0&0&0&0&1&0&0&0&0&0&1&0&0&1&0&0&0&0&0&1&0\\
0&0&0&0&0&1&0&0&1&0&0&0&0&0&1&0&0&0&0&0&0&1&0\\
0&0&0&0&0&0&0&0&0&0&1&0&0&0&0&1&0&0&0&0&1&0&0\\
0&0&0&1&0&0&0&0&0&1&0&0&1&0&0&0&0&0&0&0&1&0&0\\
0&0&1&0&0&0&0&0&1&0&0&1&0&0&0&0&0&0&0&0&1&0&0\\
0&1&0&0&0&0&0&0&0&0&0&0&0&1&1&0&0&0&0&0&1&0&0\\
0&0&0&1&0&0&1&0&0&0&0&0&0&1&0&0&1&0&0&1&0&0&0\\
0&1&0&0&0&1&0&0&0&0&0&1&0&0&0&1&0&0&0&1&0&0&0\\
1&0&0&0&1&0&0&0&0&0&0&0&1&0&1&0&0&0&0&1&0&0&0\\
0&0&0&0&0&0&1&0&0&1&0&0&0&0&1&0&0&0&1&0&0&0&0\\
0&0&1&0&0&1&0&0&0&0&1&0&0&0&0&0&1&0&1&0&0&0&0\\
1&0&0&0&0&0&0&1&1&0&0&0&0&0&0&1&0&0&1&0&0&0&0\\
0&0&0&1&0&0&0&1&0&0&1&1&0&0&0&0&0&1&0&0&0&0&0\\
0&0&1&0&1&0&0&0&0&1&0&0&0&1&0&0&0&1&0&0&0&0&0\\
0&1&0&0&0&0&0&0&1&0&0&0&1&0&0&0&0&1&0&0&0&0&0\\
0&0&0&0&0&0&0&0&0&0&0&0&0&0&1&1&1&0&0&0&0&0&1\\
0&0&0&0&0&0&0&0&0&0&0&1&1&1&0&0&0&0&0&0&0&0&1\\
0&0&0&0&0&0&0&0&1&1&1&0&0&0&0&0&0&0&0&0&0&0&1\\
0&0&0&0&1&1&1&1&0&0&0&0&0&0&0&0&0&0&0&0&0&0&1\\
1&1&1&1&0&0&0&0&0&0&0&0&0&0&0&0&0&0&0&0&0&0&1\\
0&0&0&0&0&0&0&0&0&0&0&0&0&0&0&0&0&1&1&1&1&1&0\\
\end{array}
$
\end{tiny}

\vspace{.25cm}
\noindent
and degree sequence: $$(3,4,4,4,3,4,4,4,5,5,5,4,5,5,5,5,4,4,4,4,5,5,5)$$ Edge set is:
\begin{eqnarray*}
(0,10), (0,13), (0,21), (1,7), (1,9), (1,16), (1,21), (2,6), (2,12), \\
(2,15), (2,21), (3,5), (3,8), (3,14), (3,21), (4,10), (4,15), (4,20),\\ (5,9), 
(5,12), (5,20), (6,8), (6,11), (6,20), (7,13), (7,14), (7,20),\\ (8,13), (8,16), 
(8,19), (9,11), (9,15), (9,19), (10,12), (10,14), (10,19),\\ (11,14), (11,18), (12,16), 
(12,18), (13,15), (13,18), (14,17),\\ (15,17), (16,17), (17,22), (18,22), (19,22), (20,22), 
(21,22)
\end{eqnarray*}

\vspace{.25cm}
\noindent
Sets $S_{\Gamma_{5},3}$ and $S_{\Gamma_{5},4}$ are given below. Vertices of degree $5$ are written in bold.
\begin{eqnarray*}
S_{\Gamma_{5},3} &=&\{(0,5,17), (0,6,17), (0,11,16), (0,16,{\bf 20}), (1,4,18),\\&& (1,6,{\bf 10}), (2,7,19), (3,4,18), (4,11,16), (4,11,{\bf 21})\}\\
S_{\Gamma_{5},4}&= &\{\}
\end{eqnarray*}

\vspace{.25cm}
\noindent
$\Gamma_{6}$ has incidence array:

\vspace{.25cm}
\noindent
\begin{tiny}
$
\begin{array}{ccccccccccccccccccccccccc}
0&0&0&0&0&0&0&1&0&0&0&1&0&0&1&0&0&0&0&0&0&1&0\\
0&0&0&0&0&0&1&0&0&0&1&0&0&1&0&0&1&0&0&0&0&1&0\\
0&0&0&0&0&1&0&0&0&1&0&0&0&0&0&1&0&0&0&0&0&1&0\\
0&0&0&0&1&0&0&0&1&0&0&0&1&0&0&0&0&0&0&0&0&1&0\\
0&0&0&1&0&0&0&0&0&0&0&1&0&0&0&0&1&0&0&0&1&0&0\\
0&0&1&0&0&0&0&0&0&0&1&0&0&0&1&0&0&0&0&0&1&0&0\\
0&1&0&0&0&0&0&0&1&0&0&0&0&0&0&1&0&0&0&0&1&0&0\\
1&0&0&0&0&0&0&0&0&1&0&0&1&0&0&0&0&0&0&0&1&0&0\\
0&0&0&1&0&0&1&0&0&0&0&0&0&0&1&0&0&0&0&1&0&0&0\\
0&0&1&0&0&0&0&1&0&0&0&0&0&0&0&0&1&0&0&1&0&0&0\\
0&1&0&0&0&1&0&0&0&0&0&0&1&0&0&0&0&0&0&1&0&0&0\\
1&0&0&0&1&0&0&0&0&0&0&0&0&1&0&1&0&0&0&1&0&0&0\\
0&0&0&1&0&0&0&1&0&0&1&0&0&0&0&1&0&0&1&0&0&0&0\\
0&1&0&0&0&0&0&0&0&0&0&1&0&0&0&0&0&0&1&0&0&0&0\\
1&0&0&0&0&1&0&0&1&0&0&0&0&0&0&0&1&0&1&0&0&0&0\\
0&0&1&0&0&0&1&0&0&0&0&1&1&0&0&0&0&1&0&0&0&0&0\\
0&1&0&0&1&0&0&0&0&1&0&0&0&0&1&0&0&1&0&0&0&0&0\\
0&0&0&0&0&0&0&0&0&0&0&0&0&0&0&1&1&0&0&0&0&0&1\\
0&0&0&0&0&0&0&0&0&0&0&0&1&1&1&0&0&0&0&0&0&0&1\\
0&0&0&0&0&0&0&0&1&1&1&1&0&0&0&0&0&0&0&0&0&0&1\\
0&0&0&0&1&1&1&1&0&0&0&0&0&0&0&0&0&0&0&0&0&0&1\\
1&1&1&1&0&0&0&0&0&0&0&0&0&0&0&0&0&0&0&0&0&0&1\\
0&0&0&0&0&0&0&0&0&0&0&0&0&0&0&0&0&1&1&1&1&1&0\\
\end{array}
$
\end{tiny}

\vspace{.25cm}
\noindent
and degree sequence: $$(4,5,4,4,4,4,4,4,4,4,4,5,5,3,5,5,5,3,4,5,5,5,5)$$ Edge set is:
\begin{eqnarray*}
(0,7), (0,11), (0,14), (0,21), (1,6), (1,10), (1,13), (1,16), (1,21),\\ 
(2,5), (2,9), (2,15), (2,21), (3,4), (3,8), (3,12), (3,21), (4,11),\\ (4,16), 
(4,20), (5,10), (5,14), (5,20), (6,8), (6,15), (6,20), (7,9),\\ (7,12), (7,20), 
(8,14), (8,19), (9,16), (9,19), (10,12), (10,19), (11,13),\\ (11,15), (11,19), (12,15), 
(12,18), (13,18), (14,16), (14,18),\\ (15,17), (16,17), (17,22), (18,22), (19,22), (20,22), 
(21,22)
\end{eqnarray*}

\vspace{.25cm}
\noindent
The sink nodes are vertices $19$, $21$ and $22$, which have sets of children$\{8,9,10,11,22\}$, $\{0,1,2,3,22\}$ and $\{17,18,19,20,21\}$.

\vspace{.25cm}
\noindent
Sets $S_{\Gamma_{6},3}$ and $S_{\Gamma_{6},4}$ are given below.  Neither set has elements containing vertices of degree $5$. All elements of $S_{\Gamma_{6},4}$ contain the two vertices of degree $3$ and at least one vertex of degree $4$ that is the child of a sink node.
\begin{eqnarray*}
S_{\Gamma_{6},3} &=&\{(0,10,17), (2,4,18), (2,8,13), (3,5,13), (3,5,17),\\&& (3,9,13), (3,13,17), (5,13,17), (6,9,18), (7,8,13), \\&&(7,8,17), (7,13,17), (8,13,17), \}\\
S_{\Gamma_{6},4}&= &\{(3,5,13,17), (7,8,13,17)\}
\end{eqnarray*}

\section*{Appendix E}\label{appendix24}
The extremal graph $(24,54)$. Here $(deg_{4},deg_{5}) = (12,12)$

\vspace{.25cm}
\noindent
$\Gamma_{0}$ has incidence array:

\vspace{.25cm}
\noindent
\begin{tiny}
$
\begin{array}{cccccccccccccccccccccccccc}
0&0&0&0&0&0&0&0&0&0&0&1&0&0&1&0&0&0&0&1&0&0&0&1\\
0&0&0&0&0&1&0&0&1&0&0&0&0&1&0&0&0&0&0&1&0&0&0&0\\
0&0&0&0&1&0&0&1&0&0&1&0&0&0&0&0&0&0&0&1&0&0&1&0\\
0&0&0&0&0&0&0&0&0&0&1&0&0&1&0&0&0&0&1&0&0&1&0&1\\
0&0&1&0&0&0&0&0&0&0&0&1&1&0&0&0&0&0&1&0&0&0&0&0\\
0&1&0&0&0&0&1&0&0&1&0&0&0&0&1&0&0&0&1&0&0&0&0&0\\
0&0&0&0&0&1&0&0&0&0&0&1&0&0&0&0&0&1&0&0&0&1&1&0\\
0&0&1&0&0&0&0&0&0&1&0&0&0&1&0&0&0&1&0&0&0&0&0&0\\
0&1&0&0&0&0&0&0&0&0&1&0&1&0&0&0&0&1&0&0&0&0&0&0\\
0&0&0&0&0&1&0&1&0&0&0&0&1&0&0&0&1&0&0&0&0&0&0&1\\
0&0&1&1&0&0&0&0&1&0&0&0&0&0&1&0&1&0&0&0&0&0&0&0\\
1&0&0&0&1&0&1&0&0&0&0&0&0&1&0&0&1&0&0&0&0&0&0&0\\
0&0&0&0&1&0&0&0&1&1&0&0&0&0&0&1&0&0&0&0&0&1&0&0\\
0&1&0&1&0&0&0&1&0&0&0&1&0&0&0&1&0&0&0&0&0&0&0&0\\
1&0&0&0&0&1&0&0&0&0&1&0&0&0&0&1&0&0&0&0&0&0&0&0\\
0&0&0&0&0&0&0&0&0&0&0&0&1&1&1&0&0&0&0&0&1&0&1&0\\
0&0&0&0&0&0&0&0&0&1&1&1&0&0&0&0&0&0&0&0&1&0&0&0\\
0&0&0&0&0&0&1&1&1&0&0&0&0&0&0&0&0&0&0&0&1&0&0&0\\
0&0&0&1&1&1&0&0&0&0&0&0&0&0&0&0&0&0&0&0&1&0&0&0\\
1&1&1&0&0&0&0&0&0&0&0&0&0&0&0&0&0&0&0&0&1&1&0&0\\
0&0&0&0&0&0&0&0&0&0&0&0&0&0&0&1&1&1&1&1&0&0&0&0\\
0&0&0&1&0&0&1&0&0&0&0&0&1&0&0&0&0&0&0&1&0&0&0&0\\
0&0&1&0&0&0&1&0&0&0&0&0&0&0&0&1&0&0&0&0&0&0&0&1\\
1&0&0&1&0&0&0&0&0&1&0&0&0&0&0&0&0&0&0&0&0&0&1&0\\
\end{array}
$
\end{tiny}

\vspace{.25cm}
\noindent
and degree sequence: $$(4,4,5,5,4,5,5,4,4,5,5,5,5,5,4,5,4,4,4,5,5,4,4,4)$$ Edge set is:
\begin{eqnarray*}
(0,11), (0,14), (0,19), (0,23), (1,5), (1,8), (1,13), (1,19), (2,4), 
(2,7), \\(2,10), (2,19), (2,22), (3,10), (3,13), (3,18), (3,21), (3,23), (4,11), 
(4,12), \\(4,18), (5,6), (5,9), (5,14), (5,18), (6,11), (6,17), (6,21), (6,22), 
(7,9),\\ (7,13), (7,17), (8,10), (8,12), (8,17), (9,12), (9,16), (9,23), (10,14),\\ 
(10,16),(11,13), (11,16), (12,15), (12,21), (13,15), (14,15), (15,20),\\ (15,22), (16,20), 
(17,20), (18,20), (19,20), (19,21), (22,23), 
\end{eqnarray*}

\vspace{.25cm}
\noindent
Set $S_{\Gamma_{0},3}$ is given below. Elements of $S_{\Gamma_{0},3}$ contain only vertices of degree $4$.
\begin{eqnarray*}
S_{\Gamma_{0},3} &=&\{(0,7,18), (0,8,18), (1,4,23), (1,16,22), (4,14,17), \\
&&(4,17,23), (7,14,21), (8,18,22)\}\\
\end{eqnarray*}

\vspace{.25cm}
\noindent
The subgraph on the vertices of degree $4$,  $\Gamma_{4}$,  has vertex set $$V_{4}=\{0,1,4,7,8,14,16,17,18,21,22,23\}$$ and edge set 
$E_{4}=\{(0,14), (0,23), (1,8), (4,18), (7,17), (8,17), (22,23)\}$.

\vspace{.25cm}
\noindent
 The only set of $4$ non-intersecting elements of $S_{\Gamma_{0},3}$ is $$X=\{(0,8,18),(1,16,22),(4,17,23),(7,14,21)\}$$  The only pair $X_{1},X_{2}\in X$ such that for all $x_{1}\in X_{1}$, $x_{2}\in X_{2}$ there is no edge $(x_{1},x_{2}) \in \Gamma$ is $X_{1}=(1,16,22)$ and $X_{2}=(7,14,21)$. 

\vspace{.25cm}
\noindent
There are $10$ sets of $3$ non-intersecting elements of  $S_{\Gamma,3}$. These are shown in Table \ref{missedTable}. In each case, at least $2$ of the vertices of degree $4$ not contained in any set are in an edge with a vertex that is in one of the sets. 

\begin{table}
\begin{tabular}{|cc|}
\hline
Non-intersecting sets & vertices missed\\
\hline
$(0,7,18)$, $(1,16,22)$, $(4,14,17)$ &  $8$, $21$, $23$\\ 
$(0,7,18)$, $(1,16,22)$, $(4,17,23)$ &  $8$, $14$, $21$\\ 
$(0,8,18)$, $(1,4,23)$, $(7,14,21)$ &  $16$, $17$, $22$\\ 
$(0,8,18)$, $(1,16,22)$, $(4,14,17)$ &  $7$, $21$, $23 $\\
$(0,8,18)$, $(1,16,22)$, $(4,17,23)$ &  $7$, $14$, $21$\\ 
$(0,8,18)$, $(1,16,22)$, $(7,14,21)$ &  $4$, $17$, $23$\\ 
$(0,8,18)$, $(4,17,23)$, $(7,14,21)$ &  $1$, $16$, $22$\\ 
$(1,4,23)$, $(7,14,21)$, $(8,18,22)$ &  $0$, $16$, $17$\\ 
$(1,16,22)$, $(4,17,23)$, $(7,14,21)$ &  $0$, $8$, $18$\\ 
$(4,17,23)$, $(7,14,21)$, $(8,18,22)$ &  $0$, $1$, $16$\\ 
\hline
\end{tabular}
\caption{Sets of $3$ non-intersecting elements of $S_{\Gamma,3}$ and the points missed by those elements \label{missedTable}}
\end{table}

\section*{Appendix F}\label{appendix25}
The $6$ extremal graphs on $25$ vertices. In each case there are $57$ edges. Graphs $\Gamma_{0},\ldots, \Gamma_{2}$ have $(deg_{3},deg_{4},deg_{5})=(1,9,15)$ and graphs $\Gamma_{3},\ldots,\Gamma_{5}$ have $(deg_{4},deg_{5})=(11,14)$. 

\vspace{.25cm}
\noindent
$\Gamma_{0}$ has incidence array:

\vspace{.25cm}
\noindent
\begin{tiny}
$
\begin{array}{ccccccccccccccccccccccccccccc}
0&0&0&0&0&0&0&0&0&0&0&1&0&0&0&0&0&1&0&0&0&0&1&0&1\\
0&0&0&0&0&0&0&1&0&0&1&0&0&0&1&0&0&0&0&0&0&0&1&0&0\\
0&0&0&0&0&0&1&0&0&1&0&0&0&1&0&0&0&0&0&0&0&0&1&0&0\\
0&0&0&0&0&1&0&0&1&0&0&0&1&0&0&0&1&0&0&0&0&0&1&0&0\\
0&0&0&0&0&0&0&0&0&0&0&0&0&0&1&0&1&0&0&0&0&1&0&0&0\\
0&0&0&1&0&0&0&0&0&1&0&0&0&0&0&0&0&0&0&0&0&1&0&0&1\\
0&0&1&0&0&0&0&0&0&0&0&1&1&0&0&1&0&0&0&0&0&1&0&0&0\\
0&1&0&0&0&0&0&0&1&0&0&0&0&1&0&0&0&1&0&0&0&1&0&0&0\\
0&0&0&1&0&0&0&1&0&0&0&0&0&0&0&1&0&0&0&0&1&0&0&0&0\\
0&0&1&0&0&1&0&0&0&0&0&0&0&0&1&0&0&1&0&0&1&0&0&0&0\\
0&1&0&0&0&0&0&0&0&0&0&0&1&0&0&0&0&0&0&0&1&0&0&0&1\\
1&0&0&0&0&0&1&0&0&0&0&0&0&0&0&0&1&0&0&0&1&0&0&0&0\\
0&0&0&1&0&0&1&0&0&0&1&0&0&0&0&0&0&1&0&1&0&0&0&0&0\\
0&0&1&0&0&0&0&1&0&0&0&0&0&0&0&0&1&0&0&1&0&0&0&0&1\\
0&1&0&0&1&0&0&0&0&1&0&0&0&0&0&1&0&0&0&1&0&0&0&0&0\\
0&0&0&0&0&0&1&0&1&0&0&0&0&0&1&0&0&0&1&0&0&0&0&0&1\\
0&0&0&1&1&0&0&0&0&0&0&1&0&1&0&0&0&0&1&0&0&0&0&0&0\\
1&0&0&0&0&0&0&1&0&1&0&0&1&0&0&0&0&0&1&0&0&0&0&0&0\\
0&0&0&0&0&0&0&0&0&0&0&0&0&0&0&1&1&1&0&0&0&0&0&1&0\\
0&0&0&0&0&0&0&0&0&0&0&0&1&1&1&0&0&0&0&0&0&0&0&1&0\\
0&0&0&0&0&0&0&0&1&1&1&1&0&0&0&0&0&0&0&0&0&0&0&1&0\\
0&0&0&0&1&1&1&1&0&0&0&0&0&0&0&0&0&0&0&0&0&0&0&1&0\\
1&1&1&1&0&0&0&0&0&0&0&0&0&0&0&0&0&0&0&0&0&0&0&1&0\\
0&0&0&0&0&0&0&0&0&0&0&0&0&0&0&0&0&0&1&1&1&1&1&0&0\\
1&0&0&0&0&1&0&0&0&0&1&0&0&1&0&1&0&0&0&0&0&0&0&0&0\\
\end{array}
$
\end{tiny}

\vspace{.25cm}
\noindent
and degree sequence: $$(4,4,4,5,3,4,5,5,4,5,4,4,5,5,5,5,5,5,4,4,5,5,5,5,5)$$ Edge set is:
\begin{eqnarray*}
(0,11), (0,17), (0,22), (0,24), (1,7), (1,10), (1,14), (1,22), (2,6), 
(2,9), \\
(2,13), (2,22), (3,5), (3,8), (3,12), (3,16), (3,22), (4,14), (4,16), 
(4,21),\\
 (5,9), (5,21), (5,24), (6,11), (6,12), (6,15), (6,21), (7,8), (7,13), 
(7,17),\\
 (7,21), (8,15), (8,20), (9,14), (9,17), (9,20), (10,12), (10,20),
 (10,24), \\
(11,16), (11,20), (12,17), (12,19), (13,16), (13,19), (13,24), (14,15), (14,19),\\ (15,18), 
(15,24), (16,18), (17,18), (18,23), (19,23), (20,23), (21,23), (22,23) 
\end{eqnarray*}

\vspace{.25cm}
\noindent
Sets $S_{\Gamma_{0},3}$ and $S_{\Gamma_{0},4}$ are given below.  Neither set has elements containing vertices of degree $5$.
\begin{eqnarray*}
S_{\Gamma_{0},3} &=&\{(0,4,8), (0,8,19), (1,5,11), (1,5,18), (2,4,8), \\&&(2,4,10), (2,10,18), (5,11,19)\}\\
S_{\Gamma_{0},4}&= &\{\}
\end{eqnarray*}

\vspace{.25cm}
\noindent
Note that the only sets of $3$ parallel elements from $S_{\Gamma_{0},3}$ which contain the vertex of degree $3$ (i.e. $4$) in the first set and for which neither of the other sets contain a neighbour of the vertex of degree $3$ (i.e. $14$, $16$ or $21$) are: are:
$(0,4,8)$, $(1,5,11)$, $(2,10,18)$; 
$(0,4,8)$, $(2,10,18)$, $(5,11,19)$; 
$(2,4,10)$, $(0,8,19)$, $(1,5,11)$ and 
$(2,4,10)$, $(0,8,19)$, $(1,5,18)$. In all cases there is an edge from the first set to one of the other sets (via edge $(0,11)$ or $(1,10)$).

\vspace{.25cm}
\noindent
$\Gamma_{1}$ has incidence array:

\vspace{.25cm}
\noindent
\begin{tiny}
$
\begin{array}{ccccccccccccccccccccccccccccc}
0&0&0&0&0&0&0&1&0&0&0&1&0&0&0&0&0&1&0&0&0&0&1&0&0\\
0&0&0&0&0&0&1&0&0&0&1&0&0&0&0&0&1&0&0&0&0&0&1&0&0\\
0&0&0&0&0&1&0&0&0&1&0&0&0&0&1&1&0&0&0&0&0&0&1&0&0\\
0&0&0&0&1&0&0&0&0&0&0&0&0&1&0&0&0&0&0&0&0&0&1&0&1\\
0&0&0&1&0&0&0&0&0&0&0&1&0&0&1&0&1&0&0&0&0&1&0&0&0\\
0&0&1&0&0&0&0&0&0&0&1&0&1&0&0&0&0&0&0&0&0&1&0&0&0\\
0&1&0&0&0&0&0&0&1&0&0&0&0&0&0&1&0&0&0&0&0&1&0&0&1\\
1&0&0&0&0&0&0&0&0&1&0&0&0&1&0&0&0&0&0&0&0&1&0&0&0\\
0&0&0&0&0&0&1&0&0&0&0&0&0&0&1&0&0&0&0&0&1&0&0&0&0\\
0&0&1&0&0&0&0&1&0&0&0&0&0&0&0&0&1&0&0&0&1&0&0&0&1\\
0&1&0&0&0&1&0&0&0&0&0&0&0&1&0&0&0&1&0&0&1&0&0&0&0\\
1&0&0&0&1&0&0&0&0&0&0&0&1&0&0&1&0&0&0&0&1&0&0&0&0\\
0&0&0&0&0&1&0&0&0&0&0&1&0&0&0&0&0&0&0&1&0&0&0&0&1\\
0&0&0&1&0&0&0&1&0&0&1&0&0&0&0&1&0&0&0&1&0&0&0&0&0\\
0&0&1&0&1&0&0&0&1&0&0&0&0&0&0&0&0&1&0&1&0&0&0&0&0\\
0&0&1&0&0&0&1&0&0&0&0&1&0&1&0&0&0&0&1&0&0&0&0&0&0\\
0&1&0&0&1&0&0&0&0&1&0&0&0&0&0&0&0&0&1&0&0&0&0&0&0\\
1&0&0&0&0&0&0&0&0&0&1&0&0&0&1&0&0&0&1&0&0&0&0&0&1\\
0&0&0&0&0&0&0&0&0&0&0&0&0&0&0&1&1&1&0&0&0&0&0&1&0\\
0&0&0&0&0&0&0&0&0&0&0&0&1&1&1&0&0&0&0&0&0&0&0&1&0\\
0&0&0&0&0&0&0&0&1&1&1&1&0&0&0&0&0&0&0&0&0&0&0&1&0\\
0&0&0&0&1&1&1&1&0&0&0&0&0&0&0&0&0&0&0&0&0&0&0&1&0\\
1&1&1&1&0&0&0&0&0&0&0&0&0&0&0&0&0&0&0&0&0&0&0&1&0\\
0&0&0&0&0&0&0&0&0&0&0&0&0&0&0&0&0&0&1&1&1&1&1&0&0\\
0&0&0&1&0&0&1&0&0&1&0&0&1&0&0&0&0&1&0&0&0&0&0&0&0\\
\end{array}
$
\end{tiny}

\vspace{.25cm}
\noindent
and degree sequence: $$(4,4,5,4,5,4,5,4,3,5,5,5,4,5,5,5,4,5,4,4,5,5,5,5,5)$$ Edge set is:
\begin{eqnarray*}
(0,7), (0,11), (0,17), (0,22), (1,6), (1,10), (1,16), (1,22), (2,5), 
(2,9), \\
(2,14), (2,15), (2,22), (3,4), (3,13), (3,22), (3,24), (4,11), (4,14), 
(4,16),\\
 (4,21), (5,10), (5,12), (5,21), (6,8), (6,15), (6,21), (6,24), (7,9), 
(7,13), \\
(7,21), (8,14), (8,20), (9,16), (9,20), (9,24), (10,13), (10,17), (10,20),\\ 
(11,12), (11,15), (11,20), (12,19), (12,24), (13,15), (13,19), (14,17), (14,19), \\(15,18), 
(16,18), (17,18), (17,24), (18,23), (19,23), (20,23), (21,23), (22,23) 
\end{eqnarray*}
%

\vspace{.25cm}
\noindent
The sink nodes are vertices $2$ and $15$ and $22$, which have sets of children $\{5,9,14,15,22\}$ and $\{2,6,11,13,18\}$.

\vspace{.25cm}
\noindent
Sets $S_{\Gamma_{1},3}$ and $S_{\Gamma_{1},4}$ are given below. Vertices of degree $5$ are written in bold. Every element of $S_{\Gamma_{1},4}$ contains the vertex of degre $3$ ($8$), a child of a sink node, and two other vertices of degree $4$.
\begin{eqnarray*}
S_{\Gamma_{1},3} &=&\{(0,5,8), (0,5,16), (0,{\bf 6},19), (0,8,16), (0,16,19), (1,7,12), (1,7,{\bf 14}),\\ 
&&(3,5,8), (3,5,18), (3,8,18), (5,8,16), (5,8,18), (7,8,12), (7,8,18), \\&&(7,12,18), (8,12,16), (8,12,18), (8,12,{\bf 22}), (8,{\bf 13},16)\}\\
S_{\Gamma_{1},4}&= &\{(0,5,8,16)(3,5,8,18)(7,8,12,18)\}
\end{eqnarray*}

\vspace{.25cm}
\noindent
Note that the only sets of $3$ parallel elements from $S_{\Gamma_{0},3}$ which contain the vertex of degree $3$ (i.e. $8$)  in the first set and for which neither of the other sets contain a neighbour of the vertex of degree $3$ (i.e. $6$, $14$ or $20$) are:
$(3,8,18)$, $(0,5,16)$, $(1,7,12)$; $(0,8,16)$, $(1,7,12)$, $(3,5,18)$; 
$(3,5,8)$, $(0,16,19)$, $(1,7,12)$; 
$(3,8,18)$, $(0,16,19)$, $(1,7,12)$;
$(5,8,18)$, $(0,16,19)$, $(1,7,12)$;  
$(3,5,8)$, $(0,16,19)$, $(7,12,18)$; 
$(7,8,12)$, $(0,16,19)$, $(3,5,18)$; 
$(8,12,22)$, $(0,16,19)$, $(3,5,18)$; and 
$(8,13,16)$, $(1,7,12)$, $(3,5,18)$. 

 In all cases there is an edge from the first set to one of the other sets (via edge $(16,18)$ or $(5,12)$).

\vspace{.25cm}
\noindent
$\Gamma_{2}$ has incidence array:

\vspace{.25cm}
\noindent
\begin{tiny}
$
\begin{array}{ccccccccccccccccccccccccccccc}
0&0&0&0&0&0&0&1&0&0&0&1&0&0&0&0&0&1&0&0&0&0&1&0&0\\
0&0&0&0&0&0&1&0&0&0&0&0&0&0&1&0&0&0&0&0&0&0&1&0&1\\
0&0&0&0&0&1&0&0&0&0&1&0&0&1&0&0&1&0&0&0&0&0&1&0&0\\
0&0&0&0&1&0&0&0&0&1&0&0&0&0&0&1&0&0&0&0&0&0&1&0&0\\
0&0&0&1&0&0&0&0&0&0&1&0&0&0&0&0&0&0&0&0&0&1&0&0&1\\
0&0&1&0&0&0&0&0&0&0&0&0&1&0&0&1&0&0&0&0&0&1&0&0&0\\
0&1&0&0&0&0&0&0&0&1&0&0&0&1&0&0&0&1&0&0&0&1&0&0&0\\
1&0&0&0&0&0&0&0&1&0&0&0&0&0&1&0&1&0&0&0&0&1&0&0&0\\
0&0&0&0&0&0&0&1&0&0&0&0&1&0&0&0&0&0&0&0&1&0&0&0&0\\
0&0&0&1&0&0&1&0&0&0&0&0&0&0&0&0&1&0&0&0&1&0&0&0&0\\
0&0&1&0&1&0&0&0&0&0&0&0&0&0&1&0&0&1&0&0&1&0&0&0&0\\
1&0&0&0&0&0&0&0&0&0&0&0&0&1&0&1&0&0&0&0&1&0&0&0&1\\
0&0&0&0&0&1&0&0&1&0&0&0&0&0&0&0&0&1&0&1&0&0&0&0&1\\
0&0&1&0&0&0&1&0&0&0&0&1&0&0&0&0&0&0&0&1&0&0&0&0&0\\
0&1&0&0&0&0&0&1&0&0&1&0&0&0&0&1&0&0&0&1&0&0&0&0&0\\
0&0&0&1&0&1&0&0&0&0&0&1&0&0&1&0&0&0&1&0&0&0&0&0&0\\
0&0&1&0&0&0&0&1&0&1&0&0&0&0&0&0&0&0&1&0&0&0&0&0&1\\
1&0&0&0&0&0&1&0&0&0&1&0&1&0&0&0&0&0&1&0&0&0&0&0&0\\
0&0&0&0&0&0&0&0&0&0&0&0&0&0&0&1&1&1&0&0&0&0&0&1&0\\
0&0&0&0&0&0&0&0&0&0&0&0&1&1&1&0&0&0&0&0&0&0&0&1&0\\
0&0&0&0&0&0&0&0&1&1&1&1&0&0&0&0&0&0&0&0&0&0&0&1&0\\
0&0&0&0&1&1&1&1&0&0&0&0&0&0&0&0&0&0&0&0&0&0&0&1&0\\
1&1&1&1&0&0&0&0&0&0&0&0&0&0&0&0&0&0&0&0&0&0&0&1&0\\
0&0&0&0&0&0&0&0&0&0&0&0&0&0&0&0&0&0&1&1&1&1&1&0&0\\
0&1&0&0&1&0&0&0&0&0&0&1&1&0&0&0&1&0&0&0&0&0&0&0&0\\
\end{array}
$
\end{tiny}

\vspace{.25cm}
\noindent
and degree sequence: $$(4,4,5,4,4,4,5,5,3,4,5,5,5,4,5,5,5,5,4,4,5,5,5,5,5)$$ Edge set is:
\begin{eqnarray*}
(0,7), (0,11), (0,17), (0,22), (1,6), (1,14), (1,22), (1,24), (2,5), 
(2,10),\\ (2,13), (2,16), (2,22), (3,4), (3,9), (3,15), (3,22), (4,10), (4,21), 
(4,24), \\(5,12), (5,15), (5,21), (6,9), (6,13), (6,17), (6,21), (7,8), (7,14), 
(7,16),\\ (7,21), (8,12), (8,20), (9,16), (9,20), (10,14), (10,17), (10,20), (11,13), \\
(11,15), (11,20), (11,24), (12,17), (12,19), (12,24), (13,19), (14,15), (14,19), \\(15,18), 
(16,18), (16,24), (17,18), (18,23), (19,23), (20,23), (21,23), (22,23)
\end{eqnarray*}
%

\vspace{.25cm}
\noindent
The sink node is vertex $10$ which has set of children $\{2,4,14,17,20\}$.

\vspace{.25cm}
\noindent
Sets $S_{\Gamma_{2},3}$ and $S_{\Gamma_{2},4}$ are given below. Vertices of degree $5$ are written in bold. Every element of $S_{\Gamma_{2},4}$ contains the vertex of degre $3$ ($8$), a child of a sink node, and two other vertices of degree $4$.
\begin{eqnarray*}
S_{\Gamma_{2},3} &=&\{(0,4,19), (0,5,9), (0,9,19), (1,5,{\bf 20}), (1,8,18), (3,{\bf 7},13),\\&& (3,8,13), (4,8,13), (4,8,18), (4,13,18), (6,8,{\bf 15}), (8,13,18), \}\\
S_{\Gamma_{2},4}&= &\{(4,8,13,18)\}
\end{eqnarray*}

\vspace{.25cm}
\noindent
Note that the only sets of $3$ parallel elements from $S_{\Gamma_{0},3}$ which contain the vertex of degree $3$ (i.e. $8$)  in the first set and for which neither of the other sets contain a neighbour of the vertex of degree $3$ (i.e. $7$, $12$ or $20$) are:
$(6,8,15)$, $(0,5,9)$, $(4,13,18)$ and
$(6,8,15)$, $(0,9,19)$, $(4,13,18)$.   In both cases there is an edge from the first set to one of the other sets (via edge $(6,13)$.

\vspace{.25cm}
\noindent
$\Gamma_{3}$ has incidence array:

\vspace{.25cm}
\noindent
\begin{tiny}
$
\begin{array}{ccccccccccccccccccccccccccccc}
0&0&0&0&0&0&0&0&0&0&0&1&0&0&0&0&0&1&0&0&0&0&1&0&1\\
0&0&0&0&0&0&0&1&0&0&1&0&0&0&1&0&1&0&0&0&0&0&1&0&0\\
0&0&0&0&0&0&1&0&0&1&0&0&0&0&0&1&0&0&0&0&0&0&1&0&0\\
0&0&0&0&0&1&0&0&1&0&0&0&0&1&0&0&0&0&0&0&0&0&1&0&0\\
0&0&0&0&0&0&0&0&0&0&0&1&0&0&1&1&0&0&0&0&0&1&0&0&0\\
0&0&0&1&0&0&0&0&0&0&1&0&1&0&0&0&0&1&0&0&0&1&0&0&0\\
0&0&1&0&0&0&0&0&1&0&0&0&0&0&0&0&1&0&0&0&0&1&0&0&0\\
0&1&0&0&0&0&0&0&0&1&0&0&0&1&0&0&0&0&0&0&0&1&0&0&1\\
0&0&0&1&0&0&1&0&0&0&0&0&0&0&1&0&0&0&0&0&1&0&0&0&1\\
0&0&1&0&0&0&0&1&0&0&0&0&0&0&0&0&0&1&0&0&1&0&0&0&0\\
0&1&0&0&0&1&0&0&0&0&0&0&0&0&0&1&0&0&0&0&1&0&0&0&0\\
1&0&0&0&1&0&0&0&0&0&0&0&0&1&0&0&1&0&0&0&1&0&0&0&0\\
0&0&0&0&0&1&0&0&0&0&0&0&0&0&0&0&1&0&0&1&0&0&0&0&1\\
0&0&0&1&0&0&0&1&0&0&0&1&0&0&0&0&0&0&0&1&0&0&0&0&0\\
0&1&0&0&1&0&0&0&1&0&0&0&0&0&0&0&0&1&0&1&0&0&0&0&0\\
0&0&1&0&1&0&0&0&0&0&1&0&0&0&0&0&0&0&1&0&0&0&0&0&1\\
0&1&0&0&0&0&1&0&0&0&0&1&1&0&0&0&0&0&1&0&0&0&0&0&0\\
1&0&0&0&0&1&0&0&0&1&0&0&0&0&1&0&0&0&1&0&0&0&0&0&0\\
0&0&0&0&0&0&0&0&0&0&0&0&0&0&0&1&1&1&0&0&0&0&0&1&0\\
0&0&0&0&0&0&0&0&0&0&0&0&1&1&1&0&0&0&0&0&0&0&0&1&0\\
0&0&0&0&0&0&0&0&1&1&1&1&0&0&0&0&0&0&0&0&0&0&0&1&0\\
0&0&0&0&1&1&1&1&0&0&0&0&0&0&0&0&0&0&0&0&0&0&0&1&0\\
1&1&1&1&0&0&0&0&0&0&0&0&0&0&0&0&0&0&0&0&0&0&0&1&0\\
0&0&0&0&0&0&0&0&0&0&0&0&0&0&0&0&0&0&1&1&1&1&1&0&0\\
1&0&0&0&0&0&0&1&1&0&0&0&1&0&0&1&0&0&0&0&0&0&0&0&0\\
\end{array}
$
\end{tiny}

\vspace{.25cm}
\noindent
and degree sequence: $$(4,5,4,4,4,5,4,5,5,4,4,5,4,4,5,5,5,5,4,4,5,5,5,5,5)$$ Edge set is:
\begin{eqnarray*}
(0,11), (0,17), (0,22), (0,24), (1,7), (1,10), (1,14), (1,16), (1,22), 
(2,6),\\
 (2,9), (2,15), (2,22), (3,5), (3,8), (3,13), (3,22), (4,11), (4,14), 
(4,15),\\
 (4,21), (5,10), (5,12), (5,17), (5,21), (6,8), (6,16), (6,21), (7,9), 
(7,13),\\
 (7,21), (7,24), (8,14), (8,20), (8,24), (9,17), (9,20), (10,15), (10,20), \\
(11,13), (11,16), (11,20), (12,16), (12,19), (12,24), (13,19), (14,17), (14,19),\\
 (15,18), (15,24), (16,18), (17,18), (18,23), (19,23), (20,23), (21,23), (22,23)
\end{eqnarray*}
%

\vspace{.25cm}
\noindent
The sink node is vertex $1$ which has set of children $\{7,10,14,16,22\}$

\vspace{.25cm}
\noindent
Sets $S_{\Gamma_{3},3}$ and $S_{\Gamma_{3},4}$ are given below. Vertices of degree $5$ are written in bold. Every element of $S_{\Gamma_{3},4}$ contains a child of  the sink node, and three other vertices of degree $4$.
\begin{eqnarray*}
S_{\Gamma_{3},3} &=&\{ (0,6,10), (0,6,19), (0,10,19), (2,5,{\bf 11}), (3,4,9), (3,9,{\bf 16}),\\&&
 (4,9,12), (4,12,{\bf 22}), (6,10,13), (6,10,19), (6,13,{\bf 17})\}\\
 S_{\Gamma_{3},4} &=&\{(0,6,10,19)\}
\end{eqnarray*}

\vspace{.25cm}
\noindent
Note that 
\begin{enumerate}
\item $S_{\Gamma_{3},3}$ does not contain $3$ non-intersecting sets consisting entirely of vertices of degree $4$. 
\item The only pairs of elements of  $S_{\Gamma_{3},3}$ that are parallel and have no edges between them are:
\begin{eqnarray*}
&&
((0,6,10),(3,4,9)),\\
&& ((0,6,10),(4,9,12)), \\
&&((0,6,19),(3,4,9)), \\
&&((0,10,19),(3,4,9)), \\
&&((0,10,19), (3,9,{\bf 16})),\\
&&((3,4,9),(6,10,19)), \\
&&((4,9,12),(6,10,13)),\\
&&((4,12,{\bf 22}),(6,10,13))\; {\rm and} \\
&&((4,12,{\bf 22}),(6,13,17)). 
\end{eqnarray*}
\end{enumerate}
If there is a set $X_{1}$, $X_{2}$, $X_{3}$ of mutually parallel elements of  $S_{\Gamma_{3},3}$ where there are no edges from $X_{1}$ to $X_{2}$ or $X_{3}$, then there would be two pairs of triples from the list above which intersect in one triple, and where the other two triples are distinct. This is not the case.

\vspace{.25cm}
\noindent
$\Gamma_{4}$ has incidence array:

\vspace{.25cm}
\noindent
\begin{tiny}
$
\begin{array}{ccccccccccccccccccccccccccccc}
0&0&0&0&0&0&0&1&0&0&0&1&0&0&0&0&0&1&0&0&0&0&1&0&0\\
0&0&0&0&0&0&1&0&0&0&0&0&0&0&1&0&1&0&0&0&0&0&1&0&0\\
0&0&0&0&0&1&0&0&0&0&1&0&0&1&0&0&0&0&0&0&0&0&1&0&1\\
0&0&0&0&1&0&0&0&0&1&0&0&1&0&0&1&0&0&0&0&0&0&1&0&0\\
0&0&0&1&0&0&0&0&1&0&0&0&0&1&0&0&0&1&0&0&0&1&0&0&0\\
0&0&1&0&0&0&0&0&0&1&0&0&0&0&0&0&1&0&0&0&0&1&0&0&0\\
0&1&0&0&0&0&0&0&0&0&1&0&0&0&0&1&0&0&0&0&0&1&0&0&0\\
1&0&0&0&0&0&0&0&0&0&0&0&1&0&0&0&0&0&0&0&0&1&0&0&1\\
0&0&0&0&1&0&0&0&0&0&0&0&0&0&0&0&1&0&0&0&1&0&0&0&1\\
0&0&0&1&0&1&0&0&0&0&0&0&0&0&1&0&0&0&0&0&1&0&0&0&0\\
0&0&1&0&0&0&1&0&0&0&0&0&1&0&0&0&0&1&0&0&1&0&0&0&0\\
1&0&0&0&0&0&0&0&0&0&0&0&0&1&0&1&0&0&0&0&1&0&0&0&0\\
0&0&0&1&0&0&0&1&0&0&1&0&0&0&0&0&1&0&0&1&0&0&0&0&0\\
0&0&1&0&1&0&0&0&0&0&0&1&0&0&0&0&0&0&0&1&0&0&0&0&0\\
0&1&0&0&0&0&0&0&0&1&0&0&0&0&0&0&0&1&0&1&0&0&0&0&1\\
0&0&0&1&0&0&1&0&0&0&0&1&0&0&0&0&0&0&1&0&0&0&0&0&1\\
0&1&0&0&0&1&0&0&1&0&0&0&1&0&0&0&0&0&1&0&0&0&0&0&0\\
1&0&0&0&1&0&0&0&0&0&1&0&0&0&1&0&0&0&1&0&0&0&0&0&0\\
0&0&0&0&0&0&0&0&0&0&0&0&0&0&0&1&1&1&0&0&0&0&0&1&0\\
0&0&0&0&0&0&0&0&0&0&0&0&1&1&1&0&0&0&0&0&0&0&0&1&0\\
0&0&0&0&0&0&0&0&1&1&1&1&0&0&0&0&0&0&0&0&0&0&0&1&0\\
0&0&0&0&1&1&1&1&0&0&0&0&0&0&0&0&0&0&0&0&0&0&0&1&0\\
1&1&1&1&0&0&0&0&0&0&0&0&0&0&0&0&0&0&0&0&0&0&0&1&0\\
0&0&0&0&0&0&0&0&0&0&0&0&0&0&0&0&0&0&1&1&1&1&1&0&0\\
0&0&1&0&0&0&0&1&1&0&0&0&0&0&1&1&0&0&0&0&0&0&0&0&0\\
\end{array}
$
\end{tiny}

\vspace{.25cm}
\noindent
and degree sequence: $$(4,4,5,5,5,4,4,4,4,4,5,4,5,4,5,5,5,5,4,4,5,5,5,5,5)$$ Edge set is:
\begin{eqnarray*}
(0,7), (0,11), (0,17), (0,22), (1,6), (1,14), (1,16), (1,22), (2,5), 
(2,10),\\
 (2,13), (2,22), (2,24), (3,4), (3,9), (3,12), (3,15), (3,22), (4,8), 
(4,13),\\
 (4,17), (4,21), (5,9), (5,16), (5,21), (6,10), (6,15), (6,21), (7,12), 
(7,21),\\
 (7,24), (8,16), (8,20), (8,24), (9,14), (9,20), (10,12), (10,17), (10,20), \\
(11,13), (11,15), (11,20), (12,16), (12,19), (13,19), (14,17), (14,19), (14,24), \\(15,18), 
(15,24), (16,18), (17,18), (18,23), (19,23), (20,23), (21,23), (22,23)
\end{eqnarray*}
%

\vspace{.25cm}
\noindent
The sink nodes are vertices $3$ and $10$ which have sets of children $\{4,9,12,15,22\}$ and $\{2,6,12,17,17\}$.

\vspace{.25cm}
\noindent
Sets $S_{\Gamma_{4},3}$ and $S_{\Gamma_{4},4}$ are given below. Vertices of degree $5$ are written in bold. Every element of $S_{\Gamma_{4},4}$ contains a child of  the sink node, and three other vertices of degree $4$.
\begin{eqnarray*}
S_{\Gamma_{4},3} &=&\{ (0,5,19), (0,6,8), (0,6,9), (0,6,19), (0,8,19), (1,7,13),\\&& (1,7,{\bf 20}), (5,{\bf 15},19), (6,8,19), (6,9,13), (7,9,13), (7,9,18), \\&&(7,13,18), (9,13,18), (11,{\bf 14},{\bf 21}), \}\\
 S_{\Gamma_{4},4} &=&\{(0,6,8,19), (7,9,13,18)\}
\end{eqnarray*}
Note that  
\begin{enumerate}
\item $S_{\Gamma_{4},3}$ does not contain $3$ non-intersecting sets consisting entirely of vertices of degree $4$.\item The only pairs of elements of  $S_{\Gamma_{4},3}$ that are parallel and have no edges between them are:
\newline 
$((0,6,8),(9,13,18))$,  
$((1,7,20),(5,15,19))$ and $((6,8,19),(7,9,18))$. 
\end{enumerate}

If there is a set $X_{1}$, $X_{2}$, $X_{3}$ of mutually parallel elements of  $S_{\Gamma_{3},3}$ where there are no edges from $X_{1}$ to $X_{2}$ or $X_{3}$, then there would be two pairs of triples from the list above which intersect in one triple, and where the other two triples are distinct. This is not the case.

\vspace{.25cm}
\noindent
$\Gamma_{5}$ has incidence array:

\vspace{.25cm}
\noindent
\begin{tiny}
$
\begin{array}{ccccccccccccccccccccccccccccc}
0&0&0&0&0&0&0&1&0&0&0&1&0&0&1&0&0&0&0&0&0&0&1&0&0\\
0&0&0&0&0&0&1&0&0&0&1&0&0&1&0&0&0&1&0&0&0&0&1&0&0\\
0&0&0&0&0&1&0&0&0&0&0&0&1&0&0&0&0&0&0&0&0&0&1&0&1\\
0&0&0&0&1&0&0&0&0&1&0&0&0&0&0&0&1&0&0&0&0&0&1&0&0\\
0&0&0&1&0&0&0&0&0&0&0&1&0&1&0&1&0&0&0&0&0&1&0&0&0\\
0&0&1&0&0&0&0&0&1&0&0&0&0&0&1&0&0&1&0&0&0&1&0&0&0\\
0&1&0&0&0&0&0&0&0&1&0&0&1&0&0&0&0&0&0&0&0&1&0&0&0\\
1&0&0&0&0&0&0&0&0&0&1&0&0&0&0&0&1&0&0&0&0&1&0&0&1\\
0&0&0&0&0&1&0&0&0&0&0&0&0&1&0&0&1&0&0&0&1&0&0&0&0\\
0&0&0&1&0&0&1&0&0&0&0&0&0&0&1&0&0&0&0&0&1&0&0&0&1\\
0&1&0&0&0&0&0&1&0&0&0&0&0&0&0&1&0&0&0&0&1&0&0&0&0\\
1&0&0&0&1&0&0&0&0&0&0&0&1&0&0&0&0&1&0&0&1&0&0&0&0\\
0&0&1&0&0&0&1&0&0&0&0&1&0&0&0&0&1&0&0&1&0&0&0&0&0\\
0&1&0&0&1&0&0&0&1&0&0&0&0&0&0&0&0&0&0&1&0&0&0&0&1\\
1&0&0&0&0&1&0&0&0&1&0&0&0&0&0&1&0&0&0&1&0&0&0&0&0\\
0&0&0&0&1&0&0&0&0&0&1&0&0&0&1&0&0&0&1&0&0&0&0&0&0\\
0&0&0&1&0&0&0&1&1&0&0&0&1&0&0&0&0&0&1&0&0&0&0&0&0\\
0&1&0&0&0&1&0&0&0&0&0&1&0&0&0&0&0&0&1&0&0&0&0&0&0\\
0&0&0&0&0&0&0&0&0&0&0&0&0&0&0&1&1&1&0&0&0&0&0&1&0\\
0&0&0&0&0&0&0&0&0&0&0&0&1&1&1&0&0&0&0&0&0&0&0&1&0\\
0&0&0&0&0&0&0&0&1&1&1&1&0&0&0&0&0&0&0&0&0&0&0&1&0\\
0&0&0&0&1&1&1&1&0&0&0&0&0&0&0&0&0&0&0&0&0&0&0&1&0\\
1&1&1&1&0&0&0&0&0&0&0&0&0&0&0&0&0&0&0&0&0&0&0&1&0\\
0&0&0&0&0&0&0&0&0&0&0&0&0&0&0&0&0&0&1&1&1&1&1&0&0\\
0&0&1&0&0&0&0&1&0&1&0&0&0&1&0&0&0&0&0&0&0&0&0&0&0\\
\end{array}
$
\end{tiny}

\vspace{.25cm}
\noindent
and degree sequence: $$(4,5,4,4,5,5,4,5,4,5,4,5,5,5,5,4,5,4,4,4,5,5,5,5,4)$$ Edge set is:
\begin{eqnarray*}
(0,7), (0,11), (0,14), (0,22), (1,6), (1,10), (1,13), (1,17), (1,22), 
(2,5),\\
 (2,12), (2,22), (2,24), (3,4), (3,9), (3,16), (3,22), (4,11), (4,13), 
(4,15), \\
(4,21), (5,8), (5,14), (5,17), (5,21), (6,9), (6,12), (6,21), (7,10), 
(7,16),\\
 (7,21), (7,24), (8,13), (8,16), (8,20), (9,14), (9,20), (9,24), (10,15), \\
(10,20), (11,12), (11,17), (11,20), (12,16), (12,19), (13,19), (13,24), (14,15),\\
 (14,19), 
(15,18), (16,18), (17,18), (18,23), (19,23), (20,23), (21,23), (22,23)
\end{eqnarray*}
%

\vspace{.25cm}
\noindent
The sink node is vertex $21$ which has set of children $\{4,5,6,7,23\}$.

\vspace{.25cm}
\noindent
Sets $S_{\Gamma_{5},3}$ and $S_{\Gamma_{5},4}$ are given below. Vertices of degree $5$ are written in bold. 
\begin{eqnarray*}
S_{\Gamma_{5},3} &=&\{(0,6,8), (0,6,18), (0,{\bf 13},18), ({\bf 1},{\bf 14},{\bf 16}), (3,{\bf 5},10),\\&& (3,10,19), (3,17,19), (6,8,15), ({\bf 7},17,19), (8,15,{\bf 22})\}\\
 S_{\Gamma_{5},4} &=&\{\}
\end{eqnarray*}
Note that  
\begin{enumerate}
\item $S_{\Gamma_{5},3}$ does not contain $3$ non-intersecting sets consisting entirely of vertices of degree $4$. 
\item The only pairs of elements of  $S_{\Gamma_{3},3}$ that are parallel and have no edges between them are:
\begin{eqnarray*}
&& ((0,6,8,),(3,10,19)),\\
&&((0,6,8,),(3,17,19),\\ 
&&((0,6,18,),(3,5,10),\\
&&((0,6,18,),(3,10,19),\\ 
&&((0,13,18,),(3,5,10),\\
&&((3,17,19,),(6,8,15),\\
&&((6,8,15,),(7,17,19),\\
&&((7,17,19,),(8,15,22).
\end{eqnarray*}
\end{enumerate}

If there is a set $X_{1}$, $X_{2}$, $X_{3}$ of mutually parallel elements of  $S_{\Gamma_{3},3}$ where there are no edges from $X_{1}$ to $X_{2}$ or $X_{3}$, then there would be two pairs of triples from the list above which intersect in one triple, and where the other two triples are distinct. This is not the case.

\section*{Appendix G}\label{appendix26}
The $2$ extremal graphs on $26$ vertices. In each case there are $61$ edges and $(deg_{4},deg_{5}) = (8,18)$.

\vspace{.25cm}
\noindent
$\Gamma_{0}$ has incidence array:

\vspace{.25cm}
\noindent
\begin{tiny}
$
\begin{array}{cccccccccccccccccccccccccccc}
0&0&0&0&0&0&0&1&0&0&0&1&0&0&0&0&0&1&0&0&0&0&1&0&0&1&\\
0&0&0&0&0&0&1&0&0&0&1&0&0&0&0&0&1&0&0&0&0&0&1&0&0&0&\\
0&0&0&0&0&1&0&0&0&1&0&0&0&0&1&1&0&0&0&0&0&0&1&0&0&0&\\
0&0&0&0&1&0&0&0&0&0&0&0&0&1&0&0&0&0&0&0&0&0&1&0&1&0&\\
0&0&0&1&0&0&0&0&0&0&0&1&0&0&1&0&1&0&0&0&0&1&0&0&0&0&\\
0&0&1&0&0&0&0&0&0&0&1&0&1&0&0&0&0&0&0&0&0&1&0&0&0&1&\\
0&1&0&0&0&0&0&0&1&0&0&0&0&0&0&1&0&0&0&0&0&1&0&0&1&0&\\
1&0&0&0&0&0&0&0&0&1&0&0&0&1&0&0&0&0&0&0&0&1&0&0&0&0&\\
0&0&0&0&0&0&1&0&0&0&0&0&0&0&1&0&0&0&0&0&1&0&0&0&0&1&\\
0&0&1&0&0&0&0&1&0&0&0&0&0&0&0&0&1&0&0&0&1&0&0&0&1&0&\\
0&1&0&0&0&1&0&0&0&0&0&0&0&1&0&0&0&1&0&0&1&0&0&0&0&0&\\
1&0&0&0&1&0&0&0&0&0&0&0&1&0&0&1&0&0&0&0&1&0&0&0&0&0&\\
0&0&0&0&0&1&0&0&0&0&0&1&0&0&0&0&0&0&0&1&0&0&0&0&1&0&\\
0&0&0&1&0&0&0&1&0&0&1&0&0&0&0&1&0&0&0&1&0&0&0&0&0&0&\\
0&0&1&0&1&0&0&0&1&0&0&0&0&0&0&0&0&1&0&1&0&0&0&0&0&0&\\
0&0&1&0&0&0&1&0&0&0&0&1&0&1&0&0&0&0&1&0&0&0&0&0&0&0&\\
0&1&0&0&1&0&0&0&0&1&0&0&0&0&0&0&0&0&1&0&0&0&0&0&0&1&\\
1&0&0&0&0&0&0&0&0&0&1&0&0&0&1&0&0&0&1&0&0&0&0&0&1&0&\\
0&0&0&0&0&0&0&0&0&0&0&0&0&0&0&1&1&1&0&0&0&0&0&1&0&0&\\
0&0&0&0&0&0&0&0&0&0&0&0&1&1&1&0&0&0&0&0&0&0&0&1&0&0&\\
0&0&0&0&0&0&0&0&1&1&1&1&0&0&0&0&0&0&0&0&0&0&0&1&0&0&\\
0&0&0&0&1&1&1&1&0&0&0&0&0&0&0&0&0&0&0&0&0&0&0&1&0&0&\\
1&1&1&1&0&0&0&0&0&0&0&0&0&0&0&0&0&0&0&0&0&0&0&1&0&0&\\
0&0&0&0&0&0&0&0&0&0&0&0&0&0&0&0&0&0&1&1&1&1&1&0&0&0&\\
0&0&0&1&0&0&1&0&0&1&0&0&1&0&0&0&0&1&0&0&0&0&0&0&0&0&\\
1&0&0&0&0&1&0&0&1&0&0&0&0&0&0&0&1&0&0&0&0&0&0&0&0&0&\\
\end{array}
$
\end{tiny}

\vspace{.25cm}
\noindent
and degree sequence: $$(5,4,5,4,5,5,5,4,4,5,5,5,4,5,5,5,5,5,4,4,5,5,5,5,5,4)$$ Edge set is:
\begin{eqnarray*}
(0,7), (0,11), (0,17), (0,22), (0,25), (1,6), (1,10), (1,16), (1,22), (2,5), (2,9),\\ 
(2,14), (2,15), (2,22), (3,4), (3,13), (3,22), (3,24), (4,11), (4,14), (4,16), (4,21), \\
(5,10), (5,12), (5,21), (5,25), (6,8), (6,15), (6,21), (6,24), (7,9), (7,13), (7,21), \\
(8,14), (8,20), (8,25), (9,16), (9,20), (9,24), (10,13), (10,17), (10,20), (11,12),\\ 
 (11,15), (11,20), (12,19), (12,24), (13,15), (13,19), (14,17), (14,19),(15,18),\\ 
 (16,18), (16,25), (17,18), (17,24), (18,23), (19,23), (20,23),
(21,23), (22,23)
\end{eqnarray*}

\vspace{.25cm}
\noindent
The sink node is vertex $2$ which has set of children $\{5,9,14,15,22\}$.

\vspace{.25cm}
\noindent
Sets $S_{\Gamma_{0},3}$ and $S_{\Gamma_{0},4}$ are given below. Vertices of degree $5$ are written in bold. Note that no element of $S_{\Gamma_{0},4}$ contains a child of the sink node.
\begin{eqnarray*}
S_{\Gamma_{0},3} &=&\{({\bf 0}, {\bf 6},19), (1,7,12), (1,7, {\bf 14}), (3, {\bf 5},18), (3,8,18), (7,8,12), (7,8,18),\\
&&\; (7,12,18), (8,12,18), (8,12, {\bf 22}), \}\\
S_{\Gamma_{0},4}&= &\{(7,8,12,18)\}
\end{eqnarray*}
The pairs of elements from $S_{\Gamma_{0},3}$ that are parallel and have no edges between them are:
$(({\bf 0}, {\bf 6},19), (3, {\bf 5},18))$, $( (1,7,12), (3, 8,18))$, $((1,7, {\bf 14}), (3, {\bf 5},18))$ and $((1,7, {\bf 14}), ( {\bf 23}, {\bf 24},25))$. Only one of these pairs consists entirely of vertices of degree $4$, namely $( (1,7,12), (3, 8,18))$.

\vspace{.25cm}
\noindent
$\Gamma_{1}$ has incidence array:

\vspace{.25cm}
\noindent
\begin{tiny}
$
\begin{array}{cccccccccccccccccccccccccccc}
0&0&0&0&0&0&0&1&0&0&0&1&0&0&0&0&0&1&0&0&0&0&1&0&0&0\\
0&0&0&0&0&0&1&0&0&0&1&0&0&0&0&0&1&0&0&0&0&0&1&0&0&0\\
0&0&0&0&0&1&0&0&0&1&0&0&0&0&1&1&0&0&0&0&0&0&1&0&0&0\\
0&0&0&0&1&0&0&0&0&0&0&0&0&1&0&0&0&0&0&0&0&0&1&0&1&1\\
0&0&0&1&0&0&0&0&0&0&0&1&0&0&1&0&1&0&0&0&0&1&0&0&0&0\\
0&0&1&0&0&0&0&0&0&0&1&0&1&0&0&0&0&0&0&0&0&1&0&0&0&1\\
0&1&0&0&0&0&0&0&1&0&0&0&0&0&0&1&0&0&0&0&0&1&0&0&1&0\\
1&0&0&0&0&0&0&0&0&1&0&0&0&1&0&0&0&0&0&0&0&1&0&0&0&0\\
0&0&0&0&0&0&1&0&0&0&0&0&0&0&1&0&0&0&0&0&1&0&0&0&0&1\\
0&0&1&0&0&0&0&1&0&0&0&0&0&0&0&0&1&0&0&0&1&0&0&0&1&0\\
0&1&0&0&0&1&0&0&0&0&0&0&0&1&0&0&0&1&0&0&1&0&0&0&0&0\\
1&0&0&0&1&0&0&0&0&0&0&0&1&0&0&1&0&0&0&0&1&0&0&0&0&0\\
0&0&0&0&0&1&0&0&0&0&0&1&0&0&0&0&0&0&0&1&0&0&0&0&1&0\\
0&0&0&1&0&0&0&1&0&0&1&0&0&0&0&1&0&0&0&1&0&0&0&0&0&0\\
0&0&1&0&1&0&0&0&1&0&0&0&0&0&0&0&0&1&0&1&0&0&0&0&0&0\\
0&0&1&0&0&0&1&0&0&0&0&1&0&1&0&0&0&0&1&0&0&0&0&0&0&0\\
0&1&0&0&1&0&0&0&0&1&0&0&0&0&0&0&0&0&1&0&0&0&0&0&0&0\\
1&0&0&0&0&0&0&0&0&0&1&0&0&0&1&0&0&0&1&0&0&0&0&0&1&0\\
0&0&0&0&0&0&0&0&0&0&0&0&0&0&0&1&1&1&0&0&0&0&0&1&0&1\\
0&0&0&0&0&0&0&0&0&0&0&0&1&1&1&0&0&0&0&0&0&0&0&1&0&0\\
0&0&0&0&0&0&0&0&1&1&1&1&0&0&0&0&0&0&0&0&0&0&0&1&0&0\\
0&0&0&0&1&1&1&1&0&0&0&0&0&0&0&0&0&0&0&0&0&0&0&1&0&0\\
1&1&1&1&0&0&0&0&0&0&0&0&0&0&0&0&0&0&0&0&0&0&0&1&0&0\\
0&0&0&0&0&0&0&0&0&0&0&0&0&0&0&0&0&0&1&1&1&1&1&0&0&0\\
0&0&0&1&0&0&1&0&0&1&0&0&1&0&0&0&0&1&0&0&0&0&0&0&0&0\\
0&0&0&1&0&1&0&0&1&0&0&0&0&0&0&0&0&0&1&0&0&0&0&0&0&0
\end{array}
$
\end{tiny}

\vspace{.25cm}
\noindent
and degree sequence: $$(4,4,5,5,5,5,5,4,4,5,5,5,4,5,5,5,4,5,5,4,5,5,5,5,5,4)$$
Edge set is: 
\begin{eqnarray*}
(0,7), (0,11), (0,17), (0,22), (1,6), (1,10), (1,16), (1,22), (2,5), \\ 
(2,9), (2,14), (2,15), (2,22), (3,4), (3,13), (3,22), (3,24), (3,25),\\
 (4,11), (4,14), (4,16), (4,21), (5,10), (5,12), (5,21), (5,25), (6,8),\\
 (6,15), (6,21), (6,24), (7,9), (7,13), (7,21),(8,14), (8,20), (8,25), \\
(9,16), (9,20), (9,24), (10,13), (10,17), (10,20), (11,12), (11,15), (11,20),\\
 (12,19), (12,24), (13,15), (13,19), (14,17), (14,19), (15,18), (16,18),\\
 (17,18),(17,24), (18,23), (18,25), (19,23), (20,23) 
(21,23), (22,23) 
\end{eqnarray*}

\vspace{.25cm}
\noindent
Sets $S_{\Gamma_{1},3}$ and $S_{\Gamma_{1},4}$ are given below. Vertices of degree $5$ are written in bold.
\begin{eqnarray*}
S_{\Gamma_{1},3} &=&\{(0,{\bf 5},16), (0,{\bf 6},19), (0,8,16), (0,16,19), (0,19,25), (1,7,12), \\
&&\; (1,7,{\bf 14}),(1,7,25), (1,{\bf 11},25), (1,19,25), (7,8,12), (7,12,{\bf 18}),\\&&\; (8,12,16), (8,12,{\bf 22}), (8,13,16), ({\bf 9},19,25) \}\\
S_{\Gamma_{1},4}&= &\{\}
\end{eqnarray*}

\vspace{.25cm}
\noindent
Note that $S_{\Gamma_{1},3}$ has no pair of distinct sets that consist of vertices of degree $4$ and  have no edges between them. 

\section*{Appendix H}\label{appendix27}
The $1$ extremal graph on $27$ vertices. There are $65$ edges and $(deg_{4},deg_{5}) = (5,22)$.

\vspace{.25cm}
\noindent
$\Gamma_{0}$ has incidence array:

\vspace{.25cm}
\noindent
\begin{tiny}
$
\begin{array}{ccccccccccccccccccccccccccccc}
0&0&0&0&0&0&0&1&0&0&0&1&0&0&0&0&0&1&0&0&0&0&1&0&0&1&0\\
0&0&0&0&0&0&1&0&0&0&1&0&0&0&0&0&1&0&0&0&0&0&1&0&0&0&0\\
0&0&0&0&0&1&0&0&0&1&0&0&0&0&1&1&0&0&0&0&0&0&1&0&0&0&0\\
0&0&0&0&1&0&0&0&0&0&0&0&0&1&0&0&0&0&0&0&0&0&1&0&1&0&0\\
0&0&0&1&0&0&0&0&0&0&0&1&0&0&1&0&1&0&0&0&0&1&0&0&0&0&0\\
0&0&1&0&0&0&0&0&0&0&1&0&1&0&0&0&0&0&0&0&0&1&0&0&0&1&0\\
0&1&0&0&0&0&0&0&1&0&0&0&0&0&0&1&0&0&0&0&0&1&0&0&1&0&0\\
1&0&0&0&0&0&0&0&0&1&0&0&0&1&0&0&0&0&0&0&0&1&0&0&0&0&1\\
0&0&0&0&0&0&1&0&0&0&0&0&0&0&1&0&0&0&0&0&1&0&0&0&0&1&1\\
0&0&1&0&0&0&0&1&0&0&0&0&0&0&0&0&1&0&0&0&1&0&0&0&1&0&0\\
0&1&0&0&0&1&0&0&0&0&0&0&0&1&0&0&0&1&0&0&1&0&0&0&0&0&0\\
1&0&0&0&1&0&0&0&0&0&0&0&1&0&0&1&0&0&0&0&1&0&0&0&0&0&0\\
0&0&0&0&0&1&0&0&0&0&0&1&0&0&0&0&0&0&0&1&0&0&0&0&1&0&1\\
0&0&0&1&0&0&0&1&0&0&1&0&0&0&0&1&0&0&0&1&0&0&0&0&0&0&0\\
0&0&1&0&1&0&0&0&1&0&0&0&0&0&0&0&0&1&0&1&0&0&0&0&0&0&0\\
0&0&1&0&0&0&1&0&0&0&0&1&0&1&0&0&0&0&1&0&0&0&0&0&0&0&0\\
0&1&0&0&1&0&0&0&0&1&0&0&0&0&0&0&0&0&1&0&0&0&0&0&0&1&0\\
1&0&0&0&0&0&0&0&0&0&1&0&0&0&1&0&0&0&1&0&0&0&0&0&1&0&0\\
0&0&0&0&0&0&0&0&0&0&0&0&0&0&0&1&1&1&0&0&0&0&0&1&0&0&1\\
0&0&0&0&0&0&0&0&0&0&0&0&1&1&1&0&0&0&0&0&0&0&0&1&0&0&0\\
0&0&0&0&0&0&0&0&1&1&1&1&0&0&0&0&0&0&0&0&0&0&0&1&0&0&0\\
0&0&0&0&1&1&1&1&0&0&0&0&0&0&0&0&0&0&0&0&0&0&0&1&0&0&0\\
1&1&1&1&0&0&0&0&0&0&0&0&0&0&0&0&0&0&0&0&0&0&0&1&0&0&0\\
0&0&0&0&0&0&0&0&0&0&0&0&0&0&0&0&0&0&1&1&1&1&1&0&0&0&0\\
0&0&0&1&0&0&1&0&0&1&0&0&1&0&0&0&0&1&0&0&0&0&0&0&0&0&0\\
1&0&0&0&0&1&0&0&1&0&0&0&0&0&0&0&1&0&0&0&0&0&0&0&0&0&0\\
0&0&0&0&0&0&0&1&1&0&0&0&1&0&0&0&0&0&1&0&0&0&0&0&0&0&0
\end{array}
$
\end{tiny}

\vspace{.25cm}
\noindent
and degree sequence: $$(5,4,5,4,5,5,5,5,5,5,5,5,5,5,5,5,5,5,5,4,5,5,5,5,5,4,4)$$
Edge set is: 
\begin{eqnarray*}
(0,7), (0,11), (0,17), (0,22), (0,25), (1,6), (1,10), (1,16), (1,22), (2,5),\\ 
(2,9), (2,14), (2,15), (2,22), (3,4), (3,13), (3,22), (3,24), (4,11), (4,14), \\
(4,16), (4,21), (5,10), (5,12), (5,21), (5,25), (6,8), (6,15), (6,21), 
(6,24), \\
(7,9), (7,13), (7,21), (7,26), (8,14),(8,20), (8,25),
 (8,26), (9,16), (9,20),\\ 
(9,24), (10,13), (10,17), (10,20), (11,12), (11,15), 
(11,20), (12,19), (12,24), \\(12,26), (13,15), (13,19), (14,17), (14,19),  (15,18),
 (16,18), (16,25), \\(17,18), (17,24),  
(18,23), (18,26), (19,23), (20,23), (21,23), (22,23) 
\end{eqnarray*}

\vspace{.25cm}
\noindent
Sets $S_{\Gamma_{0},3}$ and $S_{\Gamma_{0},4}$ are given below. Vertices of degree $5$ are written in bold.
\begin{eqnarray*}
S_{\Gamma_{0},3} &=&\{({\bf 0},{\bf 6},19), (1,{\bf 7},{\bf 14}), (3,{\bf 5},{\bf 18}), ({\bf 4},{\bf 10},26), ({\bf 23},{\bf 24},25)\}\\
S_{\Gamma_{0},4}&= &\{\}
\end{eqnarray*}

All of the vertices of degree $5$ in $S_{\Gamma_{0},3}$ have $4$ neighbours of degree $5$ and one of degree $4$.

\section*{Appendix I}\label{appendix28}
The $4$ extremal graph on $28$ vertices.  There are $68$ edges. Graph $\Gamma_{0}$  has $(deg_{3},deg_{4},deg_{5}) =(1,4,21,2)$, 
$\Gamma_{1}$ has $(deg_{4},deg_{5})=(4,24)$, $\Gamma_{2}$ has $(deg_{4},deg_{5},deg_{6}) = (6,20,2)$ and  $\Gamma_{3}$ has $(deg_{4},deg_{5},deg_{6}) = (7,18,3)$.

\vspace{.25cm}
\noindent
$\Gamma_{0}$ has incidence array:

\vspace{.25cm}
\noindent
\begin{tiny}
$
\begin{array}{cccccccccccccccccccccccccccc}
0&0&0&0&0&0&0&1&0&0&0&1&0&0&0&0&0&1&0&0&0&0&1&0&0&1&0&1\\
0&0&0&0&0&0&1&0&0&0&1&0&0&0&0&0&1&0&0&0&0&0&1&0&0&0&0&0\\
0&0&0&0&0&1&0&0&0&1&0&0&0&0&1&1&0&0&0&0&0&0&1&0&0&0&0&0\\
0&0&0&0&1&0&0&0&0&0&0&0&0&1&0&0&0&0&0&0&0&0&1&0&1&0&0&0\\
0&0&0&1&0&0&0&0&0&0&0&1&0&0&1&0&1&0&0&0&0&1&0&0&0&0&0&0\\
0&0&1&0&0&0&0&0&0&0&1&0&1&0&0&0&0&0&0&0&0&1&0&0&0&1&0&0\\
0&1&0&0&0&0&0&0&1&0&0&0&0&0&0&1&0&0&0&0&0&1&0&0&1&0&0&1\\
1&0&0&0&0&0&0&0&0&1&0&0&0&1&0&0&0&0&0&0&0&1&0&0&0&0&1&0\\
0&0&0&0&0&0&1&0&0&0&0&0&0&0&1&0&0&0&0&0&1&0&0&0&0&1&1&0\\
0&0&1&0&0&0&0&1&0&0&0&0&0&0&0&0&1&0&0&0&1&0&0&0&1&0&0&0\\
0&1&0&0&0&1&0&0&0&0&0&0&0&1&0&0&0&1&0&0&1&0&0&0&0&0&0&0\\
1&0&0&0&1&0&0&0&0&0&0&0&1&0&0&1&0&0&0&0&1&0&0&0&0&0&0&0\\
0&0&0&0&0&1&0&0&0&0&0&1&0&0&0&0&0&0&0&1&0&0&0&0&1&0&1&0\\
0&0&0&1&0&0&0&1&0&0&1&0&0&0&0&1&0&0&0&1&0&0&0&0&0&0&0&0\\
0&0&1&0&1&0&0&0&1&0&0&0&0&0&0&0&0&1&0&1&0&0&0&0&0&0&0&0\\
0&0&1&0&0&0&1&0&0&0&0&1&0&1&0&0&0&0&1&0&0&0&0&0&0&0&0&0\\
0&1&0&0&1&0&0&0&0&1&0&0&0&0&0&0&0&0&1&0&0&0&0&0&0&1&0&0\\
1&0&0&0&0&0&0&0&0&0&1&0&0&0&1&0&0&0&1&0&0&0&0&0&1&0&0&0\\
0&0&0&0&0&0&0&0&0&0&0&0&0&0&0&1&1&1&0&0&0&0&0&1&0&0&1&0\\
0&0&0&0&0&0&0&0&0&0&0&0&1&1&1&0&0&0&0&0&0&0&0&1&0&0&0&1\\
0&0&0&0&0&0&0&0&1&1&1&1&0&0&0&0&0&0&0&0&0&0&0&1&0&0&0&0\\
0&0&0&0&1&1&1&1&0&0&0&0&0&0&0&0&0&0&0&0&0&0&0&1&0&0&0&0\\
1&1&1&1&0&0&0&0&0&0&0&0&0&0&0&0&0&0&0&0&0&0&0&1&0&0&0&0\\
0&0&0&0&0&0&0&0&0&0&0&0&0&0&0&0&0&0&1&1&1&1&1&0&0&0&0&0\\
0&0&0&1&0&0&1&0&0&1&0&0&1&0&0&0&0&1&0&0&0&0&0&0&0&0&0&0\\
1&0&0&0&0&1&0&0&1&0&0&0&0&0&0&0&1&0&0&0&0&0&0&0&0&0&0&0\\
0&0&0&0&0&0&0&1&1&0&0&0&1&0&0&0&0&0&1&0&0&0&0&0&0&0&0&0\\
1&0&0&0&0&0&1&0&0&0&0&0&0&0&0&0&0&0&0&1&0&0&0&0&0&0&0&0\\
\end{array}
$
\end{tiny}

\vspace{.25cm}
\noindent
and degree sequence: $$(6,4,5,4,5,5,6,5,5,5,5,5,5,5,5,5,5,5,5,5,5,5,5,5,5,4,4,3)$$ Edge set is:
\begin{eqnarray*}
(0,7), (0,11), (0,17), (0,22), (0,25), (0,27), (1,6), (1,10), (1,16), 
(1,22),(2,5), \\
(2,9), (2,14), (2,15), (2,22), (3,4), (3,13), (3,22), (3,24), (4,11),
 (4,14), (4,16),\\ 
(4,21), (5,10), (5,12), (5,21), (5,25), (6,8), (6,15), 
(6,21), (6,24),(6,27),\\
(7,9), (7,13), (7,21), (7,26), (8,14), (8,20), (8,25), 
(8,26),
 (9,16),(9,20), \\
(9,24), (10,13), (10,17), (10,20), (11,12), (11,15), (11,20), 
(12,19),
 (12,24),\\ (12,26), (13,15), (13,19), (14,17), (14,19), (15,18), (16,18), (16,25), 
(17,18),\\
 (17,24), (18,23),(18,26), (19,23), (19,27), (20,23), (21,23), (22,23), 
\end{eqnarray*}

\vspace{.25cm}
\noindent
 All vertices of degree $5$ are the root of an embedded $S_{5,[4,4,4,4,4]}$ star.

\vspace{.25cm}
\noindent
Sets $S_{\Gamma_{0},3}$ and $S_{\Gamma_{0},4}$ are given below. Vertices of degree $5$ are written in bold. Note that both elements of $S_{\Gamma_{0},4}$ contain the vertex of degree $3$ (i.e. $27$), which is adjacent to the vertex of degree $6$ (i.e. $0$).
\begin{eqnarray*}
S_{\Gamma_{0},3} &=&\{(1,{\bf 7},{\bf 14}),({\bf 2},26,27),(3,{\bf 5},{\bf 18}),(3,{\bf 5},27),(3,{\bf 18},27),
\\&&(3,{\bf 20},27),(3,26,27),({\bf 4},{\bf 10},26),({\bf 4},{\bf 10},27),({\bf 4},26,27),\\&&
({\bf 5},{\bf 18},27),({\bf 10},26,27),({\bf 23},{\bf 24},25)\}\\
S_{\Gamma_{0},4}&= &\{(3,{\bf 5},{\bf 18},27) ({\bf 4},{\bf 10},26,27)\}
\end{eqnarray*}

\vspace{.25cm}
\noindent
$\Gamma_{1}$ has incidence array:

\vspace{.25cm}
\noindent
\begin{tiny}
$
\begin{array}{cccccccccccccccccccccccccccc}
0&0&0&0&0&0&0&1&0&0&0&1&0&0&0&0&0&1&0&0&0&0&1&0&0&1&0&0\\
0&0&0&0&0&0&1&0&0&0&1&0&0&0&0&0&1&0&0&0&0&0&1&0&0&0&1&0\\
0&0&0&0&0&1&0&0&0&1&0&0&0&0&1&1&0&0&0&0&0&0&1&0&0&0&0&0\\
0&0&0&0&1&0&0&0&0&0&0&0&0&1&0&0&0&0&0&0&0&0&1&0&1&0&0&1\\
0&0&0&1&0&0&0&0&0&0&0&1&0&0&1&0&1&0&0&0&0&1&0&0&0&0&0&0\\
0&0&1&0&0&0&0&0&0&0&1&0&1&0&0&0&0&0&0&0&0&1&0&0&0&1&0&0\\
0&1&0&0&0&0&0&0&1&0&0&0&0&0&0&1&0&0&0&0&0&1&0&0&1&0&0&0\\
1&0&0&0&0&0&0&0&0&1&0&0&0&1&0&0&0&0&0&0&0&1&0&0&0&0&1&0\\
0&0&0&0&0&0&1&0&0&0&0&0&0&0&1&0&0&0&0&0&1&0&0&0&0&1&0&1\\
0&0&1&0&0&0&0&1&0&0&0&0&0&0&0&0&1&0&0&0&1&0&0&0&1&0&0&0\\
0&1&0&0&0&1&0&0&0&0&0&0&0&1&0&0&0&1&0&0&1&0&0&0&0&0&0&0\\
1&0&0&0&1&0&0&0&0&0&0&0&1&0&0&1&0&0&0&0&1&0&0&0&0&0&0&0\\
0&0&0&0&0&1&0&0&0&0&0&1&0&0&0&0&0&0&0&1&0&0&0&0&1&0&1&0\\
0&0&0&1&0&0&0&1&0&0&1&0&0&0&0&1&0&0&0&1&0&0&0&0&0&0&0&0\\
0&0&1&0&1&0&0&0&1&0&0&0&0&0&0&0&0&1&0&1&0&0&0&0&0&0&0&0\\
0&0&1&0&0&0&1&0&0&0&0&1&0&1&0&0&0&0&1&0&0&0&0&0&0&0&0&0\\
0&1&0&0&1&0&0&0&0&1&0&0&0&0&0&0&0&0&1&0&0&0&0&0&0&1&0&0\\
1&0&0&0&0&0&0&0&0&0&1&0&0&0&1&0&0&0&1&0&0&0&0&0&1&0&0&0\\
0&0&0&0&0&0&0&0&0&0&0&0&0&0&0&1&1&1&0&0&0&0&0&1&0&0&0&1\\
0&0&0&0&0&0&0&0&0&0&0&0&1&1&1&0&0&0&0&0&0&0&0&1&0&0&0&0\\
0&0&0&0&0&0&0&0&1&1&1&1&0&0&0&0&0&0&0&0&0&0&0&1&0&0&0&0\\
0&0&0&0&1&1&1&1&0&0&0&0&0&0&0&0&0&0&0&0&0&0&0&1&0&0&0&0\\
1&1&1&1&0&0&0&0&0&0&0&0&0&0&0&0&0&0&0&0&0&0&0&1&0&0&0&0\\
0&0&0&0&0&0&0&0&0&0&0&0&0&0&0&0&0&0&1&1&1&1&1&0&0&0&0&0\\
0&0&0&1&0&0&1&0&0&1&0&0&1&0&0&0&0&1&0&0&0&0&0&0&0&0&0&0\\
1&0&0&0&0&1&0&0&1&0&0&0&0&0&0&0&1&0&0&0&0&0&0&0&0&0&0&0\\
0&1&0&0&0&0&0&1&0&0&0&0&1&0&0&0&0&0&0&0&0&0&0&0&0&0&0&1\\
0&0&0&1&0&0&0&0&1&0&0&0&0&0&0&0&0&0&1&0&0&0&0&0&0&0&1&0\\
\end{array}
$
\end{tiny}

\vspace{.25cm}
\noindent
and degree sequence: $$(5,5,5,5,5,5,5,5,5,5,5,5,5,5,5,5,5,5,5,4,5,5,5,5,5,4,4,4)$$
Edge set is: 
\begin{eqnarray*}
(0,7), (0,11), (0,17), (0,22), (0,25), (1,6), (1,10), (1,16), (1,22), 
(1,26),
(2,5), \\ (2,9), (2,14), (2,15), (2,22), (3,4), (3,13), (3,22), (3,24), 
(3,27), 
(4,11), (4,14),\\ (4,16), (4,21), (5,10), (5,12), (5,21), (5,25), (6,8), 
(6,15),
 (6,21), \\(6,24), (7,9), (7,13), (7,21), (7,26), (8,14), (8,20), (8,25), 
(8,27),
 (9,16), \\(9,20), (9,24), (10,13), (10,17), (10,20), (11,12), (11,15), (11,20),
(12,19),\\ (12,24), (12,26), (13,15), (13,19), (14,17), (14,19), (15,18), (16,18), (16,25), \\
(17,18), (17,24), (18,23), (18,27), (19,23), (20,23), (21,23), (22,23), (26,27),
\end{eqnarray*}

\vspace{.25cm}
\noindent
Sets $S_{\Gamma_{1},3}$ and $S_{\Gamma_{1},4}$ are given below. Vertices of degree $5$ are written in bold.
\begin{eqnarray*}
S_{\Gamma_{1},3} &=&\{({\bf 0},{\bf 6},19), ({\bf 0},19,27), ({\bf 8},{\bf 12},{\bf 22}), ({\bf 9},19,27),\\&& ({\bf 15},25,26), ({\bf 23},{\bf 24},25), ({\bf 23},25,26)\}\\
S_{\Gamma_{1},4}&= &\{\}
\end{eqnarray*}

\vspace{.25cm}
\noindent
$\Gamma_{2}$ has incidence array:

\vspace{.25cm}
\noindent
\begin{tiny}
$
\begin{array}{cccccccccccccccccccccccccccc}
0&0&0&0&0&0&0&1&0&0&0&1&0&0&0&0&0&1&0&0&0&0&1&0&0&1&0&0\\
0&0&0&0&0&0&1&0&0&0&1&0&0&0&0&0&1&0&0&0&0&0&1&0&0&0&1&0\\
0&0&0&0&0&1&0&0&0&1&0&0&0&0&1&1&0&0&0&0&0&0&1&0&0&0&0&0\\
0&0&0&0&1&0&0&0&0&0&0&0&0&1&0&0&0&0&0&0&0&0&1&0&1&0&0&1\\
0&0&0&1&0&0&0&0&0&0&0&1&0&0&1&0&1&0&0&0&0&1&0&0&0&0&0&0\\
0&0&1&0&0&0&0&0&0&0&1&0&1&0&0&0&0&0&0&0&0&1&0&0&0&1&0&1\\
0&1&0&0&0&0&0&0&1&0&0&0&0&0&0&1&0&0&0&0&0&1&0&0&1&0&0&0\\
1&0&0&0&0&0&0&0&0&1&0&0&0&1&0&0&0&0&0&0&0&1&0&0&0&0&1&0\\
0&0&0&0&0&0&1&0&0&0&0&0&0&0&1&0&0&0&0&0&1&0&0&0&0&1&0&0\\
0&0&1&0&0&0&0&1&0&0&0&0&0&0&0&0&1&0&0&0&1&0&0&0&1&0&0&0\\
0&1&0&0&0&1&0&0&0&0&0&0&0&1&0&0&0&1&0&0&1&0&0&0&0&0&0&0\\
1&0&0&0&1&0&0&0&0&0&0&0&1&0&0&1&0&0&0&0&1&0&0&0&0&0&0&0\\
0&0&0&0&0&1&0&0&0&0&0&1&0&0&0&0&0&0&0&1&0&0&0&0&1&0&0&0\\
0&0&0&1&0&0&0&1&0&0&1&0&0&0&0&1&0&0&0&1&0&0&0&0&0&0&0&0\\
0&0&1&0&1&0&0&0&1&0&0&0&0&0&0&0&0&1&0&1&0&0&0&0&0&0&1&0\\
0&0&1&0&0&0&1&0&0&0&0&1&0&1&0&0&0&0&1&0&0&0&0&0&0&0&0&0\\
0&1&0&0&1&0&0&0&0&1&0&0&0&0&0&0&0&0&1&0&0&0&0&0&0&1&0&0\\
1&0&0&0&0&0&0&0&0&0&1&0&0&0&1&0&0&0&1&0&0&0&0&0&1&0&0&0\\
0&0&0&0&0&0&0&0&0&0&0&0&0&0&0&1&1&1&0&0&0&0&0&1&0&0&0&1\\
0&0&0&0&0&0&0&0&0&0&0&0&1&1&1&0&0&0&0&0&0&0&0&1&0&0&0&0\\
0&0&0&0&0&0&0&0&1&1&1&1&0&0&0&0&0&0&0&0&0&0&0&1&0&0&0&0\\
0&0&0&0&1&1&1&1&0&0&0&0&0&0&0&0&0&0&0&0&0&0&0&1&0&0&0&0\\
1&1&1&1&0&0&0&0&0&0&0&0&0&0&0&0&0&0&0&0&0&0&0&1&0&0&0&0\\
0&0&0&0&0&0&0&0&0&0&0&0&0&0&0&0&0&0&1&1&1&1&1&0&0&0&0&0\\
0&0&0&1&0&0&1&0&0&1&0&0&1&0&0&0&0&1&0&0&0&0&0&0&0&0&0&0\\
1&0&0&0&0&1&0&0&1&0&0&0&0&0&0&0&1&0&0&0&0&0&0&0&0&0&0&0\\
0&1&0&0&0&0&0&1&0&0&0&0&0&0&1&0&0&0&0&0&0&0&0&0&0&0&0&1\\
0&0&0&1&0&1&0&0&0&0&0&0&0&0&0&0&0&0&1&0&0&0&0&0&0&0&1&0\\
\end{array}
$
\end{tiny}

\vspace{.25cm}
\noindent
and degree sequence: $$(5,5,5,5,5,6,5,5,4,5,5,5,4,5,6,5,5,5,5,4,5,5,5,5,5,4,4,4)$$ Edge set is:
\begin{eqnarray*}
(0,7), (0,11), (0,17), (0,22), (0,25), (1,6), (1,10), (1,16), (1,22), 
(1,26), (2,5),\\
 (2,9), (2,14), (2,15), (2,22), (3,4), (3,13), (3,22), (3,24), 
(3,27), (4,11), (4,14),\\
 (4,16), (4,21), (5,10), (5,12), (5,21), (5,25), (5,27), 
(6,8), (6,15), (6,21), \\
(6,24), (7,9), (7,13), (7,21), (7,26), (8,14), (8,20), 
(8,25), (9,16), (9,20),\\
 (9,24), (10,13), (10,17), (10,20), (11,12), (11,15), (11,20), 
(12,19), (12,24), \\
(13,15), (13,19), (14,17), (14,19), (14,26), (15,18), (16,18), (16,25), 
(17,18),\\
 (17,24), (18,23), (18,27), (19,23), (20,23), (21,23), (22,23), (26,27)
\end{eqnarray*}

\vspace{.25cm}
\noindent
Vertices of degree $6$ are $5$ and $14$. Sets $S_{\Gamma_{2},3}$ and $S_{\Gamma_{2},4}$ are given below. Vertices of degree $5$ are written in bold.
Elements of both sets  contain only vertices of degree $4$ and $5$. Note that for every  vertex $p$ of degree $6$ ($5$ and $14$) all $3$ elements of $S_{\Gamma_{2},4}$ contain precisely one vertex (of degree $4$) adjacent to $p$. (If $p=5$, vertices adjacent to $p$ in the elements of $S_{\Gamma_{2},4}$  are $27$, $12$ and $25$ respectively, and if $p=14$, vertices adjacent to $p$ in the  elements of $S_{\Gamma_{2},4}$  are $19$, $8$ and $26$ respectively.
\begin{eqnarray*}
S_{\Gamma_{2},3} &=&\{({\bf 0},{\bf 6},19),({\bf 0},{\bf 6},27),({\bf 0},19,27),({\bf 6},19,27),({\bf 7},8,12),\\&&
({\bf 7},8,{\bf 18}),({\bf 7},12,{\bf 18}),(8,12,{\bf 18}),(8,12,{\bf 22}),({\bf 9},19,27),\\&&
({\bf 15},25,26),({\bf 23},{\bf 24},{\bf 25}),({\bf 23},{\bf 24},26),({\bf 23},25,26),({\bf 24},25,26)\}\\
S_{\Gamma_{2},4}&= &\{({\bf 0},{\bf 6},19,27) ({\bf 7},8,12,{\bf 18}) ({\bf 23},{\bf 24},25,26)\}
\end{eqnarray*}

\vspace{.25cm}
\noindent
$\Gamma_{3}$ has incidence array:

\vspace{.25cm}
\noindent
\begin{tiny}
$
\begin{array}{cccccccccccccccccccccccccccc}
0&0&0&0&0&0&0&1&0&0&0&1&0&0&0&0&0&1&0&0&0&0&1&0&0&1&1&0\\
0&0&0&0&0&0&1&0&0&0&1&0&0&0&0&0&1&0&0&0&0&0&1&0&0&0&0&0\\
0&0&0&0&0&1&0&0&0&1&0&0&0&0&1&1&0&0&0&0&0&0&1&0&0&0&0&0\\
0&0&0&0&1&0&0&0&0&0&0&0&0&1&0&0&0&0&0&0&0&0&1&0&1&0&0&1\\
0&0&0&1&0&0&0&0&0&0&0&1&0&0&1&0&1&0&0&0&0&1&0&0&0&0&0&0\\
0&0&1&0&0&0&0&0&0&0&1&0&1&0&0&0&0&0&0&0&0&1&0&0&0&1&0&1\\
0&1&0&0&0&0&0&0&1&0&0&0&0&0&0&1&0&0&0&0&0&1&0&0&1&0&1&0\\
1&0&0&0&0&0&0&0&0&1&0&0&0&1&0&0&0&0&0&0&0&1&0&0&0&0&0&0\\
0&0&0&0&0&0&1&0&0&0&0&0&0&0&1&0&0&0&0&0&1&0&0&0&0&1&0&0\\
0&0&1&0&0&0&0&1&0&0&0&0&0&0&0&0&1&0&0&0&1&0&0&0&1&0&0&0\\
0&1&0&0&0&1&0&0&0&0&0&0&0&1&0&0&0&1&0&0&1&0&0&0&0&0&0&0\\
1&0&0&0&1&0&0&0&0&0&0&0&1&0&0&1&0&0&0&0&1&0&0&0&0&0&0&0\\
0&0&0&0&0&1&0&0&0&0&0&1&0&0&0&0&0&0&0&1&0&0&0&0&1&0&0&0\\
0&0&0&1&0&0&0&1&0&0&1&0&0&0&0&1&0&0&0&1&0&0&0&0&0&0&0&0\\
0&0&1&0&1&0&0&0&1&0&0&0&0&0&0&0&0&1&0&1&0&0&0&0&0&0&0&0\\
0&0&1&0&0&0&1&0&0&0&0&1&0&1&0&0&0&0&1&0&0&0&0&0&0&0&0&0\\
0&1&0&0&1&0&0&0&0&1&0&0&0&0&0&0&0&0&1&0&0&0&0&0&0&1&0&0\\
1&0&0&0&0&0&0&0&0&0&1&0&0&0&1&0&0&0&1&0&0&0&0&0&1&0&0&0\\
0&0&0&0&0&0&0&0&0&0&0&0&0&0&0&1&1&1&0&0&0&0&0&1&0&0&0&1\\
0&0&0&0&0&0&0&0&0&0&0&0&1&1&1&0&0&0&0&0&0&0&0&1&0&0&1&0\\
0&0&0&0&0&0&0&0&1&1&1&1&0&0&0&0&0&0&0&0&0&0&0&1&0&0&0&0\\
0&0&0&0&1&1&1&1&0&0&0&0&0&0&0&0&0&0&0&0&0&0&0&1&0&0&0&0\\
1&1&1&1&0&0&0&0&0&0&0&0&0&0&0&0&0&0&0&0&0&0&0&1&0&0&0&0\\
0&0&0&0&0&0&0&0&0&0&0&0&0&0&0&0&0&0&1&1&1&1&1&0&0&0&0&0\\
0&0&0&1&0&0&1&0&0&1&0&0&1&0&0&0&0&1&0&0&0&0&0&0&0&0&0&0\\
1&0&0&0&0&1&0&0&1&0&0&0&0&0&0&0&1&0&0&0&0&0&0&0&0&0&0&0\\
1&0&0&0&0&0&1&0&0&0&0&0&0&0&0&0&0&0&0&1&0&0&0&0&0&0&0&1\\
0&0&0&1&0&1&0&0&0&0&0&0&0&0&0&0&0&0&1&0&0&0&0&0&0&0&1&0\\
\end{array}
$
\end{tiny}

\vspace{.25cm}
\noindent
and degree sequence: $$(6,4,5,5,5,6,6,4,4,5,5,5,4,5,5,5,5,5,5,5,5,5,5,5,5,4,4,4)$$ Edge set is:
\begin{eqnarray*}
(0,7), (0,11), (0,17), (0,22), (0,25), (0,26), (1,6), (1,10), (1,16), 
(1,22), (2,5),\\
 (2,9), (2,14), (2,15), (2,22), (3,4), (3,13), (3,22), (3,24), 
(3,27), (4,11), (4,14),\\
 (4,16), (4,21), (5,10), (5,12), (5,21), (5,25), (5,27), 
(6,8), (6,15), (6,21),\\
 (6,24), (6,26), (7,9), (7,13), (7,21), (8,14), (8,20), 
(8,25), (9,16), (9,20),\\
 (9,24), (10,13), (10,17), (10,20), (11,12), (11,15), (11,20), 
(12,19), (12,24),\\
 (13,15), (13,19), (14,17), (14,19), (15,18), (16,18), (16,25), (17,18), 
(17,24),\\ (18,23), (18,27), (19,23), (19,26), (20,23), (21,23), (22,23), (26,27), 
\end{eqnarray*}

\vspace{.25cm}
\noindent
Vertices of degree $6$ are $0$, $5$ and $6$. Sets $S_{\Gamma_{3},3}$ and $S_{\Gamma_{3},4}$ are given below. Vertices of degree $5$ are written in bold.
Elements of both sets  contain only vertices of degree $4$ and $5$. Note that for every  vertex $p$ of degree $6$ ($0$, $5$ and $6$) both elements of $S_{\Gamma_{3},4}$ contain precisely one vertex (of degree $4$) adjacent to $p$. (If $p=0$, vertices adjacent to $p$ in the elements of $S_{\Gamma_{3},4}$  are $7$ in both cases,  if $p=5$ or $p=6$ vertices adjacent to $p$ in the  elements of $S_{\Gamma_{3},4}$  are $1$, and $8$  respectively.

\begin{eqnarray*}
S_{\Gamma_{3},3} &=&\{ (1,7,12),(1,7,{\bf 14}),(1,7,27),(1,{\bf 11},27),(1,{\bf 14},27),
({\bf 4},{\bf 10},26),(7,8,12),\\ &&(7,8,{\bf 18}),(7,8,27),(7,12,{\bf 18}),
(7,{\bf 14},27),(8,12,{\bf 18}),(8,12,{\bf 22}),({\bf 23},{\bf 24},25)  \}\\
S_{\Gamma_{3},4}&= &\{(1,7,{\bf 14},27) (7,8,12,{\bf 18})  \}
\end{eqnarray*}

\section*{Appendix J}\label{appendix29}
The $1$ extremal graph on $29$ vertices.  There are $72$ edges and  $(deg_{4},deg_{5},deg_{6}) =(5,20,4)$.  

\vspace{.25cm}
\noindent
$\Gamma_{0}$ has incidence array:

\vspace{.25cm}
\noindent
\begin{tiny}
$
\begin{array}{ccccccccccccccccccccccccccccc}
0&0&0&0&0&0&0&1&0&0&0&1&0&0&0&0&0&1&0&0&0&0&1&0&0&1&0&1&0\\
0&0&0&0&0&0&1&0&0&0&1&0&0&0&0&0&1&0&0&0&0&0&1&0&0&0&0&0&0\\
0&0&0&0&0&1&0&0&0&1&0&0&0&0&1&1&0&0&0&0&0&0&1&0&0&0&0&0&0\\
0&0&0&0&1&0&0&0&0&0&0&0&0&1&0&0&0&0&0&0&0&0&1&0&1&0&0&0&1\\
0&0&0&1&0&0&0&0&0&0&0&1&0&0&1&0&1&0&0&0&0&1&0&0&0&0&0&0&0\\
0&0&1&0&0&0&0&0&0&0&1&0&1&0&0&0&0&0&0&0&0&1&0&0&0&1&0&0&1\\
0&1&0&0&0&0&0&0&1&0&0&0&0&0&0&1&0&0&0&0&0&1&0&0&1&0&0&1&0\\
1&0&0&0&0&0&0&0&0&1&0&0&0&1&0&0&0&0&0&0&0&1&0&0&0&0&1&0&0\\
0&0&0&0&0&0&1&0&0&0&0&0&0&0&1&0&0&0&0&0&1&0&0&0&0&1&1&0&0\\
0&0&1&0&0&0&0&1&0&0&0&0&0&0&0&0&1&0&0&0&1&0&0&0&1&0&0&0&0\\
0&1&0&0&0&1&0&0&0&0&0&0&0&1&0&0&0&1&0&0&1&0&0&0&0&0&0&0&0\\
1&0&0&0&1&0&0&0&0&0&0&0&1&0&0&1&0&0&0&0&1&0&0&0&0&0&0&0&0\\
0&0&0&0&0&1&0&0&0&0&0&1&0&0&0&0&0&0&0&1&0&0&0&0&1&0&1&0&0\\
0&0&0&1&0&0&0&1&0&0&1&0&0&0&0&1&0&0&0&1&0&0&0&0&0&0&0&0&0\\
0&0&1&0&1&0&0&0&1&0&0&0&0&0&0&0&0&1&0&1&0&0&0&0&0&0&0&0&0\\
0&0&1&0&0&0&1&0&0&0&0&1&0&1&0&0&0&0&1&0&0&0&0&0&0&0&0&0&0\\
0&1&0&0&1&0&0&0&0&1&0&0&0&0&0&0&0&0&1&0&0&0&0&0&0&1&0&0&0\\
1&0&0&0&0&0&0&0&0&0&1&0&0&0&1&0&0&0&1&0&0&0&0&0&1&0&0&0&0\\
0&0&0&0&0&0&0&0&0&0&0&0&0&0&0&1&1&1&0&0&0&0&0&1&0&0&1&0&1\\
0&0&0&0&0&0&0&0&0&0&0&0&1&1&1&0&0&0&0&0&0&0&0&1&0&0&0&1&0\\
0&0&0&0&0&0&0&0&1&1&1&1&0&0&0&0&0&0&0&0&0&0&0&1&0&0&0&0&0\\
0&0&0&0&1&1&1&1&0&0&0&0&0&0&0&0&0&0&0&0&0&0&0&1&0&0&0&0&0\\
1&1&1&1&0&0&0&0&0&0&0&0&0&0&0&0&0&0&0&0&0&0&0&1&0&0&0&0&0\\
0&0&0&0&0&0&0&0&0&0&0&0&0&0&0&0&0&0&1&1&1&1&1&0&0&0&0&0&0\\
0&0&0&1&0&0&1&0&0&1&0&0&1&0&0&0&0&1&0&0&0&0&0&0&0&0&0&0&0\\
1&0&0&0&0&1&0&0&1&0&0&0&0&0&0&0&1&0&0&0&0&0&0&0&0&0&0&0&0\\
0&0&0&0&0&0&0&1&1&0&0&0&1&0&0&0&0&0&1&0&0&0&0&0&0&0&0&0&0\\
1&0&0&0&0&0&1&0&0&0&0&0&0&0&0&0&0&0&0&1&0&0&0&0&0&0&0&0&1\\
0&0&0&1&0&1&0&0&0&0&0&0&0&0&0&0&0&0&1&0&0&0&0&0&0&0&0&1&0\\
\end{array}
$
\end{tiny}

\vspace{.25cm}
\noindent
and degree sequence: $$(6,4,5,5,5,6,6,5,5,5,5,5,5,5,5,5,5,5,6,5,5,5,5,5,5,4,4,4,4)$$ Edge set is:
\begin{eqnarray*}
(0,7), (0,11), (0,17), (0,22), (0,25), (0,27), (1,6), (1,10), (1,16), (1,22),\\
 (2,5), (2,9), (2,14), (2,15), (2,22), (3,4), (3,13), (3,22), (3,24), 
(3,28),\\(4,11), 
(4,14), (4,16), (4,21), (5,10), (5,12), (5,21), (5,25), 
(5,28),(6,8),\\(6,15), (6,21), (6,24), (6,27), (7,9), (7,13), (7,21), 
(7,26), (8,14), \\
(8,20), (8,25), (8,26), (9,16), (9,20), (9,24), (10,13), 
(10,17), (10,20),\\ (11,12), 
(11,15), (11,20), (12,19), (12,24), (12,26), (13,15), (13,19), \\(14,17), (14,19),(15,18), 
(16,18), (16,25), (17,18), (17,24), (18,23), \\(18,26), (18,28), (19,23), (19,27), (20,23), 
(21,23), (22,23), (27,28)
\end{eqnarray*}

\vspace{.25cm}
\noindent
Sets $S_{\Gamma,3}$ and $S_{\Gamma,4}$ are given below. Vertices of degree $5$ are written in bold, all other elements have degree $4$..
\begin{eqnarray*}
S_{\Gamma_{0},3} &=&\{(1,{\bf 7},{\bf 14}),(1,{\bf 7},28),(1,{\bf 11},28),(1,{\bf 14},28),\\&&({\bf 2},26,27),
({\bf 4},{\bf 10},26),({\bf 4},{\bf 10},27),({\bf 4},26,27),\\&&({\bf 7},{\bf 14},28),(10,26,27)
({\bf 23},{\bf 24},25)\}\\
S_{\Gamma_{0},4}&= &\{(1,{\bf 7},{\bf 14},28) ({\bf 4},{\bf 10},26,27)\}
\end{eqnarray*}

\section*{Appendix K}\label{appendix30}
The $1$ extremal graph on $30$ vertices.  There are $76$ edges and  $(deg_{4},deg_{5},deg_{6}) =(4,20,6)$.  

\vspace{.25cm}
\noindent
$\Gamma_{0}$ has incidence array:

\vspace{.25cm}
\noindent
\begin{tiny}
$
\begin{array}{ccccc ccccc ccccc ccccc ccccc ccccc}
0&0&0&0&0&0&0&1&0&0&0&1&0&0&0&0&0&1&0&0&0&0&1&0&0&1&0&1&0&0\\
0&0&0&0&0&0&1&0&0&0&1&0&0&0&0&0&1&0&0&0&0&0&1&0&0&0&0&0&0&1\\
0&0&0&0&0&1&0&0&0&1&0&0&0&0&1&1&0&0&0&0&0&0&1&0&0&0&0&0&0&0\\
0&0&0&0&1&0&0&0&0&0&0&0&0&1&0&0&0&0&0&0&0&0&1&0&1&0&0&0&1&0\\
0&0&0&1&0&0&0&0&0&0&0&1&0&0&1&0&1&0&0&0&0&1&0&0&0&0&0&0&0&0\\
0&0&1&0&0&0&0&0&0&0&1&0&1&0&0&0&0&0&0&0&0&1&0&0&0&1&0&0&1&0\\
0&1&0&0&0&0&0&0&1&0&0&0&0&0&0&1&0&0&0&0&0&1&0&0&1&0&0&1&0&0\\
1&0&0&0&0&0&0&0&0&1&0&0&0&1&0&0&0&0&0&0&0&1&0&0&0&0&1&0&0&1\\
0&0&0&0&0&0&1&0&0&0&0&0&0&0&1&0&0&0&0&0&1&0&0&0&0&1&1&0&0&0\\
0&0&1&0&0&0&0&1&0&0&0&0&0&0&0&0&1&0&0&0&1&0&0&0&1&0&0&0&0&0\\
0&1&0&0&0&1&0&0&0&0&0&0&0&1&0&0&0&1&0&0&1&0&0&0&0&0&0&0&0&0\\
1&0&0&0&1&0&0&0&0&0&0&0&1&0&0&1&0&0&0&0&1&0&0&0&0&0&0&0&0&0\\
0&0&0&0&0&1&0&0&0&0&0&1&0&0&0&0&0&0&0&1&0&0&0&0&1&0&1&0&0&0\\
0&0&0&1&0&0&0&1&0&0&1&0&0&0&0&1&0&0&0&1&0&0&0&0&0&0&0&0&0&0\\
0&0&1&0&1&0&0&0&1&0&0&0&0&0&0&0&0&1&0&1&0&0&0&0&0&0&0&0&0&1\\
0&0&1&0&0&0&1&0&0&0&0&1&0&1&0&0&0&0&1&0&0&0&0&0&0&0&0&0&0&0\\
0&1&0&0&1&0&0&0&0&1&0&0&0&0&0&0&0&0&1&0&0&0&0&0&0&1&0&0&0&0\\
1&0&0&0&0&0&0&0&0&0&1&0&0&0&1&0&0&0&1&0&0&0&0&0&1&0&0&0&0&0\\
0&0&0&0&0&0&0&0&0&0&0&0&0&0&0&1&1&1&0&0&0&0&0&1&0&0&1&0&1&0\\
0&0&0&0&0&0&0&0&0&0&0&0&1&1&1&0&0&0&0&0&0&0&0&1&0&0&0&1&0&0\\
0&0&0&0&0&0&0&0&1&1&1&1&0&0&0&0&0&0&0&0&0&0&0&1&0&0&0&0&0&0\\
0&0&0&0&1&1&1&1&0&0&0&0&0&0&0&0&0&0&0&0&0&0&0&1&0&0&0&0&0&0\\
1&1&1&1&0&0&0&0&0&0&0&0&0&0&0&0&0&0&0&0&0&0&0&1&0&0&0&0&0&0\\
0&0&0&0&0&0&0&0&0&0&0&0&0&0&0&0&0&0&1&1&1&1&1&0&0&0&0&0&0&0\\
0&0&0&1&0&0&1&0&0&1&0&0&1&0&0&0&0&1&0&0&0&0&0&0&0&0&0&0&0&0\\
1&0&0&0&0&1&0&0&1&0&0&0&0&0&0&0&1&0&0&0&0&0&0&0&0&0&0&0&0&0\\
0&0&0&0&0&0&0&1&1&0&0&0&1&0&0&0&0&0&1&0&0&0&0&0&0&0&0&0&0&0\\
1&0&0&0&0&0&1&0&0&0&0&0&0&0&0&0&0&0&0&1&0&0&0&0&0&0&0&0&1&0\\
0&0&0&1&0&1&0&0&0&0&0&0&0&0&0&0&0&0&1&0&0&0&0&0&0&0&0&1&0&1\\
0&1&0&0&0&0&0&1&0&0&0&0&0&0&1&0&0&0&0&0&0&0&0&0&0&0&0&0&1&0\\

\end{array}
$
\end{tiny}

\vspace{.25cm}
\noindent
and degree sequence: $$(6,5,5,5,5,6,6,6,5,5,5,5,5,5,6,5,5,5,6,5,5,5,5,5,5,4,4,4,5,4)$$ Edge set is:
\begin{eqnarray*}
(0,7), (0,11), (0,17), (0,22), (0,25), (0,27), (1,6), (1,10),\\
 (1,16), (1,22), (1,29), (2,5), (2,9), (2,14), (2,15), (2,22),\\ 
(3,4), (3,13), (3,22), (3,24), (3,28), (4,11), (4,14), (4,16), \\
(4,21), (5,10), (5,12), (5,21), (5,25), (5,28), (6,8), (6,15), \\
(6,21), (6,24), (6,27), (7,9), (7,13), (7,21), (7,26), (7,29),\\
(8,14), (8,20), (8,25), (8,26), (9,16), (9,20), (9,24), (10,13), \\
(10,17), (10,20), (11,12), (11,15), (11,20), (12,19), (12,24), \\
(12,26), (13,15), (13,19), (14,17), (14,19), (14,29), (15,18), \\
(16,18), (16,25), (17,18), (17,24), (18,23), (18,26), (18,28), \\
(19,23), (19,27), (20,23), (21,23), (22,23), (27,28), (28,29)
\end{eqnarray*}

\vspace{.25cm}
\noindent
Sets $S_{\Gamma,3}$ and $S_{\Gamma,4}$ are given below. Vertices of degree $5$ are written in bold.
Elements of $S_{\Gamma,4}$ contain two vertices of degree $4$ and two of degree $5$.
\begin{eqnarray*}
S_{\Gamma_{0},3} &=&\{({\bf 2},26,27),({\bf 4},{\bf 10},26),({\bf 4},{\bf 10},27),({\bf 4},26,27),({\bf 10},26,27),\\&&
({\bf 15},25,29),({\bf 23},{\bf 24},25),({\bf 23},{\bf 24},29),({\bf 23},25,29),({\bf 24},25,29)\}\\
S_{\Gamma_{0},4}&= &\{({\bf 4},{\bf 10},26,27) ({\bf 23},{\bf 24},25,29)\}
\end{eqnarray*}

Let $s_{1}=({\bf 4},{\bf 10},26,27)$ and $s_{2}=({\bf 23},{\bf 24},25,29)$ (the two elements of $S_{\Gamma_{0},4}$ respectively). 
Consider the embedded $S_{6,[4,4,4,4,4,3]}$  stars with roots $6$ and $14$ respectively. The first star has children $27$ and $24$ which are from $s_{1}$ and $s_{2}$ respectively, with degrees $4$ and $5$, and the second star has children $4$ and $29$ which are from $s_{1}$ and $s_{2}$ respectively, with degrees $5$ and $4$ respectively. 

\section*{Appendix L}\label{appendix31}
The $2$ extremal graph on $31$ vertices.  There are $80$ edges. Graph $\Gamma_{0}$ has  $(deg_{4},deg_{5},deg_{6}) =(3,20,8)$ and $\Gamma_{1}$ has $(deg_{5},deg_{6})=(26 ,5)$.

\vspace{.25cm}
\noindent
$\Gamma_{0}$ has incidence array:

\vspace{.25cm}
\noindent
\begin{tiny}
$
\begin{array}{ccccc ccccc ccccc ccccc ccccc ccccc c}
0&0&0&0&0&0&0&1&0&0&0&1&0&0&0&0&0&1&0&0&0&0&1&0&0&1&0&1&0&0&0\\
0&0&0&0&0&0&1&0&0&0&1&0&0&0&0&0&1&0&0&0&0&0&1&0&0&0&0&0&0&1&0\\
0&0&0&0&0&1&0&0&0&1&0&0&0&0&1&1&0&0&0&0&0&0&1&0&0&0&0&0&0&0&0\\
0&0&0&0&1&0&0&0&0&0&0&0&0&1&0&0&0&0&0&0&0&0&1&0&1&0&0&0&1&0&0\\
0&0&0&1&0&0&0&0&0&0&0&1&0&0&1&0&1&0&0&0&0&1&0&0&0&0&0&0&0&0&1\\
0&0&1&0&0&0&0&0&0&0&1&0&1&0&0&0&0&0&0&0&0&1&0&0&0&1&0&0&1&0&0\\
0&1&0&0&0&0&0&0&1&0&0&0&0&0&0&1&0&0&0&0&0&1&0&0&1&0&0&1&0&0&0\\
1&0&0&0&0&0&0&0&0&1&0&0&0&1&0&0&0&0&0&0&0&1&0&0&0&0&1&0&0&1&0\\
0&0&0&0&0&0&1&0&0&0&0&0&0&0&1&0&0&0&0&0&1&0&0&0&0&1&1&0&0&0&0\\
0&0&1&0&0&0&0&1&0&0&0&0&0&0&0&0&1&0&0&0&1&0&0&0&1&0&0&0&0&0&0\\
0&1&0&0&0&1&0&0&0&0&0&0&0&1&0&0&0&1&0&0&1&0&0&0&0&0&0&0&0&0&1\\
1&0&0&0&1&0&0&0&0&0&0&0&1&0&0&1&0&0&0&0&1&0&0&0&0&0&0&0&0&0&0\\
0&0&0&0&0&1&0&0&0&0&0&1&0&0&0&0&0&0&0&1&0&0&0&0&1&0&1&0&0&0&0\\
0&0&0&1&0&0&0&1&0&0&1&0&0&0&0&1&0&0&0&1&0&0&0&0&0&0&0&0&0&0&0\\
0&0&1&0&1&0&0&0&1&0&0&0&0&0&0&0&0&1&0&1&0&0&0&0&0&0&0&0&0&1&0\\
0&0&1&0&0&0&1&0&0&0&0&1&0&1&0&0&0&0&1&0&0&0&0&0&0&0&0&0&0&0&0\\
0&1&0&0&1&0&0&0&0&1&0&0&0&0&0&0&0&0&1&0&0&0&0&0&0&1&0&0&0&0&0\\
1&0&0&0&0&0&0&0&0&0&1&0&0&0&1&0&0&0&1&0&0&0&0&0&1&0&0&0&0&0&0\\
0&0&0&0&0&0&0&0&0&0&0&0&0&0&0&1&1&1&0&0&0&0&0&1&0&0&1&0&1&0&0\\
0&0&0&0&0&0&0&0&0&0&0&0&1&1&1&0&0&0&0&0&0&0&0&1&0&0&0&1&0&0&0\\
0&0&0&0&0&0&0&0&1&1&1&1&0&0&0&0&0&0&0&0&0&0&0&1&0&0&0&0&0&0&0\\
0&0&0&0&1&1&1&1&0&0&0&0&0&0&0&0&0&0&0&0&0&0&0&1&0&0&0&0&0&0&0\\
1&1&1&1&0&0&0&0&0&0&0&0&0&0&0&0&0&0&0&0&0&0&0&1&0&0&0&0&0&0&0\\
0&0&0&0&0&0&0&0&0&0&0&0&0&0&0&0&0&0&1&1&1&1&1&0&0&0&0&0&0&0&0\\
0&0&0&1&0&0&1&0&0&1&0&0&1&0&0&0&0&1&0&0&0&0&0&0&0&0&0&0&0&0&0\\
1&0&0&0&0&1&0&0&1&0&0&0&0&0&0&0&1&0&0&0&0&0&0&0&0&0&0&0&0&0&0\\
0&0&0&0&0&0&0&1&1&0&0&0&1&0&0&0&0&0&1&0&0&0&0&0&0&0&0&0&0&0&1\\
1&0&0&0&0&0&1&0&0&0&0&0&0&0&0&0&0&0&0&1&0&0&0&0&0&0&0&0&1&0&1\\
0&0&0&1&0&1&0&0&0&0&0&0&0&0&0&0&0&0&1&0&0&0&0&0&0&0&0&1&0&1&0\\
0&1&0&0&0&0&0&1&0&0&0&0&0&0&1&0&0&0&0&0&0&0&0&0&0&0&0&0&1&0&0\\
0&0&0&0&1&0&0&0&0&0&1&0&0&0&0&0&0&0&0&0&0&0&0&0&0&0&1&1&0&0&0\\
\end{array}
$
\end{tiny}

\vspace{.25cm}
\noindent
and degree sequence: $$(6,5,5,5,6,6,6,6,5,5,6,5,5,5,6,5,5,5,6,5,5,5,5,5,5,4,5,5,5,4,4)$$ Edge set is:
\begin{eqnarray*}
(0,7), (0,11), (0,17), (0,22), (0,25), (0,27), (1,6), (1,10), (1,16),\\ 
(1,22), (1,29), (2,5), (2,9), (2,14), (2,15), (2,22), (3,4), (3,13),\\ 
(3,22), (3,24), (3,28), (4,11), (4,14), (4,16), (4,21), (4,30), (5,10),\\ 
(5,12), (5,21), (5,25), (5,28), (6,8), (6,15), (6,21), (6,24),\\ (6,27),
 (7,9), (7,13), (7,21), (7,26), (7,29), (8,14), (8,20),\\ (8,25), (8,26),
 (9,16), (9,20), (9,24), (10,13), (10,17), (10,20),\\ (10,30), (11,12),
 (11,15), (11,20), (12,19), (12,24), (12,26), (13,15),\\ 
(13,19), (14,17),
(14,19), (14,29), (15,18), (16,18), (16,25),\\ (17,18), (17,24), (18,23), 
(18,26), (18,28), (19,23), (19,27), \\(20,23), (21,23), (22,23), (26,30), (27,28), (27,30), 
(28,29)
\end{eqnarray*}

\vspace{.25cm}
\noindent

Set $S_{\Gamma_{0},5}$ is given below. Vertices of degree $5$ are written in bold, all other vertices have degree $4$.
\begin{eqnarray*}
S_{\Gamma_{0},5} &=&\{({\bf 23},{\bf 24},25,29,30)\}
\end{eqnarray*}

\vspace{.25cm}
\noindent
$\Gamma_{1}$ has incidence array:

\vspace{.25cm}
\noindent
\begin{tiny}
$
\begin{array}{ccccc ccccc ccccc ccccc ccccc ccccc c}
0&0&0&0&0&0&0&0&0&0&0&1&0&0&0&1&0&0&0&1&0&0&0&1&0&0&0&0&0&1&0\\
0&0&0&0&0&0&0&1&0&0&1&0&0&0&1&0&0&0&1&0&0&0&0&0&0&0&0&0&0&1&0\\
0&0&0&0&0&0&1&0&0&1&0&0&0&0&0&0&0&1&0&0&0&0&1&0&0&0&0&0&0&1&0\\
0&0&0&0&0&1&0&0&1&0&0&0&0&1&0&0&0&0&0&0&0&1&0&0&0&0&0&0&0&1&0\\
0&0&0&0&0&0&0&0&0&0&1&0&0&1&0&0&1&0&0&0&0&0&1&0&0&0&0&0&1&0&0\\
0&0&0&1&0&0&0&0&0&0&0&1&0&0&0&0&0&1&0&0&1&0&0&0&0&0&0&0&1&0&0\\
0&0&1&0&0&0&0&0&1&0&0&0&0&0&0&1&0&0&1&0&0&0&0&0&0&0&0&0&1&0&0\\
0&1&0&0&0&0&0&0&0&1&0&0&1&0&0&0&0&0&0&0&0&0&0&1&0&0&0&0&1&0&0\\
0&0&0&1&0&0&1&0&0&0&0&0&0&0&1&0&1&0&0&0&0&0&0&1&0&0&0&1&0&0&0\\
0&0&1&0&0&0&0&1&0&0&0&0&0&1&0&0&0&0&0&1&1&0&0&0&0&0&0&1&0&0&0\\
0&1&0&0&1&0&0&0&0&0&0&0&0&0&0&1&0&1&0&0&0&1&0&0&0&0&0&1&0&0&0\\
1&0&0&0&0&1&0&0&0&0&0&0&1&0&0&0&0&0&1&0&0&0&1&0&0&0&0&1&0&0&0\\
0&0&0&0&0&0&0&1&0&0&0&1&0&0&0&0&1&0&0&0&0&1&0&0&0&0&1&0&0&0&0\\
0&0&0&1&1&0&0&0&0&1&0&0&0&0&0&0&0&0&1&0&0&0&0&0&0&0&1&0&0&0&0\\
0&1&0&0&0&0&0&0&1&0&0&0&0&0&0&0&0&0&0&1&0&0&1&0&0&0&1&0&0&0&0\\
1&0&0&0&0&0&1&0&0&0&1&0&0&0&0&0&0&0&0&0&1&0&0&0&0&0&1&0&0&0&0\\
0&0&0&0&1&0&0&0&1&0&0&0&1&0&0&0&0&0&0&0&1&0&0&0&0&1&0&0&0&0&0\\
0&0&1&0&0&1&0&0&0&0&1&0&0&0&0&0&0&0&0&0&0&0&0&1&0&1&0&0&0&0&0\\
0&1&0&0&0&0&1&0&0&0&0&1&0&1&0&0&0&0&0&0&0&0&0&0&0&1&0&0&0&0&0\\
1&0&0&0&0&0&0&0&0&1&0&0&0&0&1&0&0&0&0&0&0&1&0&0&0&1&0&0&0&0&0\\
0&0&0&0&0&1&0&0&0&1&0&0&0&0&0&1&1&0&0&0&0&0&0&0&1&0&0&0&0&0&0\\
0&0&0&1&0&0&0&0&0&0&1&0&1&0&0&0&0&0&0&1&0&0&0&0&1&0&0&0&0&0&0\\
0&0&1&0&1&0&0&0&0&0&0&1&0&0&1&0&0&0&0&0&0&0&0&0&1&0&0&0&0&0&0\\
1&0&0&0&0&0&0&1&1&0&0&0&0&0&0&0&0&1&0&0&0&0&0&0&1&0&0&0&0&0&0\\
0&0&0&0&0&0&0&0&0&0&0&0&0&0&0&0&0&0&0&0&1&1&1&1&0&0&0&0&0&0&1\\
0&0&0&0&0&0&0&0&0&0&0&0&0&0&0&0&1&1&1&1&0&0&0&0&0&0&0&0&0&0&1\\
0&0&0&0&0&0&0&0&0&0&0&0&1&1&1&1&0&0&0&0&0&0&0&0&0&0&0&0&0&0&1\\
0&0&0&0&0&0&0&0&1&1&1&1&0&0&0&0&0&0&0&0&0&0&0&0&0&0&0&0&0&0&1\\
0&0&0&0&1&1&1&1&0&0&0&0&0&0&0&0&0&0&0&0&0&0&0&0&0&0&0&0&0&0&1\\
1&1&1&1&0&0&0&0&0&0&0&0&0&0&0&0&0&0&0&0&0&0&0&0&0&0&0&0&0&0&1\\
0&0&0&0&0&0&0&0&0&0&0&0&0&0&0&0&0&0&0&0&0&0&0&0&1&1&1&1&1&1&0\\
\end{array}
$
\end{tiny}

\vspace{.25cm}
\noindent
and degree sequence: $$(5,5,5,5,5,5,5,5,6,6,6,6,5,5,5,5,5,5,5,5,5,5,5,5,5,5,5,5,5,5,6)$$ Edge set is:
\begin{eqnarray*}
(0,11), (0,15), (0,19), (0,23), (0,29), (1,7), (1,10), (1,14), (1,18), \\
(1,29), (2,6), (2,9), (2,17), (2,22), (2,29), (3,5), (3,8), (3,13),\\
 (3,21), (3,29), (4,10), (4,13), (4,16), (4,22), (4,28), (5,11), (5,17),\\
 (5,20), (5,28), (6,8), (6,15), (6,18), (6,28), (7,9), (7,12),\\ 
(7,23),(7,28), (8,14), (8,16), (8,23), (8,27), (9,13), (9,19),\\
 (9,20), (9,27), (10,15), (10,17), (10,21), (10,27), (11,12), (11,18),\\
 (11,22), (11,27), (12,16), (12,21), (12,26), (13,18), (13,26), (14,19),\\ 
(14,22), (14,26), (15,20), (15,26), (16,20), (16,25), (17,23), \\
(17,25), (18,25), (19,21), (19,25), (20,24), (21,24), (22,24),\\
 (23,24), (24,30), (25,30), (26,30), (27,30), (28,30), (29,30)
\end{eqnarray*}

\vspace{.25cm}
\noindent
Set $S_{\Gamma_{1},5}$ is given below. Vertices of degree $5$ are written in bold.
\begin{eqnarray*}
S_{\Gamma_{1},5} &=&\{({\bf 3},{\bf 7},{\bf 15},{\bf 22},{\bf 25})\}
\end{eqnarray*}

\section*{Appendix M}\label{appendix32}
The $1$ extremal graph on $32$ vertices.  There are $85$ edges and $(deg_{5},deg_{6})=(22 ,10)$.

\vspace{.25cm}
\noindent
$\Gamma_{0}$ has incidence array:

\vspace{.25cm}
\noindent
\begin{tiny}
$
\begin{array}{ccccc ccccc ccccc ccccc ccccc ccccc cc}
0&0&0&0&0&0&0&1&0&0&0&1&0&0&0&0&0&1&0&0&0&0&1&0&0&1&0&1&0&0&0&0\\
0&0&0&0&0&0&1&0&0&0&1&0&0&0&0&0&1&0&0&0&0&0&1&0&0&0&0&0&0&1&0&0\\
0&0&0&0&0&1&0&0&0&1&0&0&0&0&1&1&0&0&0&0&0&0&1&0&0&0&0&0&0&0&0&0\\
0&0&0&0&1&0&0&0&0&0&0&0&0&1&0&0&0&0&0&0&0&0&1&0&1&0&0&0&1&0&0&0\\
0&0&0&1&0&0&0&0&0&0&0&1&0&0&1&0&1&0&0&0&0&1&0&0&0&0&0&0&0&0&1&0\\
0&0&1&0&0&0&0&0&0&0&1&0&1&0&0&0&0&0&0&0&0&1&0&0&0&1&0&0&1&0&0&0\\
0&1&0&0&0&0&0&0&1&0&0&0&0&0&0&1&0&0&0&0&0&1&0&0&1&0&0&1&0&0&0&0\\
1&0&0&0&0&0&0&0&0&1&0&0&0&1&0&0&0&0&0&0&0&1&0&0&0&0&1&0&0&1&0&0\\
0&0&0&0&0&0&1&0&0&0&0&0&0&0&1&0&0&0&0&0&1&0&0&0&0&1&1&0&0&0&0&0\\
0&0&1&0&0&0&0&1&0&0&0&0&0&0&0&0&1&0&0&0&1&0&0&0&1&0&0&0&0&0&0&0\\
0&1&0&0&0&1&0&0&0&0&0&0&0&1&0&0&0&1&0&0&1&0&0&0&0&0&0&0&0&0&1&0\\
1&0&0&0&1&0&0&0&0&0&0&0&1&0&0&1&0&0&0&0&1&0&0&0&0&0&0&0&0&0&0&0\\
0&0&0&0&0&1&0&0&0&0&0&1&0&0&0&0&0&0&0&1&0&0&0&0&1&0&1&0&0&0&0&0\\
0&0&0&1&0&0&0&1&0&0&1&0&0&0&0&1&0&0&0&1&0&0&0&0&0&0&0&0&0&0&0&0\\
0&0&1&0&1&0&0&0&1&0&0&0&0&0&0&0&0&1&0&1&0&0&0&0&0&0&0&0&0&1&0&0\\
0&0&1&0&0&0&1&0&0&0&0&1&0&1&0&0&0&0&1&0&0&0&0&0&0&0&0&0&0&0&0&0\\
0&1&0&0&1&0&0&0&0&1&0&0&0&0&0&0&0&0&1&0&0&0&0&0&0&1&0&0&0&0&0&0\\
1&0&0&0&0&0&0&0&0&0&1&0&0&0&1&0&0&0&1&0&0&0&0&0&1&0&0&0&0&0&0&0\\
0&0&0&0&0&0&0&0&0&0&0&0&0&0&0&1&1&1&0&0&0&0&0&1&0&0&1&0&1&0&0&0\\
0&0&0&0&0&0&0&0&0&0&0&0&1&1&1&0&0&0&0&0&0&0&0&1&0&0&0&1&0&0&0&0\\
0&0&0&0&0&0&0&0&1&1&1&1&0&0&0&0&0&0&0&0&0&0&0&1&0&0&0&0&0&0&0&0\\
0&0&0&0&1&1&1&1&0&0&0&0&0&0&0&0&0&0&0&0&0&0&0&1&0&0&0&0&0&0&0&0\\
1&1&1&1&0&0&0&0&0&0&0&0&0&0&0&0&0&0&0&0&0&0&0&1&0&0&0&0&0&0&0&0\\
0&0&0&0&0&0&0&0&0&0&0&0&0&0&0&0&0&0&1&1&1&1&1&0&0&0&0&0&0&0&0&1\\
0&0&0&1&0&0&1&0&0&1&0&0&1&0&0&0&0&1&0&0&0&0&0&0&0&0&0&0&0&0&0&1\\
1&0&0&0&0&1&0&0&1&0&0&0&0&0&0&0&1&0&0&0&0&0&0&0&0&0&0&0&0&0&0&1\\
0&0&0&0&0&0&0&1&1&0&0&0&1&0&0&0&0&0&1&0&0&0&0&0&0&0&0&0&0&0&1&0\\
1&0&0&0&0&0&1&0&0&0&0&0&0&0&0&0&0&0&0&1&0&0&0&0&0&0&0&0&1&0&1&0\\
0&0&0&1&0&1&0&0&0&0&0&0&0&0&0&0&0&0&1&0&0&0&0&0&0&0&0&1&0&1&0&0\\
0&1&0&0&0&0&0&1&0&0&0&0&0&0&1&0&0&0&0&0&0&0&0&0&0&0&0&0&1&0&0&1\\
0&0&0&0&1&0&0&0&0&0&1&0&0&0&0&0&0&0&0&0&0&0&0&0&0&0&1&1&0&0&0&1\\
0&0&0&0&0&0&0&0&0&0&0&0&0&0&0&0&0&0&0&0&0&0&0&1&1&1&0&0&0&1&1&0\\
\end{array}
$
\end{tiny}

\vspace{.25cm}
\noindent
and degree sequence: $$(6,5,5,5,6,6,6,6,5,5,6,5,5,5,6,5,5,5,6,5,5,5,5,6,6,5,5,5,5,5,5,5)$$ Edge set is:
\begin{eqnarray*}
(0,7), (0,11), (0,17), (0,22), (0,25), (0,27), (1,6), (1,10), (1,16), \\
(1,22), (1,29), (2,5), (2,9), (2,14), (2,15), (2,22), (3,4), (3,13), \\
(3,22), (3,24), (3,28), (4,11), (4,14), (4,16), (4,21), (4,30), (5,10),\\
 (5,12), (5,21), (5,25), (5,28), (6,8), (6,15), (6,21), (6,24), (6,27),\\
 (7,9), (7,13), (7,21), (7,26), (7,29), (8,14), (8,20),\\ (8,25),
(8,26), (9,16), (9,20), (9,24), (10,13), (10,17),\\ (10,20), (10,30),
(11,12), (11,15), (11,20), (12,19), (12,24), \\(12,26), (13,15), (13,19),
 (14,17), (14,19), (14,29), (15,18),\\ (16,18), (16,25), (17,18), (17,24),
 (18,23), (18,26),(18,28),\\ (19,23), (19,27), (20,23), (21,23), (22,23),
 (23,31), (24,31),\\ (25,31), (26,30), (27,28), (27,30), (28,29), (29,31), (30,31)
\end{eqnarray*}

\vspace{.25cm}
\noindent
Sets $S_{\Gamma_{0},j}$, for $j>2$ are empty.

\begin{eqnarray*}
S_{\Gamma_{0},2} &=&\{
(1,11), (1,12), (1,19), (1,26), (2,26), (2,27),\\ &&
(2,30), (2,31), (3,8), (3,20), (3,25), (3,26),\\&& (8,13), (8,22), (8,28), (9,19), (9,27), (9,28), \\&&(9,30), (11,28), 
(11,29), (11,31), (12,16), (12,22),\\&& (12,29), (13,16), (13,25), (13,31), (15,25), (15,29), \\&&
(15,30), (15,31), (16,19), (16,27), (17,21), (19,25), \\&&(20,27), (20,28), (20,29), (22,26),
(22,30)\}
\end{eqnarray*}

\end{document}